\documentclass[10pt]{book}
\usepackage{fancyhdr}
\fancypagestyle{plain}{%
\fancyhead{} 
  
}
\usepackage{amsmath,sissa,amssymb,amsfonts,amsthm,mathrsfs,graphics,graphicx,psfrag}
\def\sqr#1#2{\vbox{\hrule height .#2pt
\hbox{\vrule width .#2pt height #1pt \kern #1pt
\vrule width .#2pt}\hrule height .#2pt }}
\def\square{\sqr74}
\def\endproof{\hphantom{MM}\hfill\llap{$\square$}\goodbreak}
\def\n{\noindent}

\def\R{{\mathbb R}}
\def\Z{{\mathbb Z}}
\def\L{{\bf L}}
\def\T{{\bf T}}
\def\D{{\cal D}}
\def\I{{\cal I}}
\def\forall{\hbox{for all}~}
\def\per{{\rm per}}
\def\sgn{{\rm sign}}
\def\vp{\varphi}
\def\M{{\cal M}}
\def\F{{\cal F}}
\def\O{{\cal O}}
\def\C{{\cal C}}

\def\wto{\rightharpoonup}
\def\meas{\hbox{meas}}
\def\v{\vskip 1em}
\def\vs{\vskip 2em}
\def\vsk{\vskip 3em}
\def\ve{\varepsilon}

\def\sign{\,{\rm sign}}
\def \dis{\displaystyle}
\def \H{{H^1(\R)}}
\def \Hper{{H^1_{\rm per}}} 
\def\eps{\varepsilon}
\def\vphi{\varphi}
\def\rank{{\rm \, rank\,}}
\def\pr{{\rm \, pr}}
\def \N{\mathbb N}
\def\st{{\rm \ \ s.t.\ \ }}

\newtheorem{example}{Example}[chapter]
\newtheorem{dhef}{Definition}[chapter]
\newtheorem{lemma}{Lemma}[chapter]
\newtheorem{oss}{Remark}[chapter]
\newtheorem{prop}{Proposition}[chapter]
\newtheorem{theorem}{Theorem}[chapter]
\newtheorem{cor}{Corollary}[chapter]
\begin{document}
\begin{titlepage}


\begin{center}
{\Huge\bf Analysis of singular solutions for two nonlinear wave equations}\\[60ex]
\end{center}
\begin{minipage}[t]{0.45\textwidth}
\Large
{\sc Candidate}
\vspace{10pt}\\
Massimo Fonte
\end{minipage}
\hfill
\begin{minipage}[t]{0.45\textwidth}\raggedleft
\Large
{\sc Supervisor}
\vspace{10pt}\\
\noindent Prof. Alberto Bressan
\end{minipage}

\vfill 
\begin{center}
{\vrule width 5cm depth 0pt height .3mm}\\
\large Thesis submitted for the degree of {\it Doctor Philosophiae}
\\
Academic Year 2004/2005\\[6ex]
\end{center}
\nonumber
\end{titlepage}

\thispagestyle{empty}
%
%
\thispagestyle{empty}
\frontmatter
\pagenumbering{arabic}
\tableofcontents
\chapter{Introduction}

This thesis deals with two strongly nonlinear evolution \emph{Partial Differential Equation} (in the following named P.D.E.) arising from mathematical physics. 
The first one was introduced first by Fokas and Fuchssteiner \cite{FF} as a bi-Hamiltonian equation, and then was rediscovered by R. Camassa and D.D. Holm \cite{CH} as an higher order level of approximation of the unidirectional shallow water wave equation than the \emph{Korteweg-de Vries} equation \cite{KdV}. It can be written as
\begin{equation}
\label{kdiv0}
\left\{
\begin{array}{l}
\dis u(t,x): \R\times \R\mapsto \R 
\\
\\
\dis u_t+2\kappa u_x - u_{xxt}+3u  u_x = 2  u_x  u_{xx}+ u  u_{xxx},
\end{array}
\right.
\end{equation} 
here the unknown $u(t,x)$ represents the water's free surface over a flat bed and $\kappa$ is a constant related to the critical shallow-water wave speed (see also \cite{Jo} for an alternative derivation as an hyperelastic-rod wave equation). We refer to this equation as to the \emph{Camassa-Holm} equation, in honour to the first two authors which found a physical meaning stemming from the Euler equation.

The second PDE we want to study is a system of hyperbolic equations with quadratic source
\begin{equation}
\left\{
\begin{array}{l}
\dis (u_1,\dots,u_N)(t,x):\R^+\times \R^2\mapsto \R^N
\\
\\
\dis (u_i)_t + {\bf c}_i\cdot \nabla_x u_i= \sum_{j k} a_{i j k} u_j u_k, \quad \forall i=1\dots N
\end{array}
\right.
\label{i-dbzm-eq}
\end{equation}
which is a discretization of the velocities in the plane $\R^2$
for the Boltzmann equation 
$$
\left\{
\begin{array}{l}
\dis f(t,x,\xi):\R^+\times \R^3\times \R^3\mapsto \R
\\
\\
\dis \partial_t f(t,x,\xi)+ \xi\cdot \nabla_x f(t,x,\xi)=Q(f,f)(t,x,\xi)\,.
\end{array}
\right.
$$

The nonlinear nature of these equations leads to the possibility of \emph{blow up} in finite time either for the solution itself, or for the gradient of the solution. The typical situation of blow up in finite time is given by the following \emph{Ordinary Differential Equation} (O.D.E.)
$$
\frac d{dt} v= -v^2,\qquad v(0)=v_0
$$
where $v$ should be either the solution $u$ or the gradient $u_x$, possibly computed along the characteristic curves of the equation. It is well known that the solution of this equation has the behaviour $\approx \frac 1{t-T}$, where $T$ depends on the initial data and it is the time of blow-up whenever $v_0<0$. 

This situation can occur for the solutions of the Camassa-Holm equation in the limit case $\kappa=0$. With this condition, the equation (\ref{kdiv0}) may be rewritten in nonlocal form as
\begin{equation}
\label{prblCH}
 u_t +u\, u_x=- \frac 12 \left[e^{-|x|}*\left(u^2+\frac{u_x^2}2\right)\right]_x.
\end{equation}
Since the $H^1$norm is conserved for regular solutions, the $L^\infty$ norm of $u$ is bounded, namely $\|u\|_{L^\infty}\leq \|u\|_{H^1}\leq \sqrt{E}$.
Arguing as in the \emph{Steepening Lemma} (see \cite{CHH}), in \cite[theorem 4.1]{CE1} the authors prove that smooth solution to (\ref{prblCH}) may not be globally defined. Let start from an odd initial data $\bar u\in H^3(\R)$ which has an inflection point in $0$, $\bar u(0)<0$ and consider the evolution of the slope at the inflection point $s(t)\doteq u_x(t,0)<0$. The computation of the function $s$ gives the differential inequality 
$$
\frac{ds}{dt}\leq -\frac 12 s^2 \qquad s(0)=\bar u_x(0).
$$
If $\bar u_x(0)$ is sufficiently small, we get
$$
\frac 1{s(t)}\geq  \frac 1{s(0)}+\frac t2
$$
and then the slope becomes vertical at finite time. 
However, the singularities thus can occur only in form of wave breaking (see also \cite{McK,C1}), in fact even if its slope can become unbounded at finite time, the solution remains bounded, because of the inequality $\|u(t)\|_{L^\infty}\leq \sqrt E$.

\vskip 5pt
The aim of Part I of this thesis (see also \cite{BF2,F}) is to construct a continuous semigroup of global solutions in two main cases:
\begin{enumerate}
\item on the space $H^1_{\rm per}$ of spatially periodic functions, locally in $\H$;
\item on a domain of $H^1$ functions with a certain exponential decay at $x\to \pm \infty$. 
\end{enumerate}
Result of existence of solutions can be found in \cite{XZ1, XZ2}, and \cite{CHK1, CHK2} where the authors added a small diffusion term to the right hand side of \ref{prblCH} and obtained solution of the original equation as a vanishing viscosity limit. On the other hand, in \cite{BC2} was developed an alternative technique, which relies on a new set of dependent and independent variables with the specific purpose to resolve all the singularities. With this change of variable the solution can be obtained as the unique fixed point of a contractive transformation. 
In Chapter \ref{multipeak} we present yet another approach based on the Hamiltonian structure of the Camassa-Holm equation. We shall construct the semigroup of global solution starting from explicit solutions of the Camassa-Holm equation with initial condition in form of \emph{multipeakon} function
$$
u_0(x)=\sum_{j=1}^{N} p_j e^{-|x-q_j|}\,.
$$
The motivation of this choice is given by the form of traveling wave solution (see \cite{CH,CE1, CE2, CM1}). Looking for solution of the equation (\ref{prblCH}) in the traveling wave form $u(t,x)=U(x-ct)$, with a function $U$ that vanishes at infinity, one obtains the function $U= c e^{-|x-ct|}$, which is a peaked  \emph{soliton} (from this fact derives the shortened term \emph{peakon}). The multipeakon functions are stable, in fact not only a single peakon subject to (\ref{prblCH}) evolves with this form, but also the evolution of a superposition of traveling wave (e.g. initial data like $u_0$) remains of the same shape
$$
u(t,x)=\sum_{j=1}^{N} p_j(t) e^{-|x-q_j(t)|}\,.
$$
The reader can see also \cite{BSS} for a recursive reconstruction of the multipeakon solutions, and \cite{CS}  which prove that multipeakon solutions are orbitally stable, i.e. stable under a general nature of perturbations.

In \cite{HR} the authors prove the existence of a global multipeakon solution when the strengths $p_i$ are positive for all $i=1\dots N$. In this case the crucial fact is that no interaction between the peakons occurs, and then the gradient remains bounded, which yields existence and uniqueness of the coefficients $p_1(t),\dots, p_N(t)$ and $q_1(t),\dots, q_N(t)$ that are solutions of the Hamiltonian system
\begin{equation}
\left\{
\begin{array}{l}
\label{HamSys}
\dis \dot q_i = \sum_{j=1}^{N} p_j e^{-|q_i-q_j|},
\\
\dis \dot p_i = p_i \sum_{j=1}^{N} p_j \sign (q_i -q_j) e^{-|q_i-q_j|},
\end{array}
\right.
\end{equation}
with Hamiltonian $H=\sum_{i,j}p_i p_j e^{-|q_i-q_j|}$.
 
\begin{figure}
\psfrag{0}{$0$}
\psfrag{p}{$p$}
\psfrag{q}{$q$}
\psfrag{-p}{$-p$}
\psfrag{-q}{$-q$}
\psfrag{x}{$x$}
\psfrag{t<tau}{$t<\tau$}
\psfrag{t=tau}{$t=\tau$}
\centerline{\includegraphics[width=12cm]{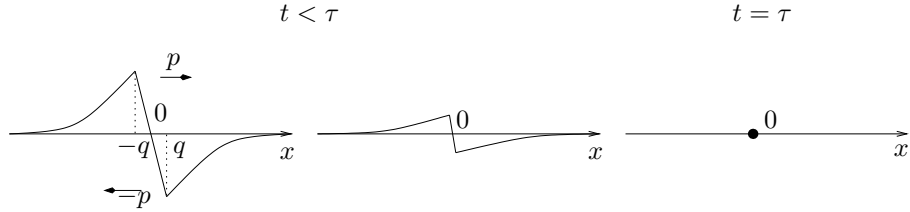}}
\caption{Peakon-antipeakon interaction.\label{figintro}}
\end{figure}

However, a general initial data contains both positive and negative peakons, as in the example of the peakon-antipeakon interaction: one positive peakon with strength $p$, centered in $-q$, moves forward and one negative anti-peakon in $q$, with strength $-p$ moves backward. The evolution of the system produces the overlapping of the two peakons at finite time $t=\tau$, so that $q\to 0$ (see Figure \ref{figintro}).

There are infinitely many ways to extend the solution after the time $\tau$ of the interaction, for example the vanishing viscosity approach in \cite{CHK1,CHK2} singles out the \emph{dissipative} solutions. As far as the example of peakon-antipeakon interaction is concerned, the vanishing viscosity approach selects the solution that, after the time $\tau$, is $\equiv 0$: all the energy $E$ is lost. In section \ref{3-2} we shall construct a \emph{conservative} solution, i.e. a solution for which the quantity $E$ is constant for a.e. time $t$. At the interaction time the energy $E\doteq \|u\|_{H^1}^2$ is described by a Dirac measure entirely concentrated at the single interaction point. After the interaction, a positive and a negative peakon emerge, whose strengths are uniquely determined by imposing the conservation of the total energy.

In Chapter \ref{chap3} we shall discuss the issue of the uniqueness and the stability.  

\vskip 5pt
\noindent\emph{Stability.}\quad
The multipeakon solutions form a continuous semigroup whose domain is dense either in $\H$ or in $\Hper$. The main novel feature in our approach is the construction of a metric $J(\cdot,\cdot)$ on the space $H^1$ (or $H^1_{\rm per}$) determined by an optimal transportation problem. While the semigroup generated by (\ref{prblCH}) is not eve continuous w.r.t. the $H^1$ distance, we show that it is Lipschitz continuous w.r.t. our new distance functional $J$. The reader can see \cite{V} for earlier applications of distances defined in term of optimal transportation problems, and \cite{BC1}, in which the authors recover a semigroup of dissipative solution for the Hunter-Saxton equation \cite{HS} (see also \cite{BZZ} for a fixed-point approach for both conservative and dissipative solutions).

The main well-posedness result is provided by a \emph{Gronwall-type lemma}, stemming from the inequality 
$$
\frac d{dt} J(u,v)\leq C(t) \cdot J(u,v)
$$
whenever $u$ and $v$ are two multipeakon solutions (see Section \ref{3-5}).

\vskip 5pt
\noindent\emph{Uniqueness.}\quad
Example \ref{ch1-ex} in Section \ref{3-7} shows that a solution of \ref{prblCH} need not be unique. Roughly speaking, every shifted antipeakon-peakon couple is also a conservative solution. A conservative solution can be characterized by an additional linear transport equation, accounting for the conservation of the total energy. It can be done by the following heuristic idea. 

We can think that the absolutely continuous measure $\mu_t$, which satisfies $d\mu_t=(u^2+u_x^2) d \mathcal{L}$ tends, as $t\to \tau$, to a Dirac measure with support in $0$. We introduce thus a further equation for the measure $\mu_t$, whose absolutely continuous part is $u^2+u_x^2$, in the following way. Since whenever $u^2+u_x^2$ is regular it satisfies the equation
$$
(u^2+u_x^2)_t+ \left[u(u^2+u_x^2)\right]_x=\left[u^3- 2 u e^{-|x|}*\left(u^2+\frac{u_x^2}2\right)\right]_x\doteq f(u)
$$
it suggests that $\mu_t$ provides a measure-valued solution of 
\begin{equation}
\label{measeqa}
\partial_t \mu+(u \mu)_x= f(u)\,.
\end{equation}
Our result (Theorem \ref{theo3}) shows that every solution $(u,\mu_t)$ of (\ref{prblCH})-(\ref{measeqa}), such that the absolutely continuous part of $\mu_t$ has density $u^2+u_x^2$, must coincide with the one provided by multipeakon approach.

\vskip 10pt
Part II of the thesis is devoted to the analysis of blow-up for the discrete Boltzmann equation (\ref{i-dbzm-eq}). Such a equation is obtained by considering a rarefied gas for which is supposed that the particles can move only along a finite number of direction characterized by the vectors ${\bf c}_1,\dots, {\bf c}_N$. The unknowns $u_i$ represent the densities of particles which travel at speed ${\bf c}_i$. By a collision, a pair of incoming particles with speeds ${\bf c}_i,\,{\bf c}_j$ is replaced by a new pair of particles say ${\bf c}_k,\,{\bf c}_\ell$. The rate at which such collision occur is given by $a_{ijk} u_j u_k$ The concentration $u_i$ is thus increasing (or at least is constant) when interact particles of speed different to ${\bf c}_i$, decreasing when an $i-$particle collides with someone other. Then, the coefficients $a_{ijk}$ are non negative when $j,k \neq i$ and negative when either $j=i$ or $k=i$. 

If the initial data is suitably small, the solution remains uniformly bounded for all times \cite{B2}. For large initial data, on the other hand, the global existence and stability of solutions are known only in the one-dimensional case \cite{B1, HT, T}. 
Since the right hand side has quadratic growth, it might happen that the solution blows up in finite time.
Examples where the $\L^\infty$ norm of the solution becomes arbitrarily large as $t\to\infty$ are easy to construct \cite{I}.  

In Chapter \ref{chapblowup} (see \cite{BF1}) we focus our analysis on the \emph{two-dimensional Broadwell model} (see, for example, \cite{B3, U, CIP} for a description of the model) and examine the possibility that blow-up actually occurs in finite time. 
\begin{figure}
\psfrag{c1}{${\bf c}_1$}
\psfrag{c2}{${\bf c}_2$}
\psfrag{c3}{${\bf c}_3$}
\psfrag{c4}{${\bf c}_4$}
\centerline{
\includegraphics[width=10.5cm]{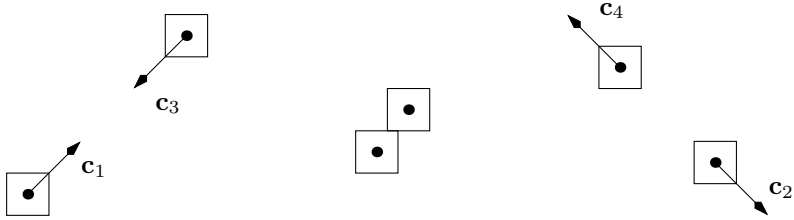}}
\caption{Two particle interaction \label{2boltz}}
\end{figure}
In this model the permitted direction are 
$$
{\bf c}_1=(1,1),~{\bf c}_2=(1,-1),~{\bf c}_3=(-1,-1),~{\bf c}_4=(-1,1)
$$ 
and the particles have a diamond-shape (see Figure \ref{2boltz}).

As we will show with the theory developed in Chapter \ref{symmgenerator}, since the equations (\ref{i-dbzm-eq}) admit a natural symmetry group (see Section \ref{symm}, and \cite{O} for a more general theory), one can perform an asymptotic rescaling of variables and ask whether there is a blow-up solution which, in the rescaled variables, converges to a steady state.  This technique has been widely used to study blow-up singularities of reaction-diffusion equations with superlinear forcing terms \cite{GV, GK}.  See also \cite{J} for an example of self-similar blow-up for hyperbolic conservation laws.
Our main results show is an \emph{a-priori} bound on the blow up rate in the $L^\infty$ norm. Namely, if blow-up occurs at time $T$, then one has
$$
\big\|u(t)\big\|_{\L^\infty}~>~ \frac 15\, \frac{\ln \big|\ln(T-t)\big|}{T-t}\,.
$$
This means that the blow-up rate must be different from the natural growth rate $\big\|u(t) \big\|_{\L^\infty} = \O(1)\cdot(T-t)^{-1}$ which would be obtained in case of a quadratic equation $\dot u=C\,u^2$.

In the final section of Chapter \ref{chapblowup} we discuss a possible scenario for blow-up.  The analysis highlights how carefully chosen should be the initial data, if blow-up is ever to happen. This suggests that finite time blow-up is a highly non-generic
phenomenon, something one would not expect to encounter in  numerical simulations.

\mainmatter
\part{The Camassa-Holm equation}
\chapter{The Camassa-Holm equation}
The Camassa-Holm equation 
$$
u_t+2\kappa u_x - u_{xxt} + 3 u u_x = 2 u_x u_{xx} + u u_{xxx}
$$
arises from a higher order level of approximation of the asymptotic expansion of the Euler's equations for a shallow water wave theory. Here we do not enter into deep details of the interpretation of such a equation, for the physical motivations we refer to \cite{CH}, \cite{CM1}, \cite{CM2}, \cite{Jo}.  

In the following we focus our attention and we refer to \emph{Camassa-Holm equation} the previous equation with $\kappa = 0$. 
\section{Non-local formulation\label{3-1}}
\v
The Camassa-Holm equation can be written as a scalar conservation law
with an additional integro-differential term:
\begin{equation}
u_t+(u^2/2)_x+P_x=0\,,
\label{ch1-equation}
\end{equation}
where $P$ is defined as a convolution:
\begin{equation}
P(x)\doteq{\frac 12} e^{-|x|} * \left(u^2+{\frac{u_x^2} 2}
\right)\,.
\label{ch1-nonlocalP}
\end{equation}
Earlier results on the
existence and uniqueness of solutions can be found in \cite{XZ1}, \cite{XZ2}.
One can regard (\ref{ch1-equation}) as an evolution equation on a space of
absolutely continuous functions with derivatives $u_x\in L^2$.
In the smooth case,
differentiating (\ref{ch1-equation}) w.r.t.~$x$ one obtains
\begin{equation}
u_{xt}+uu_{xx}+u_x^2-\left( u^2+\frac {u_x^2}2\right)+
P=0\,.
\label{ch1-nonlocalequation}
\end{equation}
Multiplying (\ref{ch1-equation}) by $u$ and (\ref{ch1-nonlocalequation}) by $u_x$ we
obtain the two conservation laws with source term
\begin{equation}
\left(\frac{u^2} 2\right)_t+\left(\frac{u^3}3+u\,P\right)_x=u_x \,P\,,
\end{equation}
\begin{equation}
\left(\frac{u^2_x} 2\right)_t+\left(\frac{uu_x^2}2-\frac{u^3}3\right)_x
=- u_x \,P\,. 
\label{ch1-uxquadro}
\end{equation}
As a consequence, for regular solutions the total energy
$$
E(t)\doteq \int \big[ u^2(t,x)+u_x^2(t,x)\big]\,dx
$$
remains constant in time.

As in the case of scalar conservation laws, by
the strong nonlinearity of the equations, solutions
with smooth initial data can lose regularity
in finite time.
For the Camassa-Holm equation (\ref{ch1-equation}), however,
the uniform bound on $\|u_x\|_{L^2}$
guarantees that only the $L^\infty$ norm of the
gradient can blow up, while the solution $u$ itself remains
H\"older continuous at all times.

In order to construct global in time solutions, two
main approaches have recently been introduced.
On one hand, one can add a small diffusion term in the
right hand side of (\ref{ch1-equation}), and recover solutions of the
original equations as a vanishing viscosity limit \cite{CHK1, CHK2}.
An alternative technique, developed in \cite{BC2}, relies on
a new set of independent and dependent variables, specifically
designed with the aim of
``resolving'' all singularities.
In terms of these new variables, the solution
to the Cauchy problem becomes regular for all times, and can
be obtained as the unique fixed point of a contractive transformation.

In the present chapter, we implement yet another approach
to the Camassa-Holm equation.
As a starting point we consider all multi-peakon solutions, of the form
\begin{equation}
u(t,x)=\sum_{i=1}^N p_i(t) e^{-|x-q_i(t)|}\,.
\label{ch1-mpeak}
\end{equation}
These are obtained by solving the system of O.D.E's
\begin{equation}
\left
\{
\begin{array}{rcl}
\displaystyle
 \dot q_i&=&\displaystyle
\sum_j p_j\,e^{-|q_i-q_j|}\,,
 \\
\displaystyle
 \dot p_i&=&
 \displaystyle
\sum_{j\not= i} p_i p_j\,\sgn(q_i-q_j)\,e^{-|q_i-q_j|}\,.
\end{array}
\right.
\label{ch1-Hsys}
\end{equation}
It is well known that this can
be written in hamiltonian form:
$$\left\{
\begin{array}{rcl}
\dot q_i&=&
\displaystyle
\frac\partial{\partial p_i} H(p,q)\,,
\\
\dot p_i&=&
\displaystyle
-\frac\partial{\partial q_i} H(p,q)\,,
\end{array}
\right.\qquad\qquad H(p,q)\doteq \frac 12 \sum_{i,j} p_ip_j
e^{-|q_i-q_j|}\,.
$$

If all the coefficients $p_i$ are initially positive, then
they remain positive and bounded for all times.
The solution $u=u(t,x)$ is thus uniformly Lipschitz continuous.
We stress, however, that
here we are not making any assumption about the signs of the
$p_i$.   In a typical situation,
two peakons can
cross at a finite time $\tau$. As $t\to\tau-$
their strengths $p_i,p_j$ and positions $q_i,q_j$
will satisfy
\begin{equation}
p_i(t)\to +\infty\,,\qquad p_j(t)\to -\infty\,,
\qquad p_i(t)+p_j(t)\to \bar p\,,
\label{ch1-condp}
\end{equation}
\begin{equation}
q_i(t)\to \bar q\,,
\quad q_j(t)\to \bar q\,,\qquad q_i(t)<q_j(t)
~~\mbox{for $t<\tau$},
\label{ch1-condcoeff}
\end{equation}
for some $\bar p,\bar q\in\R$.
Moreover, $\big\|u_x(t)\big\|_{L^\infty}\to\infty$.
In this case, we will show that
there exists
a unique way to extend the multi-peakon solution beyond
the interaction time, so that the total energy is conserved.

Having constructed a set of ``multi-peakon solutions'',
our main goal is to show that these solutions form a
continuous semigroup, whose domain is dense in the space
$H^1(\R)$.  Taking the unique continuous extension,
we thus obtain a continuous semigroup of solutions of
(\ref{ch1-equation}), defined on the entire space $H^1$.

One easily checks that the flow map
$\Phi_t: u(0)\mapsto u(t)$ cannot be continuous as
a map from $H^1$ into itself, or from $L^2$ into itself.
Distances defined in terms of convex norms perform well
in connection with linear problems, but occasionally fail
when nonlinear features become dominant.
In the present setting, we construct a new distance
$J(u,v)$ between functions $u,v\in H^1$, defined
by a problem of optimal transportation. Roughly speaking,
$J(u,v)$ will be the minimum cost in order to transport the
mass distribution with density $1+u_x^2$ located
on the graph of $u$ onto the
mass distribution with density $1+v_x^2$ located
on the graph of $v$.  See Section \ref{3-3} for details.
With this definition of distance, our main result shows that
$$
\left|\frac d{dt} J\big(u(t),\,v(t)\big)\right|~\leq~ C \cdot
J\big(u(t),\,v(t)\big)
$$
for some constant $C$ and any couple of multi-peakon solutions $u,v$.
Moreover, $J(u_n,u)\to 0$ implies the uniform convergence
$\|u_n- u\|_{L^\infty}\to 0$.
The distance functional $J$ thus provides the ideal tool to measure
continuous dependence on the initial data for solutions
to the Camassa-Holm equation.
Earlier applications
of distances defined in terms of optimal transportation
problems can be found in the monograph \cite{V}.
The issue of uniqueness of solutions must here be discussed in
greater detail.  For a multi-peakon solution, as long as all
coefficients $p_i$ remain bounded, the solution to
the system of ODE's  (\ref{ch1-Hsys}) is clearly unique.   For each time $t$, call
$\mu_t$ the measure having density
$u^2(t)+u_x^2(t)$ w.r.t.~Lebesgue measure.
Consider a time $\tau$
where a positive and a negative peakon collide, according
to (\ref{ch1-condp})-(\ref{ch1-condcoeff}).
As $t\to \tau- $, we have the weak convergence $\mu_t\wto
\mu_\tau$ for some positive measure $\mu_\tau$ which typically
contains a Dirac mass at the point $\bar q$.
By energy conservation, we thus have
$$\int\big[ u^2(\tau,x)+u_x^2(\tau, x)\big]\,dx +\mu_\tau\big(\{\bar q\}
\big)=\lim_{t\to\tau-}\int\big[ u^2(\tau,x)+u_x^2(\tau, x)\big]\,dx
= E(\tau-)\,.$$
There are now two natural ways to prolong the multi-peakon solution
beyond time $\tau$:   a conservative solution, such that
$$E(t)=\int\big[ u^2(t,x)+u_x^2(t, x)\big]\,dx = E(\tau -)\qquad\qquad
t>\tau\,,$$
or a dissipative solution, where all the energy concentrated
at the point $\bar q$ is lost.
In this case
$$E(t)=\int\big[ u^2(t,x)+u_x^2(t, x)\big]\,dx = E(\tau -)-
\mu_\tau\big(\{\bar q\}\big)\qquad\qquad
t>\tau\,.$$
For $t>\tau$,
the dissipative solution is obtained by simply replacing the two peakons
$p_i,p_j$ with one single peakon of strength $\bar p$, located at
$x=\bar q$.  On the other hand, as we will show in Section \ref{3-2},
the conservative solution contains two peakons emerging from
the point $\bar q$.  As $t\to \tau +$, their strengths and positions
satisfy again (\ref{ch1-condp}), while (\ref{ch1-condcoeff}) is replaced by
\begin{equation}
q_i(t)\to \bar q\,,
\quad q_j(t)\to \bar q\,,\qquad q_i(t)>q_j(t)
~~\hbox{for} ~~t>\tau\,.
\label{ch1-peakonuscenti}
\end{equation}
The vanishing viscosity approach in \cite{CHK1, CHK2} singles out
the dissipative solutions.  These can also be characterized
by the Oleinik type estimate
$$u_x(t,x)\leq C(1+t^{-1})\,,$$
valid for $t>0$ at a.e.~$x\in\R$.
On the other hand, the coordinate transformation approach
in \cite{BC2} and the present one, based on optimal transport metrics,
appear to be well suited for the study of
both conservative and dissipative solutions.

In the following chapters we focus on conservative solutions
to the Camassa-Holm equation.
We start with the study of the spatially periodicity because it allows us to concentrate
on the heart of the matter, i.e.~the uniqueness and stability
of solutions beyond the time of singularity formation.
It will spare us some technicalities, such as the analysis of the
tail decay of $u, u_x$ as $x\to\pm\infty$. In this respect we shall discuss the decay analysis of solutions in Section \ref{decayinfty} of  Chapter \ref{chap3}.

The main ingredients
can already be found in the paper \cite{BC1}, devoted to dissipative solutions
of the Hunter-Saxton equation.

As initial data, we take
\begin{equation}
u(0,x)=\bar u(x)\,,
\label{ch1-initialcond}
\end{equation}
with $\bar u$ in the space $H^1$ of
absolutely
continuous
functions $u$ with derivative $u_x\in L^2$.
To fix the ideas, we assume that the period of a spatially periodic function in the space $H^1_\per$ is $1$, so that
$$
u(x+1)=u(x)\qquad\qquad ~x\in\R\,.
$$
On $H^1_\per $ we shall use the norm
$$\big\|u\|_{H^1_\per }\doteq \left(\int_0^1\big|u(x)\big|^2\,dx
+\int_0^1\big|u_x(x)\big|^2\,dx\right)^{1/2}.$$
\subsection{The main results}
In this section we state the main results of Part I of this thesis. We shall write them for the spatially periodic case.

\begin{theorem}
\label{theo1} 
For each initial
data $\bar u\in H^1_\per$, there exists
a solution $u(\cdot)$ of the Cauchy
problem (\ref{ch1-equation}), (\ref{ch1-initialcond}). Namely, the map
$t\mapsto u(t)$ is Lipschitz continuous
from $\R$ into $L^2_\per$, satisfies (\ref{ch1-initialcond}) at time $t=0$,
and the identity
\begin{equation}
\frac d{dt} u = -uu_x -P_x
\label{ch1-diffl2}
\end{equation}
is satisfied as an equality between elements in $L^2_\per $
at a.e.~time $t\in\R$. This same map $t\mapsto u(t)$ is continuously
differentiable from $\R$ into $L^p_\per$ and satisfies (\ref{ch1-diffl2})
at a.e.~time $t\in\R$, for all $p\in [1,2[\,$.
The above solution is conservative in the sense that, for a.e.~$t\in\R$,
\begin{equation}
E(t)=\int_0^1 \big[u^2(t,x)+u_x^2(t,x)\big]\,dx=E^{\bar u}\doteq
\int_0^1 \big[\bar u^2(x)+\bar u_x^2(x)\big]\,dx\,.
\label{ch1-energy}
\end{equation}
\end{theorem}
\begin{theorem}
\label{theo2}
Conservative solutions to (\ref{ch1-equation}) can be constructed
so that they constitute a continuous flow $\Phi$. Namely, there exists
a distance functional $J$ on $H^1_\per $ such that
\begin{equation}
\frac 1C\cdot\|u-v\|_{L^1_\per}\leq J(u,v)\leq C\cdot
\|u-v\|_{H^1_\per }
\label{ch1-weaktop}
\end{equation}
for all $u,v\in H^1_\per$ and some constant $C$ uniformly valid on
bounded sets of $H^1_\per\,$.
Moreover, for any two solutions
$u(t)=\Phi_t\bar u$, $v(t)=\Phi_t\bar v$ of (\ref{ch1-equation}),
the map $t\mapsto J\big(u(t),\,v(t)\big)$
satisfies
\begin{equation}
J\big(u(t),\,\bar u\big)\leq C_1\cdot |t|\,,
\label{ch1-contu}
\end{equation}
\begin{equation}
J\big(u(t),\,v(t)\big)~\leq~ J(\bar u,\bar v)\cdot e^{C_2 |t|}
\label{ch1-stabsol}
\end{equation}
for a.e. $t\in\R\,$ and constants $C_1,C_2$, uniformly valid as
$u,v$ range on bounded sets of $H^1_{per}\,$.
\end{theorem}
\v
The previous results can be extended to the following space of functions which exponential decay: let $\alpha \in ]0,1[$, then
$$
X_\alpha \doteq \left\{u\in H^(\R) \st \int_\R[u^2(x)+u_x^2(x)]e^{\alpha|x|}\right\}.
$$ 
It is not so restrictive one can think, in fact the peakon functions (see Section \ref{3-2}), the natural solitary waves of the Camassa-Holm equation which have the \emph{soliton} properties, belong to it. 
\v
Somewhat surprisingly, all the properties stated in Theorem \ref{theo1}
are still not strong enough
to single out a unique solution.  To achieve uniqueness, an
additional condition is needed.
\begin{theorem}
\label{theo3}
Conservative solutions $t\mapsto u(t) $ of (\ref{ch1-equation})
can be constructed with the following additional property:
\v
\n For each $t\in\R$, call $\mu_t$ the absolutely continuous
measure having density $u^2+u_x^2$ w.r.t.~Lebesgue measure.
Then, by possibly redefining $\mu_t$ on a set of times of measure zero,
the map $t\mapsto\mu_t$ is continuous w.r.t.~the topology
of weak convergence of measures. It provides a measure-valued
solution to the conservation law
$$
w_t+(uw)_x=(u^3-2uP^u)_x\,.
$$
\v
The solution of the Cauchy problem (\ref{ch1-equation}), (\ref{ch1-initialcond}) satisfying the properties
stated in Theorem \ref{theo1} and this additional condition is unique.
\end{theorem}
\v
In Section \ref{3-2} we derive some
elementary properties
of multi-peakon solutions and show that any initial data can be
approximated in $H^1_\per$ by a finite sum of peakons.
In Section \ref{3-3} we introduce our distance functional $J(u,v)$
and study its
relations with other distances defined by Sobolev norms.
The continuity of the flow (\ref{ch1-equation}), together with the key
estimates (\ref{ch1-contu})-(\ref{ch1-stabsol})
are then proved in the following two sections.
The proofs of Theorems \ref{theo1} and \ref{theo2} are completed in
Section \ref{3-6}.  The uniqueness result stated in Theorem \ref{theo3} is proved
in Section \ref{3-7}.  As a corollary, we also show that in a multi-peakon
solution the only possible interactions
involve exactly two peakons: one positive and one negative.
In particular, no triple interactions can ever occur.
\v

Now we dedicate the rest of this chapter to exhibit an example which will be the start point of the technique we shall develop in Chapter \ref{chap3}. In particular, we shall see how the \emph{optimal transportation theory} can be useful nonlinear equation like Camassa-Holm equation. The key point is to define a \emph{Monge-Kantorovich} like metric for a space of Radon measures. In Section \ref{HSexample} we start from the \emph{Hunter-Saxton} equation \cite{BC1} for give a brief heuristic idea for how this technique is involved in.

\section{The Hunter-Saxton equation\label{HSexample}}
The Hunter-Saxton equation describes the propagation of waves in a massive vector field of a nematic liquid crystal. Since the physical interpretation and its derivation are beyond to the description of this thesis, we refer to \cite{HS, BC1}. It can be written in a non-local formulation as a conservation law with a source term:
\begin{equation} 
u_t+\left(\frac {u^2}2\right)_x=\frac 14 \left(\int_{-\infty}^x-\int^{+\infty}_x\right)u_x^2(t,y)\,dy\doteq Q^u(t,x)
\label{HunteSaxton}
\end{equation}
where
\begin{itemize}
\item $t\geq 0$ is the time variable,
\item $x\in \R$ is the space variable in a reference frame, 
\item $u(t,x)\in\R$ is related to the orientation of the liquid crystal molecules in the position $x$ at time $t$.
\end{itemize}
Suppose that there exists a smooth solution to (\ref{HunteSaxton}). To the Hunter-Saxton equation we can associate the following two conservation laws  
\begin{eqnarray}
&&(u_x)_t+(u u_x)_x=\frac{u_x^2}2,\label{clxux}
\\
&&(u^2_x)_t+ (u u_x^2)_x=0,\label{clxux2} 
\end{eqnarray}
which are obtained by computing the derivative of the equation (\ref{HunteSaxton}) w.r.t. the $x$ variable and then multiplying it by $u_x$ to achieve the second one.

A further conservation law is satisfied by the source term $Q^u(t,x)$. Since for smooth solutions (\ref{clxux2}) yields the conservation of the energy 
$$
E_0 \equiv E(t)\doteq \int_\R u_x^2(t,x)\,dx 
$$
the function $Q^u$ can be expressed in the following way
$$
Q(t,x)=-\frac{E_0}4 + \frac 12 \int_{-\infty}^x u_x^2(t,y)\,dy.
$$ 
By deriving w.r.t. $t$ we obtain
\begin{equation}
\label{conseqQ}
Q_t+u Q_x = 0.
\end{equation}
The function $Q^u$ is constant along the characteristic curves $\xi_u(t,y)$ defined by 
\begin{equation}
\label{characteristics}
\frac \partial{\partial t} \xi_u(t,y)=u(t,\xi_u(t,y)),\qquad \xi_u(0,y)=y.
\end{equation}
Let us remark that the previous equations holds for all the time $t$ in which $u$ is a classical solution. It can be seen by the method of characteristics that if $u_0\not\equiv 0$ is a smooth initial data and for some $x_0$ we have ${u_0}_x(x_0)<0$, along its outgoing characteristic the gradient blows up in finite time.  Since the quantity $E(t)$ remains bounded also at the time of blow up, we can think that a finite amount of energy will be concentrated at the point of blow-up.
In \cite{BC1} the authors focus their attention on solutions that \emph{dissipate} this quantity of energy. As far as the conservative solution is concerned, equation (\ref{clxux2}) will be satisfied in sense of measures, i.e. it means that thinking at a measure $\mu_t$ with absolutely continuous part which satisfies $d\mu_t = u_x^2(t,\cdot) d\mathcal L$ (here with $\mathcal L$ we indicate the Lebesgue measure), it satisfies
$$
\partial_t \mu_t +\partial_x(u \mu_t) = 0 \qquad {\rm in }\ \mathcal D'.
$$
To find a \emph{conservative solution} to the Cauchy problem (\ref{HunteSaxton}) with finite energy smooth initial condition $u_0$ means then to find a couple $(u(t),\mu_t)$ which satisfies the following system of conservation laws: let $\mu_0$ be the absolutely continuous measure w.r.t. Lebesgue measure defined by $d\mu_0 = {u_0}_x^2 d\mathcal L$, then
\begin{equation}
\left\{
\begin{array}{ll}
\dis \partial_t u + \partial_x \left(\frac {u^2}2\right) =- \frac 14 \mu_t(\R) + \frac 12 \mu_t(]-\infty,x])
\qquad &u(0,x)=\bar u(x),
\\
\dis \partial_t \mu_t +\partial_x(u \mu_t) = 0 \qquad &\left. 
\mu_t\right|_{t=0}=\mu_0.
\end{array}
\right.
\end{equation}
Due to the nonlinearity of the problem, as shows \cite[Example 2]{BC1} for the dissipative solution of the Hunter-Saxton equation, we can aspect that the usual ``strong'' distance stemming from convex norm is not useful in order to construct a continuous semigroup of solutions. In the following section, we give a sketch of the construction of a metric which yields continuity of solution with respect to the initial data.
\subsection{A transportation map\label{trapsit}}
Let $u_0$ and $v_0$ be two initial data whose associated measures $\mu_0^1$ and $\mu_0^2$ have the same total mass $\mu_0^1(\R)=\mu_0^2(\R).$  Suppose that such a initial data are not constant in any interval of $\R$, so that the functions $Q^{u_0}$ and $Q^{v_0}$ are absolutely continuous and increasing.

\psfrag{x}{$x$}
\psfrag{psix}{$\Psi_0(x)$}
\psfrag{qu}{$Q^{u_0}$}
\psfrag{qv}{$Q^{v_0}$}
\psfrag{E0}{$\frac 14 E_0$}
\psfrag{-E0}{$-\frac 14 E_0$}
\begin{figure}[ht]
\centerline{\includegraphics[width=7cm]{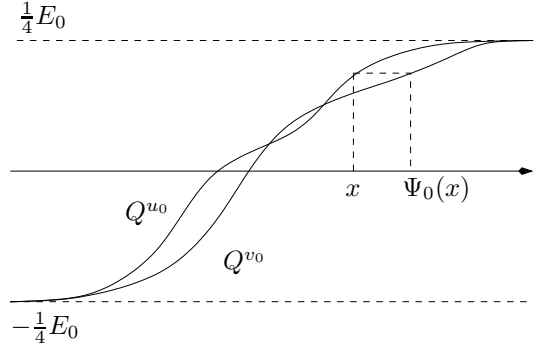}}
\caption{The map $\Psi_0$}\label{figQuQv}
\end{figure}

We can thus define a continuous map $\Psi_0$ for which at every $x\in \R$ it associates the unique point $\Psi_0(x)$ such that 
\begin{equation}
Q^{u_0}(x) = Q^{v_0}(\Psi_0(x))
\label{depsi}
\end{equation} 
(see Figure \ref{figQuQv}). 
Let us Remark that the $\Psi_0$ is an increasing function, in fact since $Q^{u_0}$ and $Q^{v_0}$ are increasing, then $x<y$ implies 
$$
Q^{v_0}(\Psi_0(x))=Q^{u_0}(x)<Q^{u_0}(y)=Q^{v_0}(\Psi_0(y)),
$$ 
so $\Psi_0(x)<\Psi_0(y)$ holds.

Now we want to see how the map $\Psi_0$ evolves in time.  Since by (\ref{conseqQ}) $Q^u$, $Q^v$ are conserved along the caracteristic curves (\ref{characteristics}), the equalities
$$
Q^u(t,\xi_u(t,x)) = Q^{u_0}(x),\qquad  Q^v(t,\xi_v(t,x))=Q^{v_0}(x)\qquad \mbox{for all $x\in \R$}
$$
yields the definition of the \emph{transportation map} 
\begin{equation}
\label{depsit}
\Psi_t(\xi_u(t,x))\doteq \xi_v(t,\Psi_0(x)).
\end{equation}
\subsection{The stability issue of a system of ODE}
As in \cite[Section 3]{BC2}, we compute a change of variables in order to obtain a system of ODE with Lipschitz vector field. Let suppose that the initial data $u_0$ is in $H^1(\R)$. 
Let set
$$ 
\omega\doteq 2\arctan (u_x),
$$
then $\omega$ belong to the unit circle $\T\doteq [0,2\pi]$, with $0$ and $2\pi$ identified.
Computing the derivative of $\omega$ along the characteristics, having in mind the equation (\ref{clxux}) for $u_x$ we have
$$
\frac d{dt} \omega(t,\xi(t,y))= 2 \frac {(u_x)_t + u (u_x)_x}{1+u_x^2} =\frac {-u_x^2}{1+u_x^2} = \frac{-\tan^2(\omega/2)}{1+\tan^2(\omega/2)} = -\sin^2(\omega/2).
$$
It is thus natural to consider the new unknowns which take values in the space $\R^2\times \T$, 
$$
\mathbf X^u=\mathbf X^u(t,y)\doteq \left(
\begin{array}{c}
\xi_u(t,y)
\\
u(t,\xi_u(t,y))
\\
\omega_u(t,\xi_u(t,y))
\end{array}
\right)
$$
and write the corresponding Cauchy problem
\begin{equation}
\label{CPX}
\frac d{dt}
\mathbf X^u(t,y)=
\left(
\begin{array}{c}
u(t,\xi_u(t,y))
\\
Q^u(t,\xi_u(t,y))
\\
\omega_u(t,\xi_u(t,y))
\end{array}
\right)
=
f(\mathbf X^u(t,y))
\end{equation}
with initial data
\begin{equation}
\label{IDX}
\mathbf X^u(0,y)=
\left(
\begin{array}{c}
y
\\
u_0(y)
\\
2\arctan ({u_0}_x(y))
\end{array}
\right).
\end{equation}
We remark that by definition, $Q^u$ is far from to be Lipschitz continuous, then we cannot suppose that in the previous system of ODE the function $f$ is a Lipschitz vector field. Hence, to overcome this lack of Lipschitz continuity, we shall make use of the function $\Psi_t$ introduced in \ref{trapsit}. Let $\mathbf X^u$ and $\mathbf X^v$ be two solution of the Cauchy problem (\ref{CPX})-(\ref{IDX}), corresponding to the initial data $u_0$, $v_0$ respectively.  We allow ourselves to make an abuse of notation by defining the map $\Psi_t$ in the following way:
$$
\Psi_t(\mathbf X^u(t,y))\doteq \mathbf X^v(t,\Psi_0(y)).
$$
We gain a sort of Lipschitz continuity for the function $f$ if we restrict it on the manifold located by $\Psi_t$.  
Let us compute the difference of the vector field $f$ evaluated in the points $\mathbf X^u(t,y)$ and $\Psi_t(\mathbf X^u(t,y))$. 
$$
|f(\mathbf X^u(t,y))-f(\Psi_t(\mathbf X^u(t,y)))|=
\left|
\left(
\begin{array}{c}
u(t,\xi_u(t,y))-v(t,\xi_v(t,\Psi_0(y)))
\\
Q^u(t,\xi_u(t,y))-Q^v(t,\xi_v(t,\Psi_0(y)))
\\
\omega_u(t,\xi_u(t,y))-\omega_v(t,\xi_v(t,\Psi_0(y)))
\end{array}
\right)
\right|
$$
since by definition of the map $\Psi_0$ and by the equation (\ref{conseqQ}) we have the identity
$$
Q^u(t,\xi_u(t,y))-Q^v(t,\xi_v(t,\Psi_0(y)))\equiv 0
$$
we deduce the following estimate for the vector field $f$
$$
|f(\mathbf X^u(t,y))-f(\Psi_t(\mathbf X^u(t,y)))|\leq | \mathbf X^u(t,y) - \mathbf X^u(t,y)|.
$$
From the previous inequality we can prove that the difference of the two solutions $X^u(t,y)$ and $X^v(t,\Psi_0(y))$ can be estimated by the initial data. In fact, Gronwall Lemma applied to the inequality
$$
\begin{array}{rl}
\dis
|\mathbf X^u(t,y)-\mathbf X^v(t,\Psi_0(y))|
\leq & \dis|\mathbf X^u(0,y)-\mathbf X^v(0,\Psi_0(y))|
\\
&\quad \dis +\int_0^t |f(\mathbf X^u(s,y))-f(\mathbf X^v(s,\Psi_0(y)))|\,ds
\end{array}
$$
yields the estimate
$$
|\mathbf X^u(t,y)-\mathbf X^v(t,\Psi_0(y))|\leq e^{t}|\mathbf X^{u_0}(y)-\mathbf X^{u_0}(\Psi_0(y))|.
$$
The previous inequality suggests how to introduce a new distance in order to obtain a stability result for solutions of the Hunter-Saxton equation. 
For every $u\in H^1(\R)$ let define the measure 
$$
\mu^u(A)\doteq \int_{\{x\in \R\,:\, (x,u(x),\omega(x))\in A\}} u_x^2(x)\,dx
$$
for every Borel set $A\subset \R^2\times \T$.
The function $\Psi_0$ can be regarded as a transportation map which transports the measure $\mu^u$ into the measure $\mu^v$.
The distance we shall introduce in Chapter \ref{chap3} will be thus a sort of Wasserstein distance between measure \cite{V}. 

\chapter{Conservative multi-peakon solution\label{multipeak}}
\begin{figure}
\centerline{\includegraphics{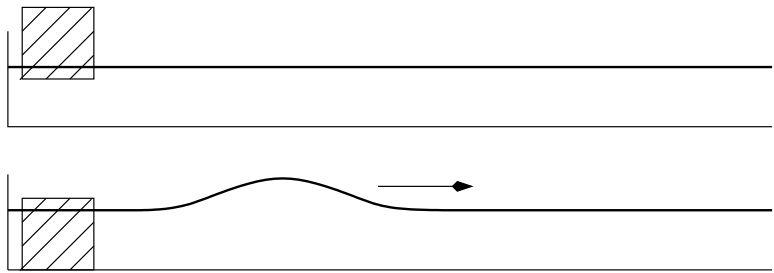}}
\caption{Diagram of Scott Russell's experiment to generate a solitary wave\label{SRussell}}
\end{figure}

\begin{quotation}
``I was observing the motion of a boat which was rapidly drawn along a narrow channel by a pair of horses, when the boat suddenly stopped - not so the mass of water in the channel which it had put in motion; it accumulated round the prow of the vessel in a state of violent agitation then suddenly leaving it behind, rolled forward with great velocity, assuming the form of large solitary elevation, a rounded, smooth and well defined heap of water, which continued its course along the channel apparently without change of form or diminution of speed. I followed it on horseback and overtook it still rolling on at a rate of some eight or nine miles an hour, preserving its original figure some thirty feet long and a foot to a foot and a half in height. its height gradually diminished and after a chase of one or two miles I lost it in the windings of the channel. Such in the month of August 1834 was my first chance interview with that singular and beautiful phenomenon which I have called the Wave of Translation...''
\vskip 10pt
\begin{minipage}[t]{9.5 cm}\raggedleft
John Scott Russell, 1844
\end{minipage}
\vskip 10pt
\end{quotation}
It was Scott Russell \cite{R} who introduces the concept of \emph{solitary waves}  to indicate no more than wave which propagate without change of form and have some localized shape (see also, \cite{DEGM,DJ}). His experiment is described in fig. \ref{SRussell}. For more than sixty years it was only a pure scientific curiosity, until Korteweg and de Vries \cite{KdV} derived the equation for the propagation of waves in one direction on the surface of a shallow canal. The profile of the travelling wave solution with permanent shape found by Miura \cite{M2} is precisely the shape of the wave which Scott Russel observed in his experiments. The term \emph{soliton} is substantially different from solitary wave, it was introduced in 1965 by Zabusky and Kruskal \cite{ZK} to indicate waves that whenever collide each other they do not break up and disperse, but remains almost identical to a solitary wave solution.
In this chapter we investigate the shape of solitary waves for the Camassa-Holm equation (\ref{ch1-equation}).
\section{Multi-peakon solution in the real line\label{3-2}}
\label{ODEsystem}
In this section we shall construct a solution of the Camassa-Holm equation starting from an initial condition $u_0$ of the form
$$
u_0(x)=\sum_{j=1}^{N} p_j e^{-|x-q_j|}\,.
$$
The motivation of this choice is given by the shape of traveling wave solution (see \cite{CH} and \cite[Example 5.2]{CE1}).  Looking for solution of the equation (\ref{ch1-equation}) in the traveling wave form $u(t,x)=U(x-ct)$, with a function $U$ that vanishes at infinity, the limit of $\kappa\to 0$ leads to the function $U= c e^{-|x-ct|}$. This is not a solitary wave in the sense introduced by Scott Russell because of the presence of the cuspid at the position $x=ct$. However, the evolution of an initial data like $u_0$ remains of the same shape \cite{CH,CS}. It is a superposition of peaked solitary waves, which evolves as
$$
u(t,x)=\sum_{j=1}^{N} p_j(t) e^{-|x-q_j(t)|}\,.
$$
Hence we term \emph{peakon} a peaked solitary wave to emphasize the soliton properties of such a function.

As long as the classical solution of the problem 
\begin{equation}
\label{HSYS}
\left\{
\begin{array}{l}
\dis \dot q_i=\sum_{j=1}^N p_j e^{-|q_i-q_j|}\,,
\\
\dis \dot p_i=p_i\sum_{j=1}^N p_j \sign(q_i-q_j)e^{-|q_i-q_j|}
\end{array}
\right.
\end{equation}
exists, the solution of this system gives the coefficients ${\bf p }(t)=(p_1,\dots, p_{N})$ and ${\bf q}(t)=(q_1,\dots, q_{N})$ for the solution $u(t,x)$ to the Camassa-Holm equation. 
Let us observe that the previous system can be viewed as an Hamiltonian system with Hamiltonian function $H({\bf q},{\bf p })=\frac 12~\sum_{i,j}p_i p_j e^{-|q_i-q_j|}$.
 
\begin{figure}
\centerline{\includegraphics[height = 7cm]{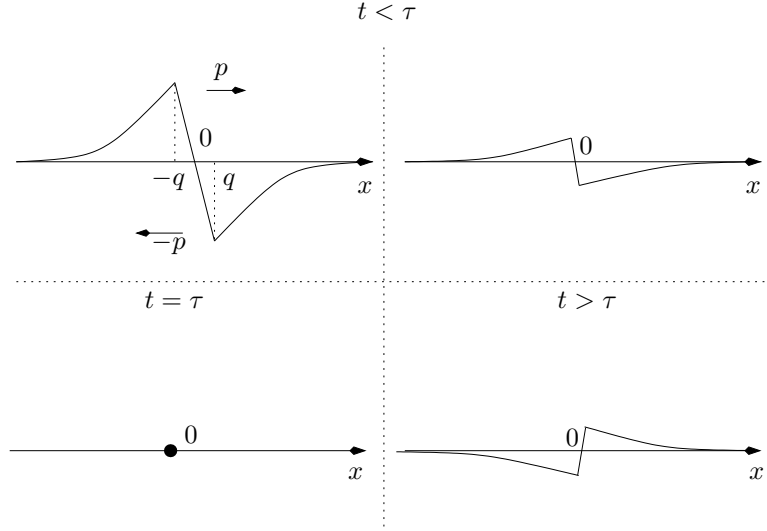}}
\caption{Peakon-antipeakon interaction\label{figpeakapeak}}
\end{figure}

In \cite{HR} the authors prove the existence of a global multipeakon solution when strengths $p_i$ are positive for all $i=1\dots N$ and the convergence of the sequence of multipeakon solution. If $u_0$ is an initial data such that the distribution $u_0~-~{u_0}_{xx}$ is a \emph{positive} Radon measure, there exists a sequence of multipeakons that converges in $L^\infty(\R,H^1_{loc}(\R))$. In this case the crucial fact is that no interaction between the peakons occurs, and then the gradient remains bounded. However, a general initial data contains both positive and negative peakons, as in the so called peakon-antipeakon interaction: one positive peakon with strength $p$, centered in $-q$, moves forward and one negative anti-peakon in $q$, with strength $p$ moves backward (fig. \ref{figpeakapeak}). The evolution of the system produces the overlapping of the two peakons at finite time $t=\tau$, so that $q\to 0$.  The conservation of the energy $E=H({\bf q}(t),{\bf p }(t))$ yields
\begin{equation}
\label{zetalimit}
E=\lim_{t\to \tau^-} p^2(1-e^{-2|q|})\,.
\end{equation}
and then the quantity $p$ blows up in finite time. At the point $(\tau,0)$ occurs thus a singularity for the solution $u$.  To extend the solution also after the interaction time with a solution which conserves the energy $E$ we can think that at the interaction point an antipeakon/peakon couple emerge, the first, negative, moving backward and the second, positive, moving forward with coefficients $(-q, -p)$ and $(q,p)$. According to the conservation of the energy, the choice of $q$ and $p$ must satisfy (\ref{zetalimit}) as $t\to \tau^+$.
\psfrag{0}{$0$}\psfrag{p}{$p$}\psfrag{-p}{$-p$}\psfrag{q}{$q$}\psfrag{-q}{$-q$}
\psfrag{t<tau}{$t<\tau$}\psfrag{t=tau}{$t=\tau$}\psfrag{t>tau}{$t>\tau$}
It  yields a change of variables which resolves the singularity at $(\tau,0)$
$$
\zeta\doteq p^2 q\qquad \omega\doteq\arctan (p)
$$ 
with this choice, the Hamiltonian system leads to the ODE
$$
\frac d{dt}
\left(
\begin{array}{c}
\zeta
\\
\omega
\end{array}
\right)
=f(\zeta,\omega)\,,
\qquad
\left(
\begin{array}{c}
\zeta
\\
\omega
\end{array}
\right)(\tau)
=
\left(
\begin{array}{c}
\frac E 2
\\
\frac \pi 2
\end{array}
\right)
$$
with 
$$
f(\zeta,\omega)= 
\left(
\begin{array}{c}
\left[1- e^{-\zeta \cot^2(\omega)}-\zeta\cot^2(\omega)e^{-\zeta \cot^2(\omega)}\right]\tan^3(\omega)
\\
\sin^2(\omega)e^{-\zeta \cot^2(\omega)}
\end{array}
\right)
$$
and $f$ is a Lipschitz vector field in a neighborhood of the point $( \frac E 2,\frac \pi 2)$. The solution $(\zeta(t),\omega(t))$ of this problem provides then the unique couple $(q(t),p(t))$ which coincides with the classical solution of the Hamiltonian system for $t<\tau$ and extends it for $t\geq \tau$.

This example suggests the way to construct the multipeakon solution whenever an interaction between peakons occurs (see also \cite{BF2} for an ``energetic'' motivation). Suppose that two or more peakons with strengths $p_1,\dots,p_k$ annihilate at the position $\bar q$ at time  $\tau$  and produce a blow up of the gradient $u_x$. The conservation of the energy yields that there exists and is positive the limit 
$$
e_{\tau}\doteq \lim_{t\to \tau^-} \int_{\xi^-(t)}^{\xi^+(t)} u_x^2(t,x)\,dx
$$ 
where $\xi^-$ and $\xi^+$ are the smallest and the largest characteristic curve passing through the point $(\tau,\bar q)$. Assume that after the interaction two peakons appear with strengths $p_1,\,p_2$ and placed at the position $q_1,\,q_2$.
Let consider the change of variables
$$
z=p_2+p_1\quad w=2\arctan (p_2-p_1)\quad \eta = q_2+q_1 \quad \zeta=(p_2-p_1)^2(q_2-q_1),
$$
then the system (\ref{HSYS}) turns out to be
$$
\begin{array}{rl}
\dot w=&\!\! -\left[\sin (w) \cosh \Big(\frac {\zeta}{2\tan^2 (w/2)}\Big)
+2z \sinh\Big(\frac {\zeta}{2\tan^2 (w/2)}\Big) \right]\cdot\sum\limits_{j\geq k+1} p_j e^{-q_j+ \eta /2}
\\
&+[z^2 \cos^2(w/2) - \sin^2 (w/2)]e^{-\frac{\zeta}{\tan^2 (w/2)}}
\\
\dot z= &\!\!-\left[\frac 12\sin(w) \sinh \Big(\frac {\zeta}{2\tan^2 (w/2)}\Big)+ z\cosh \Big(\frac {\zeta}{2\tan^2 (w/2)}\Big) \right]\cdot \sum\limits_{j\geq k+1} p_j e^{-q_j}
\\
\dot \eta= &\!\! z[1+ e^{-\frac{\zeta}{\tan^2 (w/2)}}]+ 2 \cosh \Big(\frac {\zeta}{2\tan^2 (w/2)}\Big)\cdot \sum\limits_{j\geq k+1} p_j e^{-q_j}
\\
\dot \zeta = &\!\!\frac{z^2\zeta}{\tan (w/2)}e^{-\frac{\zeta}{\tan^2 (w/2)}}-\tan^3(w/2) \left(1-e^{-\frac{\zeta}{\tan^2 (w/2)}}- \frac{\zeta}{\tan^2 (w/2)}\right) +
\\
&
+ 2\zeta\left[ \sinh \Big(\frac{\zeta}{2\tan^2 (w/2)}\Big)\cdot \left(\frac {\tan^2 (w/2)}{\zeta}-
\frac{z}{\tan (w/2)} \right)\right.+
\\
&\qquad\qquad\qquad\qquad\qquad\qquad
 \left. -\cosh\Big( \frac {\zeta}{2\tan^2 (w/2)}\Big)\right]\cdot \sum\limits_{j\geq k+1}p_j e^{-q_j+\eta/2} 
\end{array}
$$
$$
\begin{array}{rl}
\dot p_i = & p_i e^{-q_i+\eta/2}\left[ z \cosh\Big( \frac {\zeta}{2\tan^2 (w/2)}\Big)+ \tan(w/2) \sinh\Big(\frac {\zeta}{2\tan^2 (w/2)}\Big)\right] + 
\\
&+\sum\limits_{j\geq k+1} p_i p_j \sign(q_i-q_j) e^{-|q_i-q_j|}
\\
\dot q_i= & e^{-q_i+\eta/2}\left[ z \cosh\Big( \frac {\zeta}{2\tan^2 (w/2)}\Big)+ \tan(w/2) \sinh\Big(\frac {\zeta}{2\tan^2 (w/2)}\Big)\right]+
\\&+\sum\limits_{j\geq k+1} p_j e^{-|q_i-q_j|}
\end{array}
$$
which is a system of ODE with locally Lipschitz continuous right hand side that can be extended smoothly also at the value $w=\pi$. The initial data become
$$
z(\tau)=\lim_{t\to \tau^-} \sum_{i=1}^k p_i(t)\qquad w(\tau)=\pi \qquad  \eta(\tau)=2\bar q\qquad \zeta(\tau)= e_\tau 
$$
$$
p_i(\tau)=\lim_{t\to \tau^-}p_i(t)\qquad q_i(\tau)=\lim_{t\to \tau^-} q_i(t) \qquad i=k+1,\dots, N
$$
Thus there exists a unique solution of such a system  which provides a multipeakon solution defined on some interval $[\tau,\tau'[$, up to the next interaction time.
As we will show in Corollary \ref{cormax2peak}, since Camassa-Holm equation is time reversible, once we prove the uniqueness of the  solution of a Cauchy problem, we have that maximal the number of peakons interaction is actually $k=2$, one with positive strength the other with negative one.
\section{Approximation of the initial data}
In this section we shall construct an approximation with initial data with a multipeakon function. Our aim is to approximate it with a sequence $u_\varepsilon$ which has an exponential decay at infinity uniformly w.r.t $\varepsilon$.  
\begin{lemma}
Let $f\in X_\alpha$. Then for every $\varepsilon>0$ there exists a multipeakon function $g$ of the form
$$
g(x)=\sum_{i=1}^N p_i e^{-|x-q_i|}
$$
such that 
\begin{eqnarray}
&&\|f-g\|_{\H}<\varepsilon
\\
&&\int_\R [g^2(x)+{g}_x^2(x)]e^{\alpha|x|}\,dx\leq C_0 
\end{eqnarray}
for some constant $C_0>0$ which does not depend on $\eps$.
\end{lemma}
\begin{proof}
Let $\rho(x)\in \mathcal C^\infty_0$ be a cut-off function such that
\begin{itemize}
\item $\rho(x)\geq 0$
\item $\rho(x)=1$ for every $|x|\leq 1$, $\rho(x)=0$ for every $|x|> 2$
\item $\int_\R \rho(x)\,dx =1$ 
\end{itemize}
and $\rho_\eps(x)\doteq \frac 1\eps \rho(\frac x\eps) $ be a mollifiers sequence. Observe that for every $\eps>0$, $f_\eps(x)\doteq \rho_\eps*f(x)$ is a smooth function which approximates the function $f$ in $H^1-$norm
\begin{equation}
\|f- f_\eps\|_{\H}<C\eps
\label{apprcinf}
\end{equation}
 moreover it belongs to $X_\alpha$, indeed
\begin{eqnarray*}
\int_\R [f_\eps ^2(x)+{f_\eps}_x^2(x)]e^{\alpha|x|}\,dx
&\leq& \int_\R \left[\int_\R (f^2(y)+f_x^2(y)) \rho_\eps(x-y) \,dy\right]e^{\alpha|x|}\,dx
\\
&\leq& \int_\R [f^2(y)+f_x^2(y)]\int_\R \rho_\eps(x-y) e^{\alpha|x|}\,dx \, dy
\\
&\leq& \int_\R [f^2(y)+f_x^2(y)] C_0 e^{\alpha|y|}\,dy =C_0 C^{\alpha,f}<\infty
\end{eqnarray*}
and $C_0$ is a constant which does not depend on $\eps$. From the previous inequality we can assert that for every $R>0$ one has
$\|f_\eps\|_{H^1(\R\setminus [-R,R])}\leq C_0 C_\alpha e^{-\alpha R}$ uniformly in $\eps>0$. We can choose thus $R_\eps$ big enough in order to have  
\begin{equation}
\|f_\eps\|_{H^1(\R\setminus [-R_\eps,R_\eps])}<\eps/2.
\label{fuori} 
\end{equation}

In the space $H^1([-R_\eps,R_\eps])$ we can approximate $f_\eps$ with a multipeakon function.  By using the identity
$$
\frac 12\left(I-\frac{\partial^2}{{\partial x}^2}\right)e^{-|x|}=\delta_0 \qquad in \mathcal D'
$$
the function $f_\eps$ can be rewritten in convolution form
$$ 
f_\eps= e^{-|x|}*\left(\frac{f_\eps -  {f_\eps}_{xx}}2 \right)=\int_\R e^{-|x-y|}\cdot \frac{f_\eps(y) -  {f_\eps}_{xx} (y)} 2 \,dy\,.
$$
In the interval $[-R_\eps, R_\eps]$ the previous integral can now be approximated with a Riemann sum 
$$
g(x)= \sum_{i=-N}^N p_i e^{-|x-q_i|},\qquad 
\left\{
\begin{array}l
\dis q_i= \frac i N R_\eps 
\\
\dis p_i=\int_{q_{i-1}}^{q_i}\frac{ f_\eps(y) - {f_\eps}_{xx}(y)}{2}\,dy\,.
\end{array}
\right.
$$
Choosing $N$ sufficiently large we obtain $\|f_\eps- g\|_{H^1([-R_\eps,R_\eps])}<\eps$. Together with (\ref{apprcinf}) and (\ref{fuori}) this last estimate yields the result.
\end{proof}
\section{Periodic multi-peakon}
By a periodic peakon we mean a function of the form
\begin{equation}
u(x)=p\, \chi(x-q)\,,\qquad\qquad
\chi(x)\doteq\sum_{n\in\Z} e^{-|x-n|}\,.
\label{ch1-perpeak}
\end{equation}
Observe that the periodic function $\chi$ satisfies
\begin{equation}
\begin{array}{l}
\displaystyle\chi(-x)=\chi(x)=\chi(x+1)\qquad\qquad x\in\R\,,
\\
\displaystyle\chi(x)
= \frac{e^x+e^{1-x}}{e-1}\qquad\qquad\quad\qquad x\in [0,1]\,.
\end{array}
\label{ch1-formachi}
\end{equation}
For future use, we observe that for every periodic function $u$, the convolution $P$, defined in (\ref{ch1-nonlocalP}),  takes the form
\begin{equation}
P(x)=\frac 12\int_0^1 \chi(x-y)\left(u^2(y)+\frac{u_x^2(y)}2 \right)\,dy
\label{Pperiodica}
\end{equation}

We begin this section by observing that also any periodic initial data
can be approximated by periodic multi-peakons.
\v
\begin{lemma}
\label{lem1}
Let $f\in H^1_\per$.
Then for any $\ve>0$
there exists periodic multi-peakon $g$,
of the form
$$
g(x)=\sum_{i=1}^{N} p_i\sum_{n\in\Z} e^{-|x-q_i-n|}=
\sum_{i=1}^N p_i\,\chi(x-q_i)
$$
such that
$$\|f-g\|_{H^1_\per}<\ve\,.$$
\end{lemma}
\v
\begin{proof} By taking a suitable mollification, we can
approximate $f$ with a periodic function $\tilde f\in\C^\infty$, so that
\begin{equation}
\|f-\tilde f\|_{H^1_\per}<\ve/2\,.
\label{ch1-apprx}
\end{equation}
Next, we observe that
$$\frac 12\left( e^{-|x|}-\frac{\partial^2}{\partial x^2} e^{-|x|}
\right)=\delta_0\,,$$
where $\delta_0$ denotes the Dirac distribution concentrating
a unit mass at the origin.
We can thus write $\tilde f$ as a convolution:
$$
\tilde f=\delta_0*\tilde f =\frac 12
\left( e^{-|x|}-\frac{\partial^2}{\partial x^2} e^{-|x|}
\right)*\tilde f = e^{-|x|} * \left(\frac{\tilde f-\tilde f''} 2\right)\,,
$$
$$\tilde f(x)=\int_0^1\chi(x-y)\cdot \frac{\tilde f(y)-\tilde f''(y)} 2\,dy\,.$$
The above integral
can now be approximated with a Riemann sum
$$
g(x)=\sum_{i=1}^N p_i\,\chi(x-q_i)\,,
\qquad\qquad p_i=\int_{(i-1)/N}^{i/N}\frac{\tilde f(y)-\tilde f''(y)}2 \,dy\,.$$
Choosing $N$ sufficiently large we obtain
$\|\tilde f-g\|_{H^1_{per}}<\ve/2$. Together with (\ref{ch1-apprx})
this yields the result.
\end{proof}
\vs
Next, we show how to construct a unique conservative solution,
for multi-peakon initial data.  As long as the locations
$q_i$ of the peakons remain distinct, this can be obtained by
solving the Hamiltonian system of O.D.E's (\ref{ch1-Hsys}).

However, at a time $\tau$ where two or more
peakons interact, the corresponding strengths
$p_i$ become unbounded.  A suitable transformation of
variables is needed, in order to resolve the singularity
and uniquely extend the solution beyond the interaction time.
\v
\begin{lemma}
\label{lem2}
Let $\bar u$ be any periodic, multi-peakon
initial data. Then the Cauchy problem (\ref{ch1-equation}), (\ref{ch1-initialcond})
has a global, conservative multi-peakon solution defined for all $t\in\R$.
The set $\I$ of times where two or more peakons interact is at most
countable.
Moreover, for all $t\notin\I$, the energy conservation (\ref{ch1-energy}) holds.
\end{lemma}
\begin{proof} The solution can be uniquely constructed by
solving the hamiltonian system (\ref{ch1-Hsys}), up to the first time
$\tau$ where two or more peakons interact.
We now show that there exists a unique way to prolong the solution
for $t>\tau$, in terms of two outgoing peakons.
To fix the ideas, call
$$\bar q=\lim_{t\to\tau-} q_i(t)\qquad\qquad i=1,\ldots,k\,,$$
the place
where the interaction occurs,
and let $p_1(t),\ldots, p_k(t)$ be the strengths of the
interacting peakons.
Later in Section \ref{3-7} we will show that only the case $k=2$
can actually occur, but at this stage we need to consider
the more general case.
We observe that the strengths $p_{k+1},\ldots,p_N$ of the peakons
not involved in the
interaction remain continuous at time $\tau$.   Moreover,
by (\ref{ch1-Hsys})
there exists the limit
$$\bar p=\lim_{t\to \tau-} \sum_{i=1}^k p_i(t)\,.$$
We can thus write
$$u(\tau,x)~=~\lim_{t\to\tau-}
\sum_{i=1}^N p_i(t)\,e^{-|x-q_i(t)|}~=~\bar p \,
e^{-|x-\bar q|}+\sum_{i=k+1}^N
p_i(\tau)\,e^{-|x-q_i(\tau)|}\,.$$
For $t>\tau$, we shall prolong the solution with two
peakons emerging from the point $\bar q$.  The strength
of these two peakons will be uniquely determined
by the requirement of energy conservation (\ref{ch1-energy}).

Call $\xi^-(t)$, $\xi^+(t)$ respectively
the position of the smallest and largest characteristic
curves passing through the point $(\tau,\bar q)$, namely
\begin{equation}
\begin{array}{l}
\xi^-(t)\doteq \min\Big\{ \mbox{$\xi(t)$; $\xi(\tau)=\bar q$,
$\dot\xi(s)=u\big(s,\xi(s)\big)$ $\forall s\in [\tau-h,\tau+h]$}\Big\}\,,\\
\xi^+(t)\doteq \max\Big\{ \mbox{$\xi(t)$; $\xi(\tau)=\bar q$,
$\dot\xi(s)=u\big(s,\xi(s)\big)$ $\forall s\in [\tau-h,\tau+h]$}\Big\}\,.\\
\end{array}
\label{ch1-charminmax}
\end{equation}
Moreover, define
$$
e_{(\tau,\bar q)}\doteq\lim_{t\to \tau-}
\int_{\xi^-(t)}^{\xi^+(t)} u_x^2(t,x)\,dx\,.
$$
The existence of this limit follows from the balance law
(\ref{ch1-uxquadro}).
This describes how much energy is concentrated at the interaction point.
For $t>\tau$ the solution will contain the peakons
$p_{k+1},\ldots,p_N$, located at $q_{k+1},\ldots, q_N$, together with
the two outgoing peakons $p_1,p_2$,  located at $q_1<q_2$.
The behavior of $p_i,q_i$ for $i\in\{k+1,\ldots, N\}$ is
still described by a system of O.D.E's as in (\ref{ch1-Hsys}).

However, to describe the evolution of $p_1,p_2,q_1,q_2$ one has to
use a different set of variables, resolving the singularity
occurring at $(\tau,\bar q)$.
As $t\to \tau+$ we expect
(\ref{ch1-condp}), (\ref{ch1-peakonuscenti}) to hold.
To devise a suitable set of rescaled variables,
we observe that, by (\ref{ch1-nonlocalequation}),
$$
\frac d{dt}\,u_x\big(t,\xi(t)\big)= -\frac 12
\,u^2_x\big(t,\xi(t)\big) + [u^2-P]
$$
along any characteristic curve $t\mapsto \xi(t)$
emerging from the point $\bar q$.
Since $u,P$ remain uniformly bounded, one has
$$u_x(t,x)\approx \frac 2{t-\tau}\qquad\qquad t>\tau\,,~~x\in \big[
q_1(t)\,,~q_2(t)\big]\,.$$
The total amount of energy concentrated in the interval between
the two peakons is given by
$$
\begin{array}{rl}
\displaystyle
\int_{q_1(t)}^{q_2(t)} \big[ u^2(t,x)+u_x^2(t,x)\big]\,dx
&\displaystyle 
\approx \left(\frac{u(q_2)- u(q_1)}{q_2-q_1}\right)^2\cdot (q_2-q_1)
\\
& \displaystyle\approx \frac{\Big[ (p_2-p_1)\big(1-e^{-|q_2-q_1|}\big)\Big]^2}{q_2-q_1}
\\
&
\displaystyle
\approx (p_2-p_1)^2(q_2-q_1)
\approx e_{(\tau,\bar q)}\,.
\end{array}
$$

The previous heuristic analysis suggests that, in order to
resolve the singularities, we should work with the variables
$$z=p_1+p_2\,,\quad w=2\arctan (p_2-p_1)\,,\quad
\eta=q_2+q_1\,,\quad \zeta =(p_2-p_1)^2(q_2-q_1)\,,$$
together with $p_{k+1}\,,\ldots,\,p_N\,,~~~q_{k+1}\,,\ldots, \, q_N$.

To simplify the following calculations we here
assume $0<q_1<q_2<q_{k+1}<...<q_N<1$,
which is not restrictive.

Let $\chi$ defined in (\ref{ch1-perpeak}) and $\tilde \chi(x)\doteq \frac{-e^x+e^{1-x}}{e-1}$, $x\in (0,1)$.
From the original system of equations (\ref{ch1-Hsys}) it follows
$$
\begin{array}{l}
\displaystyle\dot z=
\cosh\left(\frac \zeta 2\cot^2{w/2}\right)
z\cos^2{w/2}\sum_{j=k+1}^N p_j \chi\left(q_j-\eta/2\right)
\\
\qquad \dis
-
\sinh \left(\frac\zeta 2 \cot^2 w/2\right)
\tan w/2 
\sum_{j=k+1}^N p_j \tilde\chi\left(q_j-\eta/2
\right)
\\
\displaystyle\dot w=\left(z^2\cos^2 w/2-\sin^2 w/2 \right)
\chi\left(\zeta \cot^2 w/2\right) 
\\
\qquad\dis
+ 2
\cosh \left(\frac \zeta 2 \cot^2 w/2\right)z
\sum_{j=k+1}^N p_j \chi\left(q_j- \eta/2\right)
\\
\displaystyle\qquad
+
2
\sinh \left(\frac \zeta 2 \cot^2 w/2\right)\sin w
\sum_{j=k+1}^N p_j \tilde\chi\left(q_j-\eta/2\right)
\\
\displaystyle\dot \eta=
z\left[\chi(0)+\chi\left(\zeta \cot^2 w/2\right)\right]+
\cosh\left(\frac\zeta 2\cot^2 w/2\right)\sum_{j=k+1}^N p_j \chi\left(q_j-\eta/2\right)
\\
\displaystyle
\dot \zeta=
\left[\chi(0)-\chi\left(\zeta \cot^2 w/2\right)\right]\zeta\tan^2 w/2+
\chi\left(\frac \zeta \cot^2 w/2\right) z^2\zeta\cot w/2
\\
\qquad
\dis -\sinh\left(\frac \zeta 2 \cot^2 w/2 \right)  \tan^2 w/2
\sum_{j=k+1}^N p_j\tilde\chi\left(q_j-\eta/2\right)
\\
\displaystyle\qquad+2 \zeta\cot w/2
\left[
\cosh\left(\frac \zeta 2 \cot^2 w/2\right)z \cot w/2
\sum_{j=k+1}^N p_j\chi\left(q_j-\eta/2\right)\right.
\\
\qquad\qquad\qquad\qquad\qquad\qquad\qquad \dis
\left.-
\sinh\left(\frac\zeta 2\cot^2 w/2\right)
\sum_{j=k+1}^N p_j\tilde\chi\left(q_j-\eta/2\right)
\right]
\\
\displaystyle\dot p_i =p_i
\left[
\cosh\left(\frac \zeta 2 \cot^2 w/2\right)z
\chi\left(q_i-\eta/2\right)
+
\sinh\left(\frac \zeta 2\cot^2 w/2\right)\tan w/2\tilde
\chi\left(q_i-\eta/2\right)
\right]
\\
\displaystyle\qquad+p_i\sum_{j=k+1}^N p_j \sign(q_i- q_j)\chi\left({|q_i-q_j|}\right)
\\
\displaystyle\dot q_i=
\cosh\left(\frac \zeta 2\cot^2 w/2\right)z
\chi\left(q_i- \eta/2 \right)
+
\sinh\left(\frac \zeta 2 \cot^2 w/2\right)\tan w/2\tilde \chi
\left(q_i- \eta/2\right)
\\
\qquad\dis
+\sum_{j=k+1}^N p_j\chi\left({|q_i-q_j|}\right)
\end{array}
$$
with initial data
$$z(\tau)= \bar p\,,\qquad w(\tau)=\pi\,,\qquad
\eta(\tau)= 2\bar q\,,\qquad \zeta(\tau)= e_{(\tau,\bar q)}\,,$$
$$p_i(\tau)=\lim_{t\to \tau-}p_i(t)\,,
\qquad\qquad q_i(\tau)=\lim_{t\to \tau-}q_i(t)\qquad\qquad i=k+1,\ldots,N\,.$$
For the above system of O.D.E's, a direct inspection
reveals that the right hand side can be
extended by continuity also at the value $w=\pi$, because
all singularities are removable.
This continuous extension is actually smooth,
in a neighborhood of the initial data. Therefore, our Cauchy
problem has a unique local solution.
This provides a multi-peakon solution defined on some interval of the
form $[\tau,\,\tau'[\,$, up to the next interaction time.

The case where two or more groups of peakons interact
exactly at the same time $\tau$, but at different locations within
the interval $[0,1]$,
can be treated in exactly
the same way.
Since the total number of peakons (on a unit interval in the
$x$-variable) does not increase, it is clear that the number of
interaction times is at most countable.  The solution can thus be extended
to all times $t>0$, conserving its total energy.
\end{proof}

\chapter{Distance defined by optimal transportation problem \label{chap3}}

\section{A distance functional in the spatially periodic case\label{3-3}}
In this section we shall
construct a functional $J(u,v)$ which controls the distance
between two solutions of the equation (\ref{ch1-equation}).
All functions and measures on $\R$ are assumed to be periodic with
period 1.
Let $\T$ be the unit circle, so that $\T=[0,2\pi]$ with the
endpoints $0$ and $2\pi$ identified.
The distance $|\theta-\theta'|_*$ between two points $\theta,\theta'\in\T$
is defined as the smaller between the lengths of the two arcs connecting
$\theta$ with $\theta'$ (one clockwise, the other counterclockwise).
We now consider the product space
$$X\doteq \R\times\R\times \T$$
with distance
\begin{equation}
d^\diamondsuit\Big( (x,u,w), ~(\tilde x, \tilde u,\tilde w)\Big)
\doteq \Big(|x-\tilde x|+|u-\tilde u|+ |w-\tilde w|_*\Big)\wedge 1 \,,
\label{ch1-defdistance}
\end{equation}
where $a\wedge b\doteq \min\{a,b\}$.
Let $\M(X)$ be the space of all Radon measures on $X$ which are 1-periodic
w.r.t.~the $x$-variable. To
each 1-periodic function $u\in H^1_\per $
we now associate the positive measure
$\sigma^u\in\M(X)$ defined as
\begin{equation}
\sigma^{u}(A)\doteq \int_{\big\{ x\in\R\,:~(x,\,u(x),\, 2\arctan
\,u_x(x)\,)\in A\big\}} \big(1+u_x^2(x)\big)\, dx
\end{equation}
for every Borel
set $A\subseteq\R^2\times \T\,$.
Notice that the total mass of $\sigma^{(u,\mu)}$ over one period is
$$\sigma^u\big([0,1]\times\R\times\T)=
1+\int_0^1 u_x^2(x)\,dx
\,.$$

On this family of positive, 1-periodic Radon measures, we now introduce
a kind of Kantorovich distance,
related to an optimal transportation problem.
Given the two measures $\sigma^u$ and
$\sigma^{\tilde u}$,
their distance $J(u,\tilde u)$ is defined as follows.
\v
Call $\F$ the family of all
strictly increasing absolutely continuous
maps $\psi:\R\mapsto\R$ which have an absolutely continuous
inverse and satisfy
\begin{equation}
\psi(x+n)=n+\psi(x)\qquad\qquad \hbox{for every}~~n\in \Z\,.
\label{ch1-psifunct}
\end{equation}
For a given $\psi\in\F$, we define the 1-periodic, measurable functions
$\phi_1,\phi_2:\R\mapsto [0,1]$ by setting
\begin{equation}
\begin{array}{l}
\phi_1(x)\doteq
\sup\,\bigg\{\theta\in [0,1]\,;~~ \theta\cdot
\Big( 1+u_x^2(x)\Big)\leq \Big(
1+ \tilde u_x^2
\big(\psi(x)\big)\Big)\,\psi'(x)\bigg\}\,,
\\
\phi_2(\psi(x))\doteq
\sup\,\bigg\{\theta\in [0,1]\,;~~
1+u_x^2(x)\geq \theta\cdot \Big(
1+ \tilde u_x^2
\big(\psi(x)\big)\Big)\,\psi'(x)\bigg\}\,.
\end{array}
\label{ch1-phidef}
\end{equation}
Observe that the above definitions imply
$\max\big\{ \phi_1(x),\,\phi_2(x)\big\}=1$ together with
\begin{equation}
\phi_1(x) \,\Big( 1+u_x^2(x)\Big)=\phi_2\big(\psi(x)\big)\,\Big(
1+\tilde u_x^2
\big(\psi(x)\big)\Big)\,\psi'(x)
\label{ch1-phiidentity}
\end{equation}
for a.e. $x\in\R$.
We now define
\begin{equation}
\begin{array}{rl}
\dis J^\psi(u,\tilde u)\doteq &
\!\!\!\!\!\!\dis \int_0^1 \!\!d^\diamondsuit\Big( \big(x,\,u(x), 2\arctan u_x(x)\big)\,,
\big(\psi(x),\,\tilde u(\psi(x)),2\arctan \tilde u_x(\psi(x))\Big) \cr
&\dis \qquad\qquad\qquad\qquad\qquad\qquad \cdot
\phi_1(x)\,\big(1+u_x^2(x)\big)\,dx\cr
&\dis +\int_0^1\Big| \big(1+u_x^2(x)\big)-\big(1+\tilde u_x^2(\psi(x))\big)
\,\psi'(x)\Big|\,dx\,.\cr
\end{array}
\label{ch1-defunzjei}
\end{equation}
Of course, the integral is always computed over one period.
Observe that the map $x\mapsto \psi(x)$ can be
regarded as a \emph{transportation plan}, in order to transport the
measure $\sigma^u$ onto the measure $\sigma^{\tilde u}$.
Since these two positive
measures need not have the same total mass, we allow the presence of
some excess mass, not transferred
from one place to the other. The penalty for this
excess mass is given by the second integral in (\ref{ch1-defunzjei}).
The factor $\phi_1\leq 1$ in the first integral indicates
the percentage of the mass which is actually transported.
Integrating (\ref{ch1-phiidentity})
over one period, we find
$$\int_0^1\phi_1(x) \big(1+u_x^2(x)\big)\,dx=
\int_0^1 \phi_2(y) \big(1+\tilde u_x^2(y)\big)\,dy\,.$$
We can thus transport the
measure $\phi_1 \,\sigma^u$ onto $\phi_2\,
\sigma^{\tilde u}$ by a map
$$
\Psi:\big(x,\, u(x)\, \arctan u_x(x)\big)\mapsto \big(y,
\,\tilde u(y),\,\arctan \tilde u_x(y)\big),$$ 
where $y=\psi(x)$.
The associated
cost is given by the first integral in (\ref{ch1-defunzjei}). Notice that in this
case the measure $\phi_2\, \sigma^{\tilde u}$
is obtained as the push-forward
of the measure $\phi_1 \,\sigma^u$. We recall that the \emph{push-forward} of a measure $\sigma$ by a mapping $\Psi$ is defined as
$(\Psi\#\sigma)(A)\doteq \sigma(\Psi^{-1}(A))$ for every measurable set
$A$.  Here $\Psi^{-1}(A)\doteq \big\{ z\,;~\Psi(z)\in A\big\}$.

\v Our distance functional $J$ is now obtained by optimizing over all
transportation plans, namely
\begin{equation}
J(u,\tilde u)\doteq \inf_{\psi\in\F} J^\psi
(u,\tilde u)\,.
\label{ch1-optimaljei}
\end{equation}
\v
To check
that (\ref{ch1-optimaljei}) actually defines a distance, let $u,v,w\in H^1(\R)$ be
given.
\v
\n{\bf 1.} Choosing $\psi(x)=x$, so that
$\phi_1(x)=\phi_2(x)=1$, we immediately see that $J(u,u)=0$. Moreover,
if $J(u,\tilde u)=0$, then by the definition of $d^\diamondsuit$
we have $\tilde u= u$.
\v
\n{\bf 2.} Given $\psi\in\F$,
define $\tilde\psi=\psi^{-1}$, so that $\tilde\phi_1=\phi_2$,
$\tilde\phi_2=\phi_1$. This yields
$$J^{\tilde \psi}(\tilde u, u)=
 J^\psi(u,\tilde u)\,.$$
Hence $J(\tilde u,u)=J(u,\tilde u)$.
\v
\n{\bf 3.} Finally, to prove the
triangle inequality, let $\psi^\flat,\psi^\sharp:\R\mapsto\R$ be
two increasing diffeomorphisms satisfying (\ref{ch1-psifunct}), and let
$\phi_1^\flat,\phi_1^\sharp,\phi_2^\flat,\phi_2^\sharp:\R\mapsto
[0,1]$ be the corresponding functions, defined as in (\ref{ch1-phidef}).
We now consider the
composition $\psi\doteq \psi^\sharp\circ\psi^\flat$ and define
the functions $\phi_1,\phi_2$ according to (\ref{ch1-phidef}).
Observing that
$$\phi_1(x)\geq \phi_1^\flat(x)\cdot\phi_1^\sharp\big(\psi^\flat(x)\big)
\,,$$
$$\phi_2\big(\psi(x)\big)=\phi_2\Big(\psi^\sharp\big(\psi^\flat(x)\big)\Big)
\geq
\phi_2^\flat\big(\psi^\flat(x)\big)\cdot\phi_2^\sharp
\Big(\psi^\sharp\big(\psi^\flat(x)\big)\Big)
\,,$$
and recalling that the distance $d^\diamondsuit$ at (\ref{ch1-defdistance}) is always $\leq 1$,
we conclude
$$J^\psi(u,w)\leq
J^{\psi^\flat}(u,v)+J^{\psi^\sharp}(v,w)\,.
$$
This implies the triangle inequality $J(u,v)+J(v,w)\geq
J(u,w)$.
\endproof

In the remainder of this section we study the relations between our
distance functional $J$ and the distances determined by various norms.
\begin{lemma}
\label{lem3} 
For any $u,v\in H^1_\per$ one has
$$ \frac 1C\cdot
\|u-v\|_{L^1_\per}\leq J(u,v)\leq C\cdot\|u-v\|_{H^1_\per}\,,
$$
with a constant $C$ uniformly valid on bounded subsets of $H^1_\per$.
\end{lemma}
\begin{proof}
We shall use the elementary bound
$$
\big| \arctan a-\arctan b\big|\cdot a^2\leq
4\pi\big(|a|+|b|\big)\,|a-b|\,,
$$
valid for all $a,b\in\R$.
In connection with the identity mapping $\psi(x)=x$
we now compute
$$
\begin{array}{rl}
\dis
 J^\psi(u,v)&
\dis \leq \int_0^1 \Big\{
\big| u(x)-v(x)\big|+2 \big|\arctan u_x-\arctan v_x\big|
\Big\}\, (1+u_x^2)\,dx
\\
&\qquad\dis
+\int_0^1 \big|u_x^2-v_x^2|\,dx\cr
&\dis \leq \|u-v\|_{L^\infty}\,\|1+u_x^2\|_{L^1}+ (8\pi+1)\int_0^1
|u_x+v_x|\, |u_x-v_x|\,dx
\\
&\dis \leq (8\pi + 3) \,\big(1+ \|u\|_{H^1}+\|v\|_{H^1}\big)\cdot
\|u-v\|_{H^1}\,,
\end{array}
$$
proving the second inequality in (3.8).
\psfrag{x+dx}{$x+dx$}
\psfrag{u}{$u$}
\psfrag{v}{$v$}
\psfrag{dA}{$dA$}
\psfrag{x}{$x$}
\psfrag{psi}{$\psi(x)$}
\psfrag{psidelta}{$\psi(x+dx)$}
\begin{figure}
\centerline{\includegraphics[width=10cm]{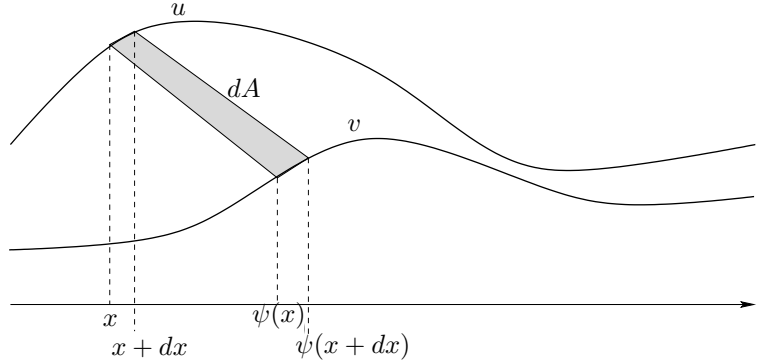}}
\caption{Infinitesimal area between $u$ and $v$ \label{infinitesimal}}
\end{figure}

To achieve the first inequality, choose any $\psi\in\F$.
For $x\in [0,1]$, call $\gamma^x$ the segment joining the
point $P^x=\big(x,u(x)\big)$ with $Q^x=\big(\psi(x), v(\psi(x))\big)$.
Clearly, the union of all these segments covers the
region between the graphs of $u$ and $v$.
Moving the base point from $x$ to $x+dx$, the corresponding segments
sweep an infinitesimal area $dA$ estimated by (fig.\ref{infinitesimal})
$$\begin{array}{rl}
|dA|&
\dis
\leq |P^x-Q^x|\cdot \big(|dP^x|+|dQ^x|\big)
\\
&
\dis
\leq
\Big(\big|x-\psi(x)\big|^2+\big|u(x)-v(\psi(x))\big|^2\Big)^{1/2}
\!\!\cdot
\!\!\Big[(1+u_x^2)^{1/2}+(1+v_x^2)^{1/2}\psi'(x)\Big]\,dx\,.
\end{array}
$$

Integrating over one period we obtain
$$
\begin{array}{rl}
\dis
\int_0^1\big|u(x)-v(x)\big|\,dx
&\leq
\dis
\int_0^1\Big(\big|x-\psi_{(t)}(x)\big|+
\big|u(x)-v(\psi_{(t)}(x))\big|\Big)
\cdot
\\
&\qquad\qquad\dis\cdot
\Big[\big(1+u_x^2(x)\big)^{1/2}+\big(1+v_x^2(\psi(x))\big)^{1/2}
\psi'_{(t)}(x)\Big]\,dx
\\
&
\dis
\leq \big(2+\|u\|_{H^1}+\|v\|_{H^1}\big)\cdot
\Big[J^\psi(u,v)+J^{\psi^{-1}}
(v,u)\Big]
\\
&\leq C \cdot J(u,v)\,,
\end{array}
$$
completing the proof of (3.8).
\end{proof}
\begin{lemma}
\label{lem4}
Let $(u_n)_{n\geq 1}$ be a
Cauchy sequence
for the distance $J$, uniformly
bounded in the $H^1_\per$ norm.
Then
\v
\begin{itemize} 
\item[(i)] There exists a limit function
$u\in H^1_\per$ such that $u_n\to u$ in $L^\infty$
and the sequence of derivatives
$u_{n,x}$ converges to $u_x$ in $L^p_\per$,
for $1\leq p <2$.
\item[(ii)] Let $\mu_n$ be the absolutely continuous measure
having density $u^2_{n,x}$ with respect to Lebesgue measure.
Then one has the weak convergence $\mu_n\wto \mu$, for
some measure $\mu$ whose absolutely continuous part
has density $u_x^2$.
\end{itemize}
\end{lemma}
\begin{proof}{\bf 1.}
By Lemma \ref{lem3} we already know the convergence
$u_n\to u$, for some limit function $u\in
L^1_\per\,$.  By a Sobolev embedding theorem,
all functions $u_n,u$  are uniformly H\"older continuous.
This implies $\|u_n-u\|_{L^\infty}\to 0$.
To establish the convergence of derivatives, we first show that
the sequence of functions
$$
v_n\doteq \exp\{ 2i\arctan u_{n,x}\}
$$
is compact in $L^1_\per\,$.

Indeed, fix $\ve>0$.  Then there exists $N$ such that
$J(u_m,u_n)<\ve$ for $m,n\geq N$.
We can now approximate $u_N$ in $H^1_\per$ with a
piecewise affine function
$\tilde u_N$ such that $J(\tilde u_N, u_N)\leq \ve$.
By assumption, choosing suitable transport maps $\psi_n$
we obtain
$$
\begin{array}{rl}
\dis
\int_0^1 \Big|\exp \big\{ 2i\arctan u_{n,x}(x)\big\}-
\exp\big\{ 2i\arctan \tilde u_{N,x}(\psi_n(x))\big\}\Big|\,dx
& \dis \!\!\!\leq~ 2\,J(u_n,\tilde u_N)		
\\ &\!\!\!\leq ~ 4 \ve
\end{array}
$$
for all $n\geq N$.
We now observe that all functions
$$
x\mapsto \exp\big\{ 2i\arctan \tilde u_{N,x}(\psi_n(x))\big\}
$$
are uniformly bounded, piecewise constant with the same
number of jumps: namely, the number of subintervals
on which $\tilde u_N$ is affine.  The set of all
such functions is compact in $L^1_\per$.
This argument shows that
the sequence $v_n\doteq \exp\{ 2i\arctan u_{n,x}\}$
eventually remains in an $\ve$-neighborhood of a compact
subset of $L^1_\per$.  Since $\ve>0$ can be taken arbitrarily
small, by possibly choosing a subsequence we
obtain the strong convergence
$v_n\to v$ for some  $v\in L^1_\per\,$.
\v
\n{\bf 2.} From the uniform $H^1$ bounds and the $L^1$ convergence
of the functions $v_n$, we now derive the $L^p$ convergence of the
derivatives.  For a given $\ve>0$, define
$$M\doteq \sup_n\|u_n\|_{H^1_\per}\,,
\qquad\qquad
A_n\doteq \Big\{ x\in [0,1]\,;~~ \big|u_{n,x}(x)\big| > M/\ve\Big\}\,.$$
The above definitions imply
$$
\meas (A_n)\leq \ve^2
$$
We now have
$$\begin{array}{rl}
\dis
\|u_{m,x}-u_{n,x}\|_{L^p}&
\dis \leq \left(\int_{A_n\cup A_m}
|u_{m,x}-u_{n,x}|^p\,dx\right)^{1/p}
\\
&\qquad+
\dis
\left(\int_{[0,1]\setminus(A_n\cup
A_m)} |u_{m,x}-u_{n,x}|^p\,dx\right)^{1/p}
\\
&
\dis
\doteq I_1+I_2\,.
\end{array}
$$

\begin{equation}
\begin{array}{rl}
\dis I_1 &\dis
\leq\left[\int_{A_m\cup A_n}
1\cdot dx\right]^{(2-p)/2p} \cdot \left[\int_{A_m\cup
A_n}\big(
|u_{m,x}|+|u_{n,x}|\big)^2\,dx\right]^{1/2}
\\
&\leq\dis  \ve^{(2-p)/p}\cdot 2M
\,.
\end{array}
\label{ch1-smallenough}
\end{equation}
Next, choosing a constant $C_\ve$ such that
\begin{equation}
\big|e^{2i\arctan a}-e^{2i\arctan b}\big|\geq C_\ve
|a-b|\qquad\qquad\hbox{whenever}~~|a|, |b|\leq
M/\ve\,,$$  we obtain
$$I_2 \leq C_\ve \,\left[\int  \Big|e^{2i\arctan u_{m,x}}
-e^{2i\arctan u_{n,x}}\big|^p\,dx\right]^{1/p}.
\label{ch1-i2piccolo}
\end{equation}
Taking $\ve>0$ small,
we can make the right hand side of (\ref{ch1-smallenough}) as small as we like.
On the other hand, choosing a subsequence such that
$v_\nu=e^{2i\arctan u_{\nu,x}}$ converges in $L^1_\per$,
the right hand side of (\ref{ch1-i2piccolo}) approaches zero.
Hence, for this subsequence,
$$\limsup_{m,n\to\infty}\|u_{m,x}-u_{n,x}\|_{L^p_\per}=0\,.$$
Since $u_n\to u$ uniformly,
in this case we must have
\begin{equation}
\|u_{n,x}-u_x\|_{L^p_\per}\to 0\,.
\label{ch1-gradconv}
\end{equation}
We now observe that from any subsequence we can extract a further
subsequence for which (\ref{ch1-gradconv}) holds.
Therefore, the whole sequence $(u_{n,x})_{n\geq 1}$
converges to $u_x$ in $L^p_\per\,$.
\v
\n{\bf 3.}   To establish (ii), we consider
the sequence of measures having density $1+u_{n,x}^2$
w.r.t.~Lebesgue measure.
This sequence converges weakly, because our distance functional is
stronger than the Kantorovich-Waserstein metric which induces
the topology of weak convergence on spaces of measures.
Therefore, $\mu_n\wto \mu$ for some positive measure $\mu$.

Since the sequence $1+u_{n,x}$ converges to $1+u_x$
in $L^1_\per$, by possibly choosing a subsequence we
achieve the pointwise convergence $u_{n,x}(x)\to u_x(x)$,
for a.e.~$x\in [0,1]$.  For any $\ve>0$, by Egorov's theorem
we have the uniform convergence $u_{n,x}(x)\to u_x(x)$
for all $x\in [0,1]\setminus V_\ve$, for some set
with $\meas(V_\ve)<\ve$.  Since $\ve>0$ can
be taken arbitrarily small, this shows that the absolutely
continuous part of the measure $\mu$ has density $u^2+u_x^2$
w.r.t.~Lebesgue measure.
\end{proof}

\subsection{Continuity in time of the distance functional\label{3-4}}
Here and in the next section we examine how the
distance functional $J(\cdot,\cdot)$ evolves in time,
in connection with multi-peakon solutions of
the Camassa-Holm equation (\ref{ch1-equation}).
We first provide estimates valid on a time interval where no
peakon interactions occur.  Then we show that the
distance functional is continuous across times of interaction.
Since the number of peakons is locally finite, this will suffice
to derive the basic estimates (\ref{ch1-contu})-(\ref{ch1-stabsol}), in the case of
multi-peakon solutions.
\begin{lemma}
\label{lem5}
 Let $t\mapsto u(t)\in H^1_\per$ be a multi-peakon
solution of (\ref{ch1-equation}).   Assume that no peakon interactions occur
within the interval $[0,\tau]$.  Then
\begin{equation}
J\big(u(s),\, u(s')\big)\leq   C\cdot |s-s'|\,,\qquad\qquad
\forall s,s'\in [0,\tau]\,,
\label{ch1-continuity}
\end{equation}
for some constant $C$, uniformly valid as $u$ ranges on bounded subsets
of $H^1_\per\,$.
\end{lemma}
\v
\begin{proof}  Assume $0\leq s<s'\leq \tau$.
By the assumptions, the solution $u=u(t,x)$
remains uniformly Lipschitz continuous on the time interval $[0,\tau]$.
Therefore, for each $s\in [0,\tau]$ and $x\in\R$, the Cauchy problem
\begin{equation}
\frac d {dt}\,\xi(t)= u\big( t,\,\xi(t)\big)\,,\qquad
\qquad \xi(s)=x\,,
\label{ch1-characteristics}
\end{equation}
determines a unique
characteristic curve $t\mapsto \xi(t; s,x)$ passing through the point $(s, x)$.
Given $s'\in [0,\tau]$, we can thus define a transportation plan by setting
\begin{equation}
\psi(x)\doteq \xi(s';s,x)\,.
\label{ch1-tranportchar}
\end{equation}
Of course, moving mass along the characteristics
is the most natural thing to do.
We then choose $\phi_1,\phi_2$ to be as large as possible,
according to (\ref{ch1-phidef}).  Namely:
$$\phi_1(x)\doteq
\sup\,\bigg\{\theta\in [0,1]\,;~~ \theta\cdot
\Big( 1+u_x^2(s,x)\Big)\leq \psi'(x)\cdot\Big(
1+ u_x^2
\big(s',\,\psi(x)\big)\Big)\bigg\}\,,$$
$$\phi_2(x)\doteq
\sup\,\bigg\{\theta\in [0,1]\,;~~ \,\theta\cdot
\Big(
1+u_x^2
\big(s',\,\psi(x)\big)\Big)\,\psi'(x)\leq
1+u_x^2(s,x)\bigg\}\,.$$
The cost of this plan is
bounded by
\begin{equation}
\begin{array}{l}
\dis
J^\psi\big( u(s)\,,~u(s')\big)\leq
\int_0^1\Big\{ \big|x-\xi(s';s,x)\big|+ \big|
u(s,x)-u(s',\xi(s';s,x))\big|
\\
\qquad\qquad+
\big|2\arctan u_x(s,x)-2\arctan  u_x(s',\xi(s';s,x))
\big|_*\Big\}\big(1+u_x^2(s,x)\big)\,dx
\\
\dis
\qquad\quad +\int_0^1 \big(1-\phi_1(x)\big)\,\big(1+u^2_x(s,x)\big)\,dx\\
\qquad\quad 
\dis +\int_0^1 \big(1-\phi_2(\psi(x))\big)\,\Big(1
+u^2_x\big(s',\psi(x)\big)\Big)\,\psi'(x)\,dx
\,.
\end{array}
\label{ch1-s1}
\end{equation}
To estimate the right hand side of (\ref{ch1-s1}), we first
observe that, for all $u\in H^1_\per$,
\begin{equation}
\begin{array}{rl}
\dis
\big\|u\big\|_{L^\infty}
&
\dis
\leq \int_0^1 \big|u(x)\big|\,dx+\int_0^1
\big|u_x(x)\big|\,dx\\
&
\dis
\leq \|u\|_{L^2}+\|u_x\|_{L^2}\leq 2\|u\|_{H^1_\per}=2(E^u)^{1/2}
.
\end{array}
\label{ch1-s2}
\end{equation}
Using (\ref{ch1-s2}) in (\ref{ch1-characteristics}) we obtain
\begin{equation}
\big|\xi(s)-\xi(s')\big|\leq 2(E^{\bar u})^{1/2}\cdot |s-s'|\,.
\label{ch1-lipchar}
\end{equation}

Next, from the definition of the source term $P$ at (\ref{ch1-nonlocalP})
it follows
\begin{equation}
\big\|P\|_{L^\infty}\leq \frac 12 \big\|e^{-|x|}\,\big\|_{L^\infty(\R)}
\cdot \left\| u^2+\frac{u_x^2}2\right\|_{L^1([0,1])}\leq
\|u\|^2_{H^1_\per}
=E^u\,.
\label{ch1-stimaperP}
\end{equation}
Similarly,
\begin{equation}
\big\|P_x\|_{L^\infty}\leq \|u\|^2_{H^1_\per}
=E^u\,.
\label{ch1-stimaperPx}
\end{equation}
Using (\ref{ch1-stimaperPx}) we obtain
\begin{equation}
\begin{array}{l}
\dis
\Big| u\big(s',\,\xi(s')\big)-u\big(s,\,\xi(s)\big)\Big|
\leq \int_s^{s'}\left| \frac d{dt} u\big(t,\,\xi(t)\big)\right|\,dt\\
\dis
\qquad\qquad=\int_s^{s'}\Big| P_x\big(t,\,\xi(t)\big)\Big|\,dt~\leq
~E^{\bar u} \cdot |s'-s|\,. 
\end{array}
\label{ch1-lipu}
\end{equation}

Concerning the term involving arctangents,
recalling (\ref{ch1-nonlocalequation}) we obtain
$$\frac d{dt}\Big[2\arctan u_x\big(t,\xi(t,x)\big)\Big]=\frac 2{1+u_x^2}\,
\left[ u^2-\frac{u_x^2}2-P\right]\,.$$
The bounds (\ref{ch1-stimaperP}) and (\ref{ch1-lipu}) thus yield
\begin{equation}
\begin{array}{l}
\dis
\Big|2\arctan u_x\big(s',\xi(s')\big)
-2\arctan u_x\big(s,\,\xi(s)\big)\Big|_*
\\
\qquad\qquad\qquad\leq
\dis
 \Big( 2\|u\|^2_{L^\infty} +1+2\|P\|_{L^\infty}\Big)\cdot |s'-s|
\\
\qquad \qquad\qquad\leq
\dis
\big(10\,E^{\bar u}+1\big)\cdot |s'-s|\,.
\end{array}
\label{ch1-liparctan}
\end{equation}
This already provides a bound on
the first integral on the right hand side of (4.4).

Next, call $I_1,I_2$ the last two integrals on the right hand side of (4.4).
Notice that
$$
I_1+I_2=\int_0^1 \bigg|\Big(1+ u_x^2(s,y)\Big)-\xi_y(s';s,y)\,
\Big( 1+ u_x^2\big(s',\xi(s';s,y)\big)\Big)\bigg|\,dy
$$
Indeed, $I_1+I_2$ measures the difference between
the measure $\big( 1+u_x^2(s',y)\big)\,dy$ and the push-forward
of the measure $\big( 1+u_x^2(s,x))\,dx$ through the mapping
$x\mapsto \xi(s';s,x)$.

Since the push-forward of the measure $u_x^2\,dy$ satisfies the
linear conservation law
\begin{equation}
w_t+(uw)_x=0\,,
\label{ch1-transpot}
\end{equation}
comparing (\ref{ch1-transpot}) with (\ref{ch1-uxquadro}) we deduce
$$
\begin{array}{l}
\dis
\int_0^1 \Big| u_x^2(s,y)-\xi_y(s,y)\,
 u_x^2\big(s',\xi(s';s,y)\big)\Big|\,dy~\leq
~\int_s^{s'} \int_0^1 2\big| (u^2-P)u_x\big|\,dx\,dt\\
\qquad
 \leq
\dis
2 \int_s^{s'}  \big( \|u\|_{L^\infty}^2+
\|P\|_{L^\infty}\big)\,\|u_x\|_{L^1}\,dt
~\leq~ 2\,\big(4E^{\bar u}+E^{\bar u}\big)\,E^{\bar u}\cdot |s'-s|\,,
\end{array}
$$
because of (\ref{ch1-s2}), (\ref{ch1-stimaperP}) and (\ref{ch1-energy}).
Finally, we need to estimate the remaining terms,
describing by how much the Lebesgue measure fails to be conserved
by the transformation  $x\mapsto \xi(s';s,x)$.
Observing that
$$
\frac \partial {\partial t}\,\xi_y(t,y)= u_x\big(t,\,\xi(t,y)\big)\,
\xi_y(t,y)\,,\qquad\qquad\xi_y(0,y)= 1\,,
$$
we find
\begin{equation}
\begin{array}{l}
\dis
\int_0^1
\big|1 -\xi_y(s';s,y)\,
\big|\,dy\leq
\int_s^{s'} \int_0^1\bigg|\frac\partial {\partial t}\big[
\xi_y(t,s,y)\big]\bigg|
\,dy\,dt
\\
\dis
\qquad\leq \int_s^{s'}\int_0^1
\xi_y(t;s,y)\,
\big|u_x(t,\xi(t;s,y))\big|
\,dy\,dt
\,.\end{array}
\label{ch1-s3}
\end{equation}
To estimate the right hand side of (\ref{ch1-s3}),
we use the decomposition
$[0,1]=Y\cup Y'\cup Y''$, where
$$Y\doteq \Big\{ y\,;~~ \xi_y(t;s,y)\in \big[(1/2)\,,~2\big]\quad
\hbox{for all}~t\in [s,s']\,\Big\}\,,$$
$$Y'\doteq \Big\{ y\,;~~ \xi_y(t;s,y) <1/2 \quad
\hbox{for some}~t\in [s,s']\,\Big\}\,,$$
$$Y''\doteq \Big\{ y\,;~~ \xi_y(t;s,y) > 2\quad
\hbox{for some}~t\in [s,s']\,\Big\}\,.$$
Integrating over $Y$ one finds
$$
\int_s^{s'} \int_Y
\xi_y(t;s,y)
\,\big|u_x(t,\xi(t;s,y))\big|
\,dy\,dt~\leq~ 2\int_s^{s'}
\big\|u_x(t)\big\|_{L^1}\,dt~\leq~ 2 E^{\bar u}\cdot |s'-s|\,.
$$
Next, if $y\in Y'$ we define
$$\tau(y)\doteq \inf\,\big\{ t>s\,;~~\xi_y(t;s,y)< 1/2\big\}\,.$$
Observe that $y\in Y'$ implies
$$\int_s^{\tau(y)} \Big|u_x\big(t,\,\xi(t;s,y)\big)\Big|\,dt
\geq \ln 2\,.$$
Therefore
$$
\begin{array}{rcl}
\dis
\int_{Y'} dy &\leq&
\dis
\frac 1{\ln 2}\int_{Y'}
\left[\int_s^{\tau(y)}
\Big|u_x\big(t,\,\xi(t;s,y)\big)\Big|\,dt
\right]\,dy\\
&\leq&
\dis
\frac 2{\ln 2}\int_s^{s'}\int_0^1
\big|u_x\big(t,x)\big|\,dx\,dt~\leq~ 4 \,E^{\bar u} \cdot |s'-s|\,.
\end{array}
$$
The estimate for the integral over $Y''$ is entirely analogous,
Indeed, the push-forward of the Lebesgue measure along characteristic
curves from $t=s$ to
$t=s'$ satisfies exactly the same type of estimates as
the pull-back of the Lebesgue measure from $t=s'$ to
$t=s$. All together, these three estimates imply
\begin{equation}
\int_{Y\cup Y'\cup Y''}
\big|1 -\xi_y(s';s,y)\,
\big|\,dy\leq 10\,E^{\bar u}\cdot |s'-s|
\label{ch1-s4}
\end{equation}
\v
Putting together the estimates (\ref{ch1-lipchar}), (\ref{ch1-lipu}), (\ref{ch1-liparctan}), (\ref{ch1-tranportchar}) and (\ref{ch1-s4}),
the distance in (4.4) can be estimated by
$$
J^\psi\big( u(s)\,,~u(s')\big)\leq\Big[
2(1+E^{\bar u}) +E^{\bar u} + (10 \,E^{\bar u} +1)+
10\,(E^{\bar u})^2+10\,E^{\bar u}\Big]\cdot |s'-s|\,.
$$
This establishes (\ref{ch1-continuity}).
\end{proof}
\v
According to Lemma \ref{lem5}, as long as no peakon interactions occur,
the map $t\mapsto u(t)$ remains
uniformly Lipschitz continuous w.r.t.~our distance functional,
with a Lipschitz constant that depends only
on the total energy $E^{\bar u}$.  Since interactions
can occur only at isolated times, to obtain a global Lipschitz estimate
it suffices to show that trajectories are continuous
(w.r.t.~the distance $J$) also at interaction times.
\begin{lemma}
\label{lem6}
Assume that the multi-peakon solution
$u(\cdot)$ contains two or more peakons which interact at
a time $\tau$.   Then
$$
\lim_{h\to 0+} J\big( u(\tau -h),\, u(\tau + h)\big)=0\,.
$$
\end{lemma}
\begin{proof} To fix the ideas, call $x=\bar q$ the place
where the interaction occurs,
and let $p_1,\ldots, p_k$ be the strengths of the peakons
that interact at time
$\tau$.  We here assume that $0<\bar q<1$.
The case where two or more groups of peakons interact
exactly at the same time $\tau$, within the interval
$[0,1]$, can be treated similarly.

For $|t-\tau|\leq h$, call $\xi^-(t)$, $\xi^+(t)$ respectively
the position of the smallest and largest characteristic
curves passing through the point $(\tau,\bar q)$, as in (\ref{ch1-charminmax}).
We observe that $u$ is Lipschitz continuous
in a neighborhood of each point $(\tau,x)$, with $x\not= \bar q$.
Hence, for $x\in [0,1]\setminus\{\bar q\}$ there exists a unique
characteristic curve $t\mapsto \xi(t; \tau,x)$ passing through
$x$ at time $\tau$.
For a fixed $h>0$, the transport map $\psi$ is defined as follows.
Consider the intervals
$I_{-h}\doteq \big[\xi^-(\tau-h),\,\xi^+(\tau-h)\big]$
and $I_h\doteq \big[\xi^-(\tau+h),\,\xi^+(\tau+h)\big]$.
On the complement $[0,1]\setminus I_{-h}$
we define
$$\psi\big(\xi(\tau-h;\tau ,x)\big)=\xi(\tau+h;\tau ,x)\,,
$$
so that transport is performed along characteristic curves.
It now remains to extends $\psi$ as a map from $I_{-h}$ onto $I_h$.
Toward this goal, we recall that
our construction of multi-peakon solutions in Section \ref{3-2}
was specifically designed in order to achieve the identity
\begin{equation}
e_{(\tau,\bar q)}=\lim_{h\to 0+} \int_{I_{-h}} u_x^2(\tau-h,x)\,dx~=~
\lim_{h\to 0+} \int_{I_h} u_x^2(\tau+h,x)\,dx\,.
\label{ch1-concentrenergy}
\end{equation}
For $h>0$ we introduce the quantities
$$E(-h)\doteq \int_{I_{-h}} \big(1+u_x^2(\tau-h,\,x)\big)\,dx\,,\qquad\qquad
E(h)\doteq \int_{I_h} \big(1+u_x^2(\tau+h,\,x)\big)\,dx\,,$$
$$e(h)\doteq 2\min\big\{E(-h),\,E(h)\big\}-\max\big\{E(-h),\,E(h)\big\}\,.$$
Notice that (\ref{ch1-concentrenergy}) implies $e(h)>0$ and
\begin{equation}
 E(-h)\to e_{(\tau,\bar q)}\,,
\qquad E(h)\to e_{(\tau,\bar q)}\,,
\qquad e(h)\to e_{(\tau,\bar q)}\,,
\label{ch1-enlimits}
\end{equation}
as $h\to 0+$.
Consider the point $x^*=x^*(h)$ inside the
interval 
$$
I_{-h}=\big[\xi^-(\tau-h),\,\xi^+(\tau-h)\big],
$$
implicitly defined by
$$\int_{\xi^-(\tau-h)}^{x^*}\big(1+u_x^2(x)\big)\,dx = e(h)\,.$$
For $x\in \big[\xi^-(\tau-h),~x^*\big]$ we define $\psi(x)$ as the
unique point such that
\begin{equation}
\int_{\xi^-(\tau+h)}^{\psi(x)}\big(1+u_x^2(\tau+h,\,x)\big)\,dx
=\int_{\xi^-(\tau-h)}^x\big(1+u_x^2(\tau-h,\,x)\big)\,dx\,.
\label{ch1-defmappsi}
\end{equation}
We then extend $\psi$ as an affine map from
$\big[x^*,\,\xi^+(\tau-h)\big]$ onto $\big[\psi(x^*),\,\xi^+(\tau+h)\big]$,
namely
$$\psi\Big(\theta\cdot\xi^+(\tau-h)+(1-\theta)
\cdot x^*\Big)=\theta\cdot\xi^+(\tau+h)+(1-\theta)
\cdot \psi(x^*)\qquad\qquad \theta\in [0,1]\,.$$
Finally, we prolong $\psi$ to the whole real line according to
(\ref{ch1-psifunct}).

As usual, the 1-periodic
functions $\phi_1,\phi_2$ are then chosen to be as large as possible,
according to (\ref{ch1-phidef}).
As $h\to 0+$, we claim that the following quantity approaches zero:
\begin{equation}
\begin{array}{l}
\dis 
J^\psi\big(u(\tau-h),~u(\tau+h)\big)
\\
\dis =
\int_0^1 d^\diamondsuit\Big( \big(x,\,u(\tau-h,x),\,
2 \arctan u_x(\tau-h,x)\big)\,,~
\\
\qquad\qquad\qquad\qquad \dis
\big(\psi(x),\,u(\tau+h, \,\psi(x)),\,2 \arctan \tilde u_x(\tau+h,\,
\psi(x))\big)\Big)\cdot \\
\dis \qquad\qquad\qquad \cdot
\phi_1(x)\,\big(1+u_x^2(\tau-h,\,x)\big)\,dx\\
\dis \qquad+\int_0^1 \big(1-\phi_1(x)\big)\,
\big( 1+u_x^2(\tau-h,\,x)\big)\,dx
\\
\dis \qquad +\int_0^1
\big(1-\phi_2(\psi(x))\big)\,\Big(1+u_x^2
\big(\tau+h,\,\psi(x)\big)\Big)\,\psi'(x)\,dx
\,.
\end{array}
\label{ch1-limsujey}
\end{equation}
It is clear that the restriction of all the above
integrals to the complement
$[0,1]\setminus I_{-h}$ approaches zero as $h\to 0$.
We now prove that their restriction to $I_{-h}$ also vanishes
in the limit.
As $h\to 0+$, for $x\in I_{-h}$ we have
$$
\begin{array}{l}
d^\diamondsuit\Big(
\big(x,\,u(\tau-h,x),\,2\arctan u_x(\tau-h,x)\big)~,~
\\
\dis \qquad \qquad\qquad
\big(\psi(x),\,u(\tau+h,\psi(x)),\,2\arctan u_x(\tau+h,\psi(x))\big)\Big)
~\to~ 0\,,
\end{array}$$
because all points approach the same limit
$\big(\bar q, \,u(\tau,\bar q),\, \pi\big)$.
The first integral in (\ref{ch1-limsujey}) thus approaches zero as $h\to 0+$.

Concerning the last two integrals, by (\ref{ch1-defmappsi}) it follows
$$\phi_1(x)=\phi_2\big(\psi(x)\big)=1\qquad\qquad \forall x\in
\big[\xi^-(\tau-h),\,x^*\big]\,.$$
Moreover, our choice of $x^*$ implies
\begin{equation}
\begin{array}{l}
\dis
\int_{x^*}^{\xi^+(\tau-h)}
\big( 1+u_x^2(\tau-h,\,x)\big)\,dx+\int_{\psi(x^*)}^{\xi^+(\tau+h)}
\Big(1+u_x^2
\big(\tau+h,\,\psi(x)\big)\Big)\,\psi'(x)\,dx
\\
\dis \qquad \leq
2\max\big\{E(-h),\,E(h)\big\}-2\min\big\{E(-h),\,E(h)\big\}\,.
\end{array}
\label{ch1-continuitaintau}
\end{equation}
By (\ref{ch1-enlimits}), as $h\to 0+$ the
right hand side of (\ref{ch1-continuitaintau}) approaches zero.
Hence the same holds for the last two integrals in (\ref{ch1-limsujey}).
This completes the proof of the lemma.
\end{proof}
\vsk

\subsection{Continuity w.r.t.~the initial data\label{3-5}}
We now consider two distinct
solutions and study how the distance $J\big(u(t)\,~
v(t)\big)$ evolves in time.
To fix the ideas, let $t\mapsto u(t)$ and $t\mapsto v(t)$
be two multi-peakon solutions of (\ref{ch1-equation}), and assume that
no interaction occurs within a given time interval $[0,T]$.
In this case, the functions $u,v$ remain Lipschitz continuous.
We can thus define the characteristic curves
$t\mapsto \xi(t,y)$ and $t\mapsto \zeta(t,\tilde y)$
as the solutions to the Cauchy problems
$$
\begin{array}{rcl}
\dis \dot \xi &=&
\dis u(t,\xi),\qquad\qquad \xi(0)=y\,,\\
\dis \dot \zeta &=&
\dis v(t,\zeta),\qquad\qquad \zeta(0)=\tilde y\,,
\end{array}
$$
respectively.
Let now $\psi_{(0)}\in\F$ be any transportation plan at time $t=0$.
For each $t\in [0,T]$ we can define a transportation plan
$\psi_{(t)}\in\F$ by setting
$$\psi_{(t)}\big(\xi(t,y)\big)\doteq \zeta\big(t,~\psi_{(0)}
(y)\big)
$$
The corresponding functions
$\phi_{1}^{(t)},\phi_2^{(t)}$
are then defined according to definitions in(\ref{ch1-phidef}), namely
$$\phi_1^{(t)}(x)\doteq
\sup\,\bigg\{\theta\in [0,1]\,;~~ \theta\cdot
\Big( 1+u_x^2(t,x)\Big)~\leq ~\Big(
1+ v_x^2
\big(t,\,\psi_{(t)}(x)\big)\Big)\,\psi_{(t)}'(x)\bigg\}\,,$$
$$\phi_2^{(t)}(x)\doteq
\sup\,\bigg\{\theta\in [0,1]\,;~~ 1+u_x^2(t,x)~\geq ~\theta\cdot
\Big(
1+v_x^2
\big(t,\,\psi_{(t)}(x)\big)\Big)\,\psi'_{(t)}(x)\bigg\}\,.$$
If
initially the point $y$ is mapped into $\tilde y\doteq \psi_{(0)}(y)$,
then at a later time $t>0$ the point $\xi(t,y)$ along the
$u$-characteristic starting from $y$ is sent to the point  $\zeta
(t,\tilde y)$ along the $v$-characteristic
starting from $\tilde y$. We thus transport
mass from the
point $\Big(\xi(t,y)\,,~u\big(t,\xi(t,y)\big)\,,~ 2\arctan
u_x\big(t,\xi(t,y)\big)\Big)$ to the corresponding point
$\Big(\zeta(t,\tilde y)\,, ~v\big(t,\zeta
(t,\tilde y)\big)\,, ~2\arctan v_x\big(t,\zeta
(t,\tilde y)\big)\Big)\,,$

In the following, our main goal is to provide an upper bound
on the
time derivative of the function
\begin{equation}
\begin{array}{l}
\dis J^{\psi(t)}\big(u(t)\,,~ v(t)\big)
\doteq
\int_0^1 d^\diamondsuit\Big( \big(x,\,u(t,x),\,
2 \arctan u_x(t,x)\big)\,,\\
\quad\qquad\qquad \dis \big(\psi_{(t)}(x),\,v(t,\,\psi_{(t)}(x)),\,
2 \arctan v_x(t,\,\psi_{(t)}(x))\big)\Big) \cdot
\phi_1^{(t)}(x)\,\big(1+u_x^2(t,x)\big)\,dx
\\
\dis
+\int_0^1 \big(1-\phi_1^{(t)}(x)\big)\,
\big( 1+u_x^2(t,x)\big)\,dx
\\ \dis
 +\int_0^1
\big(1-\phi_2^{(t)}(\psi_{(t)}(x))\big)\,\Big(1+v_x^2
\big(t,\,\psi_{(t)}(x)\big)\Big)\,\psi_{(t)}'(x)\,dx
\,.
\end{array}
\label{ch1-jeiuv}
\end{equation}
Differentiating the right hand side of (\ref{ch1-jeiuv})
one obtains several terms, due to
\begin{itemize}
\item[$\bullet$]Changes in the distance $d^\diamondsuit$
between the points
$(\xi,\, u,\, 2\arctan u_x)$ and $(\zeta,\, v, \,
2\arctan v_x)$.
\item[$\bullet$]Changes in the base measures $(1+u_x^2)\,dx$ and
$(1+v_x^2)\,dx$.
\end{itemize}
\n Throughout the following, by $\O(1)$ we denote a quantity which remains
uniformly bounded as $u,v$ range in bounded subsets of $H^1_\per\,$.
Using the elementary estimate
$$|u-v|\leq \big(1+|u|+|v|\big)\,\min\big\{ |u-v|,\,1\big\}\,,$$
we
begin by deriving the bound
$$\begin{array}{rcl}I_1
&\doteq&\dis \int_0^1 \frac d{dt}\,\big|  x-\psi_{(t)}(x)
\big|\cdot\phi_1^{(t)}(x)\,\big( 1+u_x^2(t,x)\big) \,\,dx
\\
&\leq&
\dis
 \int_0^1 \Big| u(t,x)-v\big(t,\,\psi_{(t)}(x)\big)
\Big|\cdot \phi_1^{(t)}(x)\,\big( 1+u_x^2(t,x)\big) \,\,dx
\\
&\leq& \Big(1+ \big\|u(t)\big\|_{L^\infty} +
\big\|v(t)\big\|_{L^\infty}
\Big)\cdot J^{\psi_{(t)}}\big(u(t),\,v(t)\big)\\
&=&
\dis \O(1)\cdot J^{\psi_{(t)}}\big(u(t),\,v(t)\big)\,.
\end{array}
$$
Here and in the sequel, the time derivative
is computed along characteristics.

Next,
recalling the basic equation (\ref{ch1-equation}), we consider
$$\begin{array}{rcl}
I_2 &\doteq&
\dis \int_0^1 \frac d{dt}\Big| u(t,x)- v(t,\,
\psi_{(t)}(x)
\big)\Big|\cdot\phi_1^{(t)}(x)\,\big( 1+u_x^2(t,x)\big) \,\,dx
\\
&\leq& \int_0^1 \Big| P_x^u(t,x)-P_x^v\big(t,\psi_{(t)}(x)\big)
\Big|\cdot\big( 1+u_x^2(t,x)\big) \,\,dx\,.
\end{array}
$$
In the spatially periodic case, by (\ref{ch1-nonlocalP}) and (\ref{ch1-formachi}) we
can write the source terms $P^u_x$, $P^v_x$ as
$$\begin{array}{rcl} 
\dis P^u_x(t,x)&=&
\dis \frac 12 \int_{x-1}^x \chi'(x-y)\cdot
\left[u^2(t,y)+\frac{u_x^2(t,y)}2\right]dy\,,
\\
\dis
P^v_x\big(t,\,\psi_{(t)}(x)\big)&=&
\dis \frac 12\int_0^1 \chi'\big(
\psi_{(t)}(x)-\tilde y\big)
\cdot v^2(t,\tilde y)\,d\tilde y\,\\
&&\quad 
\dis +\int_{x-1}^x \chi'\big(
\psi_{(t)}(x)-\psi_{(t)}(y)\big)
\cdot 
\frac {v_x^2\big(t,\psi_{(t)}(y)\big)} 4
\psi'_{(t)}(y)\,dy\,,\\
\end{array}
$$
where, according to (\ref{ch1-formachi}),
$$
\chi'(x)= \frac{ e^x-e^{1-x}}{e-1}\qquad 0<x<1\,,\qquad\qquad
\chi'(x)=\chi'(x+1)\qquad x\in\R\,.
$$

In the next computation, we use the estimate
$$
\int_0^1 \bigg|\big( 1+u_x^2(y)\big)- \Big( 1+v_x^2\big(\psi(y)\big)
\Big)\psi'(y)\bigg|\,dy=\O(1)\cdot J^\psi (u,v)\,.
$$
which holds because of the last two terms in the definition (\ref{ch1-defunzjei}).
Observing that $\chi'$ is Lipschitz continuous on the open interval
$]0,1[\,$, we now compute (omitting explicit references to the time
$t$)
\begin{equation}
\begin{array}{l}
\dis 
\Big|P^u_x(x) -P^v_x\big(\psi(x)\big)\Big|
\leq
\frac 12 \int_0^1 \left|\chi'(x-y)\cdot u^2(y)-\chi'\big(\psi(x)-y\big)
\cdot v^2(y)\right|\,dy
\\
\qquad\qquad\qquad \dis +\O(1)\cdot\int_{x-1}^x \Big|x-y-(\psi(x)-\psi(y)\big)\Big|\cdot
\frac{v_x^2\big(\psi(y)\big)}2\,
\psi'(y)\,dy 
\\
\dis
\qquad\qquad\qquad+\frac 14 \left|
\int_{x-1}^x
\chi'(x-y)\,\Big(u_x^2(y)-v_x^2\big(\psi(y)\big)\,\psi'(y)\Big)\,dy
\right|
\\
\qquad \dis = \O(1)\cdot \Big(\big|x-\psi(x)\big|+ \|u^2-v^2\|_{L^1}\Big)\\
\dis \quad\qquad
+\O(1)\cdot \left(\big|x-\psi(x)\big|+ \int_{x-1}^x \big|y-\psi(y)\big|\cdot
\frac{v_x^2\big(\psi(y)\big)}2\,
\psi'(y)\,dy \right)
\\
\quad\qquad \dis
+\O(1)\cdot \left(J^\psi (u,v)+
\Big|\int_{x-1}^x \chi'(x-y)\cdot \big[\psi'(y)-1\big]
\,dy\Big|\right)
\\
\dis
 \qquad= \O(1)\cdot \big|x-\psi(x)\big|+ \O(1)\cdot J^\psi (u,v)
 \\
 \dis \qquad \quad +
\O(1)\cdot \left( \big|x-\psi(x)\big|+\int_{x-1}^x
\chi''(x-y)\cdot \big[\psi(y)-y\big]
\,dy\right)
\\
\dis \qquad=\O(1)\cdot \big|x-\psi(x)\big|+ \O(1)\cdot J^\psi (u,v)\,.
\end{array}
\label{ch1-pumenopv}
\end{equation}
Integrating over one period we conclude
$$
I_2 = \O(1)\cdot J^\psi \big(u(t), \,v(t)\big)\,.
$$

For future use, we observe that a computation entirely similar to
(\ref{ch1-pumenopv}) yields
\begin{equation}
\Big|P^u(x) -P^v\big(\psi(x)\big)\Big|
=\O(1)\cdot \big|x-\psi(x)\big|+ \O(1)\cdot J^\psi (u,v)\,.
\label{ch1-Puxmenopvx}
\end{equation}

Next, we look at the term
$$
I_3\doteq \int_0^1 \frac d{dt}
\Big| 2\arctan\, u_x(t,x)
-2\arctan\, v_x\big(t,\psi_{(t)} (x)\big)
\Big|\cdot\phi_1^{(t)}(x)\,\big( 1+u_x^2(t,x)\big)\,dx\,.
$$
Along a characteristic, according to (\ref{ch1-nonlocalequation}) one has
$$\frac d{dt} 2\arctan u_x(t,\xi(t))=\frac 2{1+u_x^2}\Big[u^2-
\frac {u_x^2}2-P^u\Big]\,.$$
Call
$\theta^u\doteq 2\arctan u_x\,$,
$\theta^v\doteq 2\arctan v_x\,$,
so that
$$\frac 1{1+u_x^2}=\cos^2\frac {\theta^u}2\,,\qquad\quad
\frac {u_x}{1+u_x^2}=\frac 12\sin \theta^u\,,\qquad\quad 
\frac{u_x^2}{1+u_x^2}=\sin^2 \frac{\theta^u}2\,.$$
We now have
\begin{equation}
\begin{array}{l}
\dis 
\int_0^1 \big(1+u_x^2(x)\big)\cdot
\left| \frac{u_x^2(x)}{1+u_x^2(x)}- 
\frac{v_x^2(\psi(x))}{1+v_x^2(\psi(x))} \right|\,dx
\\
\qquad \dis
=
\int_0^1 \big(1+u_x^2(x)\big)\cdot
\Big| \sin^2\frac{\theta^u(x)} 2-\sin^2 \frac{\theta^v\big(\psi(x)\big)}2
\Big|
\,dx
\\
\qquad\dis \leq \int_0^1 \big(1+u_x^2(x)\big)\cdot \Big|\theta^u(x)-\theta^v
\big(\psi(x)\big)\Big|\,dx
\\
\qquad \dis=\O(1)\cdot J(u,v)\,.
\end{array}
\label{ch1-derarctang}
\end{equation}
Next, using (\ref{ch1-Puxmenopvx}) we compute
$$
\begin{array}{l}
\dis
\int_0^1\big(1+u_x^2(x)\big)\cdot \bigg|
\frac {u^2(x)-P^u(x)}{1+u_x^2(x)}-
\frac { v^2\big(\psi(x)\big)-P^v\big(\psi(x)\big)}
{1+v_x^2\big(\psi(x)\big)}\bigg|\,dx\\
\qquad \dis \leq
\int_0^1
\Big|u^2(x)-v^2\big(\psi(x)\big)\Big|\,dx
+\int_0^1
\Big|P^u(x)-P^v\big(\psi(x)\big)
\Big|\,dx\\
\qquad \dis \qquad\qquad+\O(1)\cdot \int_0^1 \left| \frac 1{1+u_x^2(x)}-
\frac 1{1+v_x^2\big(\psi(x)\big)}\right|\cdot \big(1+u_x^2(x)\big)\,dx
\\
\qquad\dis =\O(1)\cdot J^\psi(u,v)\,,
\end{array}
$$
where the last term was estimated by observing that
$$
\left|\frac 1{1+u_x^2}-
\frac 1{1+v_x^2}\right|\leq
\big|2\arctan u_x-2\arctan v_x\big|_*
\,.
$$
Putting together all previous estimates we conclude
$$
I_1+I_2+I_3=\O(1)\cdot J^\psi (u,v)\,.
$$
\vs

To complete the analysis, we have to consider the terms due to the change
in base measures. From (\ref{ch1-uxquadro}) it follows that
the production of new mass in the base measures is described
by the balance laws
$$\left\{
\begin{array}{l}
\dis (1+u_x^2)_t+ \big[ u(1+u_x^2)\big]_x=
\dis 
[2u^2+1-2P^u]u_x\doteq
f^u\,,
\\
\dis 
(1+v_x^2)_t+ \big[ v(1+v_x^2)\big]_x
=
\dis 
 [2v^2+1-2P^v]v_x \doteq
f^v\,.
\end{array}\right.
$$
This leads us to consider two further integrals $I_4,I_5$ :
$$\begin{array}{rcl}
I_4&=&\dis 
\int_0^1
 d^\diamondsuit\Big( \big(x,\,u(x),\, 2\arctan u_x(x)\big)\,,~
\big(\psi(x),\,v(\psi(x)),\,2\arctan v_x(\psi(x))\Big) \\
&&\dis \qquad\qquad \cdot
\big|2u^2(x)+1-2P^u(x)\big|\,\big|u_x(x)\big|\,dx\\
&=&\dis
\O(1)\cdot
\int_0^1
d^\diamondsuit\Big( \big(x,\,u(x),\, 2\arctan u_x(x)\big)\,,~
\big(\psi(x),\,v(\psi(x)),\,2\arctan v_x(\psi(x))\Big) 
\\
&&
\dis
\qquad\qquad \qquad\qquad\cdot
\big(1+u_x^2(x)\big)\,dx\cr
&=&
\dis 
\O(1)\cdot J^\psi(u,v)\,.
\end{array}
$$

$$\begin{array}{rcl}
\dis 
I_5
&=&
\dis
\int_0^1\bigg|
\big[2u^2(x)+1-2P^u(x)\big]u_x(x)-\Big[ 2v^2\big(\psi(x)\big)
+1-2P^v\big(\psi(x)\big)
\Big] v_x\big(\psi(x)\big)\,\psi'(x)\bigg|\,dx
\\
&\leq&\dis 2
\int_0^1\Bigg\{\Big| u^2(x)-v^2\big(\psi(x)\big)\Big|+\Big|P^u(x)-
P^v\big(\psi(x)\big)
\Big|\Bigg\}\, \big|u_x(x)\big|\,dx
\\
&&\dis\qquad
+\int_0^1 \Big|2v^2\big(\psi(x)\big)
+1-2P^v\big(\psi(x)\big)
\Big|\cdot \Big|u_x(x)-v_x\big(\psi(x)\big)
\psi'(x)
\Big|\,dx
\\
&=&
\dis I_5'+I_5''
\end{array}
$$

Using (\ref{ch1-Puxmenopvx}) we easily obtain
$$
\begin{array}{rcl}
I_5'&\leq&
\dis \int_0^1\Bigg\{\Big| u^2(x)-v^2\big(\psi(x)\big)
\Big|+\Big|P^u(x)-
P^v\big(\psi(x)\big)
\Big|\Bigg\}\, \big(1+u_x^2(x)
\big)\,dx\,\\
&=&\dis
\O(1)\cdot J^\psi(u,v)\,.
\end{array}
$$
On the other hand, recalling (\ref{ch1-derarctang}) and using the change of variable
$y=\psi(x)$, $x=\psi^{-1}(y)$, we find
$$\begin{array}{rcl}
I_5''&=&
\dis
\O(1)\cdot
\int_0^1 \Big|u_x(x)-v_x\big(\psi(x)\big)
\psi'(x)
\Big|\,dx
\\
&=&
\dis
\O(1)\cdot\int_0^1  \left|\frac{u_x(x)}{1+u_x^2(x)}-
\frac {v_x\big(\psi(x)\big)}{1+v_x^2\big(\psi(x)\big)}
\right|\big(1+u_x^2(x)\big)\,dx\\
&&\dis\qquad
+\O(1)\cdot\int_0^1
\Big|v_x\big(\psi(x)\big)\,\psi'(x)\Big|\cdot\left|
\frac{1+u_x^2(x)}{\big(1+v_x^2(\psi(x))\big)\,\psi'(x)}-1\right|
\\
&\leq&\dis
 \O(1)\cdot J^\psi(u,v)
+\int_0^1
\Big|\big(1+u_x^2(x)\big)-
\big(1+v_x^2(\psi(x))\big)\,\psi'(x)\Big|\,dx
\\
&=&
\dis 
\O(1)\cdot J(u,v)\,.
\end{array}
$$

All together, the previous estimates show that
\begin{equation}
\frac d{dt} J^{\psi_{(t)}}\big( u(t),\,v(t)\big)\leq
I_1+I_2+I_3+I_4+I_5'+I_5'' =\O(1)\cdot
J^{\psi_{(t)}}\big( u(t),\,v(t)\big)\,,
\label{ch1-stabilitajei}
\end{equation}
where $\O(1)$ denotes a quantity which remains uniformly bounded
as $u,v$ range on bounded sets of $H^1_\per\,$.
As an immediate consequence we obtain
\begin{lemma}
\label{lem7}
Let $t\mapsto u(t)$, $t\mapsto v(t)$
be two conservative, spatially periodic
multipeakon solutions, as in Lemma 2.
Then there exists a constant $\kappa$, depending only on
$\max\big\{\|u\|_{H^1_\per}\,,~\|v\|_{H^1_\per}\big\}$, such that
\begin{equation}
J\big( u(t)\,,~v(t)\big)\leq e^{\kappa|t-s|}\cdot
J\big( u(s)\,,~v(s)\big)\qquad\qquad s,t\in\R\,.
\label{ch1-boundsujey}
\end{equation}
\end{lemma}

\begin{proof}
For $t>s$ the estimate (\ref{ch1-boundsujey}) follows
from (\ref{ch1-stabilitajei}), taking the infimum among all transportation
plans $\psi_{(s)}$ at time $s$.   The case $t<s$
is obtained simply by observing that the Camassa-Holm
equations are time-reversible.
\end{proof}

\section{A priori bounds}\label{decayinfty}
In \cite{C3} the author discusses the finite propagation speed property for the Camassa-Holm equation. Due to the nonlocal nature of the equation \ref{ch1-equation}, it is not \emph{a priori} clear that the evolution of an initial data with compact support will remains  with compact support. On the contrary, \cite{C3} prove that the finite propagation speed property is valid only for the function $u-u_{xx}$ and not for $u$.
In this section we start from this result, and we want to estabilish what is the ``right" decay at the infinity of solutions to \ref{ch1-equation}. For this purpose, we introduce the following functional space. Let $\alpha\in \,]0,1[$, then we set
\begin{equation}
\label{decayHP}
X_\alpha\doteq \{u\in H^1(\R) \st C^{\alpha,u}\doteq \int_\R \left[u^2(x)+u_x^2(x)\right]e^{\alpha|x|}\,dx <+\infty\}.
\end{equation}
The present section is devoted to the study of some useful properties of the functions $u\in X_\alpha$. We start recalling an estimate for the $L^\infty-$norm of the $\H$ functions. We have
\begin{equation}
\label{sincichinequality}
\|f^2\|_{L^\infty}\leq \|f\|^2_{\H}.
\end{equation}
This estimate give us a bound on the $L^\infty-$norm of the conservative solution $u$ of (\ref{prblCH}), in fact the conservation of the energy yields 
\begin{equation}
\label{normainfty}
\|u(t)\|_{L^\infty}\leq \|u(t)\|_\H=\sqrt{E^{\bar u}} \qquad \mbox{for every $t\geq 0$.}
\end{equation}
Let us consider now the behaviour of the functions $u\in X_\alpha$ as $|x|$ goes to infinity. If we denote with $C^{\alpha,u}$ the constant $\int_\R  (u^2+ u_x^2)e^{\alpha|x|} \,dx$, the following holds
\begin{equation}
\label{uinfty}
\sup\limits_{x\in \R}\,u^2(x)e^{ \alpha |x|}\leq 2 C^{\alpha,u} \,.
\end{equation}
Indeed, the function
$$
f(x)\doteq u(t,x)e^{\frac\alpha 2|x| }
$$
belongs to $H^1(\R)$, moreover
$$
f_x=u_x e^{\frac\alpha 2|x| }+\frac \alpha 2 \sign(x)ue^{\frac\alpha 2|x| }
$$
and then, by using (\ref{sincichinequality}), we have 
$$
|f(x)|^2\leq\|f\|_\H^2\leq \int_\R [2u_x^2+ (1+\alpha)u^2]e^{\alpha |y|}\,dy \leq 2 C^{\alpha,u}\,.
$$
Now we study the behaviour at infinity of the multipeakon solutions of the Camassa-Holm equation. 
\begin{lemma}({\rm A-priori} bounds)
Let $u$ be a multi-peakon solution to (\ref{ch1-equation}), with initial data $\bar u$ which belongs to the space $X_\alpha$  defined in (\ref{decayHP}). Then for every $t\in \R$ there exists a continuous function $C(t)$, which depends on $C^{\alpha, \bar u}$ and on the energy $E^{\bar u}$, such that
\begin{eqnarray}
 \displaystyle
&\bullet&\int_{\R}[u^2(t,x)+u_x^2(t,x)]e^{\alpha|x|}\,dx\leq C(t)\,,
\label{weightedenergy}
\\&
\bullet&\displaystyle
\label{Pinfty}
\sup\limits_{x\in \R}\, \big|P_x^u(t,x)\big| e^{ \alpha |x|}\leq C(t)\,,
\\&\bullet&
\|u_x\|_{L^1(\R)}\leq C(t).\label{uxexpdecay}
\end{eqnarray} 
\end{lemma}
\begin{proof}
Since $|P_x^u|=P^u$, it is sufficient to prove the second inequality with $P^u_x$ replaced bu $P_u$.
Setting 
$$
I(t)\doteq\int_{\R}[u^2(t,x)+u_x^2(t,x)]e^{\alpha|x|}\,dx\,,
$$
we want to achieve a differential inequality of the form
$$
\frac d{dt}I(t)\leq A+B\cdot I(t)\,,
$$
for some constants $A$ and $B$ which depend on the initial data $\bar u$. 
We start the discussion proving a preliminary estimate for the function $P^u$. By applying the Fubini theorem to the identity
\begin{equation}
\label{stimaP}
\int_{\R}P^u(t,x)e^{\alpha|x|}\,dx= \frac 12 \int_\R e^{\alpha |x|}\,dx 
\int_\R e^{-|x-y|}\left[u^2(t,y)+\frac{u_x^2(t,y)}{2}\right]dy 
\end{equation}
we have to compute the following integral 
\begin{equation}
\int_\R e^{\alpha|x|}e^{-|x-y|}\,dx=\frac {2\alpha}{1-\alpha^2}e^{-|y|}+\frac 2{1-\alpha^2}e^{\alpha |y|}\qquad \mbox{for every $y\in\R$}\,.
\label{expint}
\end{equation}
For future use, we observe that the equality (\ref{expint}) holds for $\alpha \in (-1,1)$. 
Substituting (\ref{expint}) in (\ref{stimaP}) and using the definition of the energy $E^{\bar u}$ we have
$$
\int_{\R}P^u(t,x)e^{\alpha|x|}dx\leq \frac{E^{\bar u}}{(1-\alpha^2)} +
\frac {1}{1-\alpha^2}\,I(t)\,.
$$
Having in mind the previous inequality, we are able to estimate the time derivative of the function $I$. From the equations (\ref{ch1-equation}) and (\ref{ch1-uxquadro}) we have
$$
\begin{array}{rl}
\displaystyle
\frac d{dt}I(t)&
\displaystyle
\!\!\!= \int_{\R}\left[2u u_t+ \left(u_x^2\right)_t\right]\,e^{\alpha |x|}\,dx 
\\
&\!\!\!=\dis
\int_{\R}\left[ -2u(u_x+P^u_x)+\frac 23 (u^3)_x- (uu_x^2)_x- 2 u_x P^u\right]\,e^{\alpha |x|}\,dx  
\\
&\displaystyle
\!\!\!\leq
-2\int_\R
(u P^u +u u_x^2)_x\,e^{\alpha|x|}\,dx \leq
\left. -2u (P^u + u_x^2)\,e^{\alpha|x|}\right|_{-\infty}^{\infty}
\\&\!\!\!\quad
\dis + 2\alpha \int_\R |u| (P^u+ u_x^2)\,e^{\alpha|x|}\,dx
\\
&\displaystyle\!\!\!\leq
\frac{2\alpha\|u\|_{L^\infty}}{1-\alpha^2}\left(E^{\bar u}+ 2 I\right)
\leq 
\frac{2\sqrt{E^{\bar u}}}{1-\alpha^2}\left[E^{\bar u}+ 2 I(t)\right]
\end{array}
$$
the previous inequality gives then a bound on the function $I$, that is
$$
I(t)\leq (C^{\alpha,\bar u}+E^{\bar u}/2)\exp\Big(\frac{4\sqrt{E^{\bar u}}}{1-\alpha^2}\, t \Big)\,.
$$
To achieve the estimate (\ref{Pinfty}), set
$$
K(t)\doteq \left\|P^u(t,\cdot)e^{\alpha|\cdot|}\right\|_{L^\infty}\,.
$$
Proceeding as before, fixed $x\in \R$ we compute the derivative w.r.t the time $t$ of the function $e^{-|x|}*{u^2_x}$. 
$$
\begin{array}{rl}
\displaystyle\frac \partial{\partial t} \left(e^{-|x|}*\frac{u^2_x}4\right)
&\displaystyle =\frac 14 \frac \partial{\partial t}\int_\R e^{-|x-y|}u_x^2(t,y)\,dy
\\
&\displaystyle
=\frac 12\int_\R e^{-|x-y|}\left[\Big( \frac {u^3}3- \frac{u u_x^2}2 - u P^u\Big)_x+ u P_x^u\right] dy
\\
&\displaystyle\leq \|u\|_{L^\infty} P^u+\frac  {\|u\|_{L^\infty}}2 \int_\R e^{-|x-y|}P^u(t,y)\,dy
\\
&\displaystyle\leq
 \|u\|_{L^\infty}P^u(t,x)+\frac{\|u\|_{L^\infty}}2 K(t)\int_\R e^{-\alpha|x|}e^{-|x-y|}\,dy
\\
&\leq \|u\|_{L^\infty} P^u(t,x)+\frac 1{1-\alpha^2}e^{-\alpha|x|}\|u\|_{L^{\infty}}K(t)
\end{array}
$$
in the same way, the derivative of $e^{-|x|}*u^2$ is
$$
\begin{array}{rl}
\frac \partial{\partial t}\left(e^{-|x|}*\frac{u^2}2\right)
&\displaystyle \leq \int_\R e^{-|x-y|}u \left(|P^u_x|+|u u_x|\right)\,dy
\\
&\displaystyle \leq\|u\|_{L^\infty} \left( 2 P^u(t,x)+\int_\R e^{-|x-y|}P^u(t,y)\,dy\right) 
\\
&\displaystyle \leq \|u\|_{L^\infty}\left( 2P^u(t,x)+\frac 1{1-\alpha^2} e^{-\alpha|x|}K(t)\right)\,.
\end{array}
$$
Multiplying the previous two inequalities with $e^{\alpha|x|}$ we get
$$
\frac d{dt}K(t) \leq \left(3+\frac{ 2}{1-\alpha^2} \right) \sqrt{E^{\bar u}} K(t)
$$
which yields (\ref{Pinfty}). 

To achieve the last inequality, we write
\begin{eqnarray*}
\int_\R|u_x(y)|\,dy
&=&\!\!\!\!\int\limits_{ \{y :|u_x(y)|e^{\alpha|y|}<1 \} } \!\!\!\!|u_x(y)|\,dy
+\!\!\!\!\int \limits_{\{y:|u_x(y)|e^{\alpha|y|}>1\}} \!\!\!\!|u_x(y)|\,dy
\\
&\leq& \int_\R e^{-\alpha|y|}\,dy+\int_\R u_x^2(y)e^{\alpha|y|}\,dy\leq \frac 2\alpha + I(t)
\end{eqnarray*}
where the last estimate is given by (\ref{weightedenergy}).
\end{proof}

\section{Definition of the distance in the real line\label{metricsection}}
In this section we define a metric in order to control the distance between two solutions of the equation (\ref{prblCH}) whenever we take initial data in the subspace $X_\alpha$ of $H^1(\R)$. It is constructed as in the spatially periodic case by resolving an appropriate optimal transportation problem. Let $\T=[0,2\pi]$ be the unit circle with the end points $0$ and $\pi$ identified. \begin{figure}[ht]
\centerline{\includegraphics[height =5cm]{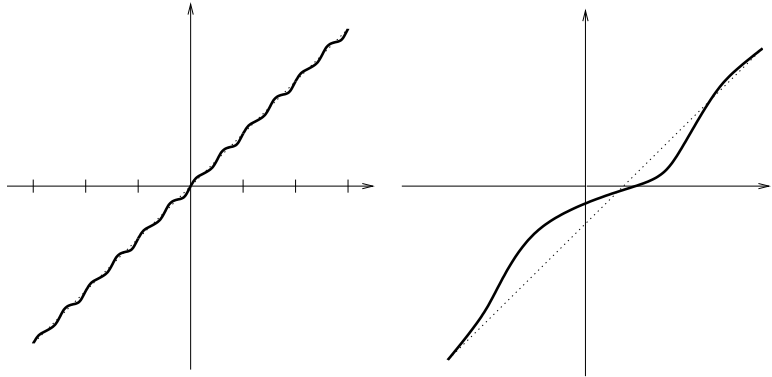}}
\caption{``Periodic'' vs. ``non-periodic'' transportation map \label{confronti}}
\end{figure}

\noindent Consider the metric space $(\R^2\times \T,d^\diamondsuit)$, with distance 
$$
d^\diamondsuit((x,u,\omega),(x',u',\omega'))\doteq \min\{|x-x'|+|u-x'|+|\omega-\omega'|_*,1\}
$$
and for every function $u\in X_\alpha$, let us define the Radon measure on $\R^2\times \T$
$$
\sigma^u(A)\doteq \int_{\{x\in \R: (x,u(x),2\arctan u_x(x))\in A\}} [1+u_x^2(x)]\,dx 
$$
for every Borel set $A$ of  $\R^2\times \T$.
\begin{dhef}\label{deftranmas}
\rm
 The set $\mathcal F$ of transportation plans consists of the functions $\psi$ with the following properties:
\begin{enumerate}
\item \label{conduno} $\psi$ is absolutely continuous, is increasing with its inverse;
\item \label{condue} $\sup\limits_{x\in\R}|x-\psi(x)|e^{\alpha/2|x|}<\infty$;
\item \label{contre} $\int_{\R} |1-\psi'(x)|\,dx<\infty$.
\end{enumerate}
\end{dhef}
\psfrag{x}{$x$}
\psfrag{u}{$u$}
\psfrag{v}{$v$}
\begin{figure}[ht]
\centerline{
\includegraphics[width=10cm]{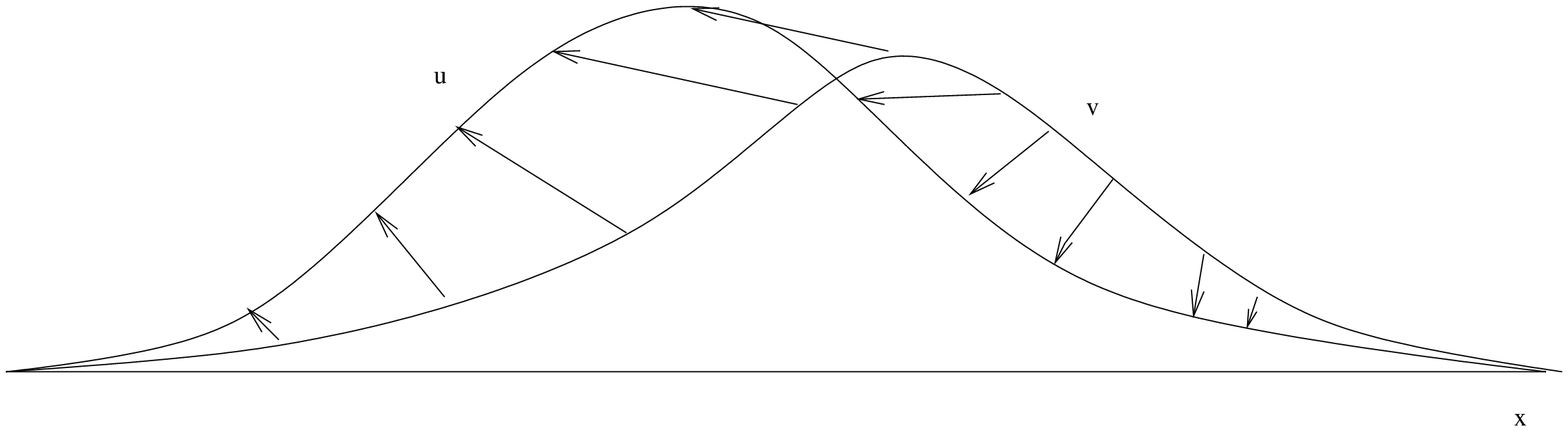}
}\caption{Transportation plan.\label{ftmass}}
\end{figure}  
The conditions \ref{condue} and \ref{contre} are not restrictive. Indeed, thanks to the exponential decay of functions $u,v \in X_\alpha$, the measures $\sigma^u$ and $\sigma^v$ located on the graph of $u$ and $v$ respectively, have small mass at the infinity, and then a transportation plan which transports mass from one to the other can be almost the identity $\psi(x) \approx x$ (see fig. \ref{ftmass}).
In order to define a distance in the space $X_\alpha$, we consider an optimization problem over all possible transportation plans. Given two functions $u,\,v$ in $X_\alpha$, we introduce two further measurable functions, related to a transportation plan $\psi$:
\begin{eqnarray}
\label{non-phi1}\phi_1(x)\doteq \sup \big\{ \theta \in [0,1] \st \theta \cdot(1+u_x^2(x))\leq \left(1+v_x^2(\psi(x))\right)\psi'(x)\big\},
\\
\label{non-phi2}\phi_2(\psi(x))\doteq \sup \big\{ \theta \in [0,1] \st  1+u_x^2(x)\leq \theta\cdot \left(1+v_x^2(\psi(x))\right)\psi'(x)\big\}.\phantom{cc}
\end{eqnarray}
The functions $\phi_1, \phi_2$ can be seen as weights that take into account the difference of the masses of the measure $\sigma^u$ and $\sigma^v$. In fact, from the definitions (\ref{non-phi1})-(\ref{non-phi2}) one has
$$
\phi_1(x) (1+u_x^2(x))=\phi_2(\psi(x)) (1+v_x^2(\psi(x)))\psi'(x)\qquad \mbox{for a.e. $x\in\R$}.
$$
According to the definitions, the identity $\max\{\phi_1(x),\phi_2(x)\}\equiv 1$ holds.
Altough the two measures $\phi_1\sigma^u$ and $\phi_2\sigma^v$ have not finite mass, they satisfy $\phi_1\sigma^u(A)=\phi_2\sigma^v(A)$ for every bounded Borel set $A\subset \R^2\times \mathbb T$. Thus, the functions $\phi_1$ and $\phi_2$ represent the percentage of mass actually transported from one measure to the other. 
A distance between the two functions $u,\,v$ in $X_\alpha$ can be characterized in the following way.

For every transport map $\psi\in \mathcal F$, let define ${\bf X}^u=(x,u(x),2\arctan u_x(x))$ and ${\bf X}^v=(\psi(x),v(\psi(x)),2\arctan v_x(\psi(x)))$ and consider the functional
$$
J^\psi(u,v)=\int_\R d^\diamondsuit ({\bf X}^u,{\bf X}^v)\phi_1(x)(1+u_x^2(x))\,dx+ \int_\R \left|1+u_x^2(x)-(1+v_x^2(\psi(x)))\psi'(x)	\right|\,dx\,.
$$
Since the above functional is well defined for every $\psi\in \mathcal F$, we can define 
$$
J(u,v)\doteq \inf_{\psi \in \mathcal F} J^\psi (u,v).
$$
The functional $J$ here defined is thus a metric on the space $X_\alpha$ (see the previous Section \ref{3-3}). 

\medskip
\subsection{Comparison with other topologies}
\begin{lemma}\label{L1J}
For every $u,\,v\in X_\alpha$ one has 
\begin{equation}
\label{h1el1}
\frac 1C\cdot \|u-v\|_{L^1(\R)}\leq J(u,v)\leq C\cdot \|u-v\|_{\H}.
\end{equation}
Let $(u_n)$ be a Cauchy sequence for the distance $J$ such that $C^{\alpha, u_n}\leq C_0$ for every $n\in \N$. Then
\begin{enumerate}
\item[i)] There exists a limit function $u\in X_\alpha$ such that $u_n\to u$ in $L^\infty$ and the sequence of derivatives ${u_n}_x$ converges to $u_x$ in $L^p(\R)$ for $p \in [1,2[$.
\item[ii)] Let $\mu_n$ be the absolutely continuous measure having density ${u_n}_x^2$ with respect to Lebesgue measure. Then there exists a measure $\mu$ whose absolutely continuous part has density $u_x^2$ such that $\mu_n \rightharpoonup \mu$. 
\end{enumerate}
\end{lemma}
\begin{proof}
The first inequality of (\ref{h1el1}) can be achieved by estimating the area between the two functions $u$ and $v$. For every $\psi\in \F$ we can write
$$
\int_\R |u-v|\, dx =\int_{S_1}|u-v|\,dx+\int_{S_2}|u-v|\,dx
$$ 
where the two subsets $S_1$ and $S_2$ are
\begin{itemize}
\item $S_1= \{x: |x-\psi(x)|\leq~1~\}=\cup_j [x_{2j-1},x_{2j}]$, where in this union we have to take into account that these intervals may be either finite or infinite, possibly having $x_j=\pm \infty$ for some $j$, 
\item $S_2=\{x: |x-\psi(x)|>1\}$. 
\end{itemize}
The integral over $S_2$ can be estimate in the following way:
\begin{equation}
\label{psigrosso}
\int_{S_2}|u(x)-v(x)|\,dx\leq (\|u\|_{L^\infty}+\|v\|_{L^\infty}) \int_\R|x-\psi(x)|\,dx\leq (E^{\bar u}+ E^{\bar v})J(u,v).
\end{equation}
The last inequality is given by the definition of the functional $J$.
\psfrag{An}{$A_{n}$}
\psfrag{An+1}{$A_{n+1}$}
\psfrag{An+k}{$A_{n+k-1}$}
\psfrag{An+j}{$A_{n+k}$}
\psfrag{u}{$u$}
\psfrag{v}{$v$}
\psfrag{x0}{$x_0$}
\psfrag{x1}{$x_1$}
\psfrag{x2k}{$x_{2k}$}
\psfrag{x2k+1}{$x_{2k+1}$}
\psfrag{xk}{$x_k$}
\psfrag{x-psix<1}{$S_1$}
\psfrag{x-psix>1}{$S_2$}
\psfrag{n}{$n$}
\psfrag{n+1}{$n+1$}
\psfrag{n+j}{$n+j$}
\psfrag{n+k}{$n+k$}
\psfrag{...}{$...$}
\psfrag{PQ}{$\overline {Q_u(n)\,Q_v(n)}$}
\begin{figure}[ht]
\centerline{\includegraphics[width=8cm]{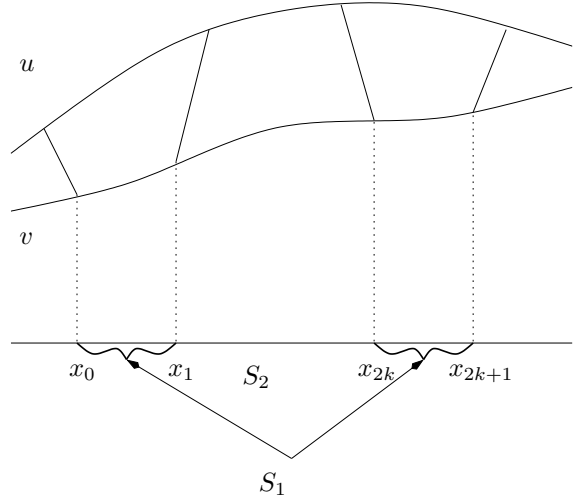}}
\caption{$L^1-$distance between two functions.\label{L1norma}}
\end{figure}

As far as the integral over $S_1$ is concerned, the integral over $S_1$ can be viewed as a sum of the area of the regions $A_j$ in the plane $\R^2$, bounded by the graph of the curves $u$, $v$ and by the segments with slope $\pm 1$ that join the points $Q_u(x_{2j-1})=(x_{2j-1},u(x_{2j-1}))$ and $Q_v(x_{2j})=(\psi(x_{2j}),v(\psi(x_{2j})))$, where $\{x_i\}=\partial S_1$. 
We have
$$
\int_{S_1}|u(x)-v(x)|\,dx\leq \sum_j \meas(A_j)\,.
$$
The measure of the subset $A_j$ is the area sweeped by the segment $\overline {Q_u(x)\,Q_v(x)}$. Recalling that in every set $A_j$ the function $\psi$ satisfies $ |x|-1\leq |\psi(x)|\leq |x|+1$, a bound on this area is given by
$$
\meas (A_j) \leq
\displaystyle
\int_{x_{2j-1}}^{x_{2j}} (|x-\psi(x)|+|u(x)-v(\psi)|)[1+u_x^2+(1+v_x^2(\psi ))\psi']\,dx
$$
and then 
$$
\begin{array}{rl}
\dis
\int_{S_1}|u(x)-v(x)|\,dx
&\dis 
\leq \int_{S_1}(|x-\psi(x)|+|u(x)-v(\psi)|)[1+u_x^2+(1+v_x^2(\psi ))\psi']\,dx
\\
&\dis\leq J^\psi(u,v) + J^{\psi^{-1}} (u,v)
\end{array}
$$
this inequality, together with (\ref{psigrosso}), yields to 
$$
\|u-v\|_{L^1(\R)}\leq C(\bar u,\bar v)\, J(u,v).                                                   
$$   
Concerning the second part of the lemma, let us observe that even if the embedding of $H^1(\R)$ in $L^2(\R)$ is not compact, the uniform exponential decay of the function $u_n$ ($C^{\alpha,u_n}\leq C_0$ uniformly in $n$) allows us to extract a subsequence which converges to a function $u$ in $L^1$-norm and, by the unformly H\"oder continuity od such functions, $\|u_n-u\|_{L^\infty}\to 0$. This allow us to prove the property ii) by dealing with analogous arguments to those developed in the periodic case. The proof of the second part of the lemma is thus perfectly similar to the one of the periodic case, once we take into account the exponential decay of the sequence $u_n$.   
\end{proof}
\subsection{Continuity of solutions w.r.t initial data}\label{stabcor}
Let $u_0$ and $v_0$ be two multipeakon initial data. The technique developed in Section \ref{ODEsystem} ensures the existence of two multipeakon solutions $u(t),v(t)$ for (\ref{ch1-equation}) which conserve the energy unless interaction of peakon occurs. Suppose then that within a given interval $[0,T]$ no interaction occurs neither for $u(t)$ nor for $v(t)$. 
The aim of this section is to prove the continuity of the functional $J$ w.r.t. the initial data, namely we prove that there exists a continuous, positive function $C(t)$ such that for $t\in [0,T]$ one has
$$
J(u(t),v(t))\leq C(t) J(u_0,v_0).
$$
\begin{lemma}
If $u(t)$ and $v(t)$ are two multipeakon solutions defined in the interval $[0,T]$ in which no interaction occurs, then there exists a positive, continuous function $c(t)$ which depends only on the energies $E^u$, $E^v$ of the two solutions, such that
\begin{equation}
\frac d{dt} J(u(t),v(t))\leq c(t) J(u(t),v(t)) \qquad \mbox{for all $t\in[0,T]$}. 
\end{equation}
\end{lemma}
\begin{proof}
We compute the time derivative of the function $J^\psi(u(t),v(t))$ with a particular choice of the transportation plan $\psi=\psi_{(t)}$. Given any $\psi_0\in \mathcal F$, at every time $t\in [0,T]$ we construct $\psi_{(t)}$ by transporting the function $\psi_0$ along the characteristic curves. More precisely, since no interaction between peakon occurs in the interval $[0,T]$, the functions $u(t,\cdot),\,v(t,\cdot)$ are Lipschitz continuous, then the flows $\vphi^t_u$, $\vphi^t_v$ solutions of the Cauchy problems
\begin{eqnarray*}
&&\frac d{dt} \vphi^t_u(x)=u(t,\vphi_u^t(x))\qquad \vphi^0_u(x)=x,
\\
&&\frac d{dt}\vphi^t_v(y)=v(t,\vphi_v^t(y))\qquad \vphi^0_v(y)=y,
\end{eqnarray*}
which are the characteristics curves associated to the equation (\ref{ch1-equation}), are well defined. Now, let $x\in \R$. $\psi_{(t)}$ is defined as the composition
\begin{equation}
\psi_{(t)}(x)\doteq \vphi_v^t\circ \psi_0\circ \left(\vphi^t_u \right)^{-1}(x),
\end{equation}
that is 
$$
\psi_{(t)}(\vphi_u^t(y))=\vphi_v^t(\psi_0(y)).
$$
The function $\psi_{(t)}$ belongs to $\mathcal F$, and hence $J^{\psi_{(t)}}$ is well defined, in fact
\begin{enumerate}
\item By the Property \ref{conduno} of Definition \ref{deftranmas} for the function $\psi_0$, and uniqueness of solution of ODE, the function $\psi_{(t)}$ is an  increasing function.
\item Let $x\in \R$ and $\vphi_u^t(y)$ be the characteristic curve passing through $x$ at time $t$. Evaluating $|x-\psi_ {(t)}(x)|e^{\alpha/2|x|}$ along this characteristic curve, and computing the derivative w.r.t. $t$ we obtain
\begin{eqnarray*}
\dis
\frac d{dt}|&&\!\!\!\!\!\!\!\!\!\!\!\!\vphi_u^t(y)-\vphi_v^t(\psi_0(y))|e^{\alpha/2|\vphi_u^t(y)|}
\\
\dis &&\!\!\!\!\!\!\leq 
\left[|u(t,x)-v(t,\psi_{(t)}(x))|+\frac \alpha 2 |u(t,x)|\cdot |x-\psi_{(t)}(x)|\right]e^{\alpha/2|x|}
\end{eqnarray*}
by properties (\ref{uinfty}), (\ref{weightedenergy}), and since $u,\,v$ are Lipschitz continuous in $[0,T]$, there exists two  $L^\infty$ functions $c_1(t)$, $c_2(t)$ such that 
$$
\frac d{dt}|x-\psi_{(t)}(x)|e^{\alpha/2|x|}\leq c_1(t)|x-\psi_{(t)}(x)|e^{\alpha/2|x|}+c_2(t)
$$ 
by Gronwall Lemma and the hypothesis $|x-\psi_0(x)|e^{\alpha/2|x|}\leq C_0$, the previous inequality gives the Property \ref{condue} of Definition \ref{deftranmas} for $\psi_{(t)}$
\begin{equation}
\label{stpsi}
|x-\psi_{(t)}(x)|e^{\alpha/2|x|}\leq C_1(t)\doteq \left(C_0+ \int_0^t c_2(s)\,ds\right)e^{\int_0^t c_1(s)\,ds}
\end{equation}
\item The last property can be achieved by choosing the change of integration variable $x=\vphi_u^t~(y)$ 
$$
\begin{array}{rl}
\dis 
\int_\R |1-\psi_{(t)}(x)|\,dx
&
\dis 
=\int_\R |1-\psi_{(t)}(x)|(\vphi_u^t)'(y)\,dy 
\\
&
\dis
= \int_\R |(\vphi_u^t)'(y)-(\vphi_v^t)'(\psi_0(y))\psi_0'(y)|\,dy
\\
&\leq \dis \int_\R |(\vphi_u^t)'(y)-1|\,dy +\int_\R |(\vphi_v^t)'(y)-1|\,dy
\\
&
\dis\quad+\int_\R|1-\psi_0'(y)|\,dy.
\end{array}
$$
Since
$$
|(\vphi_u^t)'(y)-1|\leq \int_0^t |u_x(s,x)|\cdot |(\vphi_u^s)'(y)-1|\,ds+\int_0^t |u_x(s,x)|\,ds
$$
(and a similar estimate for $\vphi_v^t$) and $u_x,v_x\in L^\infty$, by the Gronwall lemma the first two integrals of the previous formula are bounded by an absolutely continuous function $C(t)$ in the interval $[0,T]$ and then also 
Property \ref{contre} of Definition \ref{deftranmas} holds.
\end{enumerate}
At the transportation plan $\psi_{(t)}$ we associate the functions $\phi_1^{(t)},\,\phi_2^{(t)}$ defined according to  (\ref{non-phi1}), (\ref{non-phi2}), the functional $J^{\psi_{(t)}}$ is thus 
$$
\begin{array}{rl}
J^{\psi_{(t)}}(u(t),v(t))=&
\dis\int_\R d^\diamondsuit(\mathbf X^u(t),\mathbf X^v(t))\phi_1^{(t)}(x)(1+u_x^2(x))dx
\\
&+\dis
\int_\R \left|1+u_x^2(x)
-(1+v_x^2(\psi_{(t)}(x)))\psi_{(t)}'(x)\right|dx.
\end{array}
$$
By deriving $J^{\psi_{(t)}}(u(t),v(t))$ w.r.t. $t$ and computing the change of variables along the characteristics, the previous derivative can be estimate by the sum of the following terms (we leave out the dependence on the integrable variable when it is not essential)
$$
\begin{array}{l}
\dis
I_1=\int_\R |u(t,x)-v(t,\psi_{(t)}(x))|\phi_1^{(t)}(x)(1+u_x^2(t,x))\,dx
\\
\qquad \qquad \qquad\dis \leq (1+\|u(t)\|_{L^\infty}+\|v(t)\|_{L^\infty}) J^{\psi_{(t)}}(u(t),v(t))\,,
\\
\\
\dis
I_2=\int_\R|P^u_x(t,x)-P^v_x(t,\psi_{(t)}(x))|\phi_1^{(t)}(x)(1+u_x^2(t,x))\,dx\,,
\\
\\
I_3=
\dis
\int_\R \left|\frac {2u^2(t)-u_x^2(t)-2P^u(t)}{1+u_x^2(t)}\right.
\left.-\frac{2v^2(t,\psi_{(t)})-v_x^2(t,\psi_{(t)})-2P^v(t,\psi_{(t)})} {1+v_x^2(t,\psi_{(t)})}\right|\cdot
\\
\dis
\qquad\qquad\qquad \cdot\phi_1^{(t)} (1+u_x^2(t))\,dx\,,
\end{array}
$$
the term due to the variation of the base measure 
$$
\begin{array}{rl}
\dis
I_4=2\int_\R d^\diamondsuit(\mathbf X^u(t),\mathbf X^v(t))\cdot u_x(t)(u^2(t)- P^u(t))\,dx\,,
\end{array}
$$
and the terms due to the variation of the excess mass
$$
I_5=\frac d{dt}\int_\R\left|1+u_x^2(t)-(1+v_x^2(t,\psi_{(t)}))\psi_{(t)}'\right|\,dx\,.
$$
Let us start to estimate the term $I_2$. By definition, the difference of $P^u$ and $P^v$ is written in convolution form
$$
\begin{array}{rl}
 &\!\!\!\!\!\!\!\!\!\dis \int_\R\left \{  e^{-|x-y|}\sign(x-y)\left[u^2(t,y)+\frac{u^2_x(t,y)}{2}\right]\,dy \right.
\\
&\left.\dis  -e^{-|\psi_{(t)}(x)-\psi_{(t)}(y)|}\sign(\psi_{(t)}(x)-\psi_{(t)}(y))\left[v^2(t,\psi_{(t)}(y))+\frac{v^2_x(t,\psi_{(t)}(y))}{2}\right]\psi_{(t)}'(y)\right\}dy
\end{array}
$$
by this inequality, we can estimate the term $I_2$ by the sum of the following integrals
$$
\begin{array}{l}
\displaystyle
\dis A=\int_\R(1+u_x^2(t,x))\int_\R e^{-|x-y|}|u^2(t,y)-v^2(t,y)|\,dy\,dx
\\
\dis B=\int_\R(1+u_x^2(t,x))\left|\int_\R e^{-|x-y|}\sign(x-y)[v^2(t,y)-v^2(t,\psi_{(t)}(y))\psi_{(t)}'(y)]\,dy \,\right|\,dx
\\
\dis C= \int_\R(1+u_x^2(t,x))\int_\R\left[v^2(t,\psi_{(t)} (y))+\frac {v_x^2(t,\psi_{(t)} (y))}2\right]\cdot
\\
\dis\qquad\quad \cdot\left|e^{-|x-y|}\sign(x-y)-e^{-|\psi_{(t)}(x)-\psi_{(t)}(y)|}\sign(\psi_{(t)}(x)-\psi_{(t)}(y))\right| \cdot
\psi_{(t)}'(y)\,dy\,dx
\\
\dis D=\frac 12\int_\R(1+u_x^2(t,x)) \left|\int_\R e^{-|x-y|}\sign(x-y)[u_x^2(t,y)-v_x^2(t,\psi_{(t)}(y))\psi_{(t)}'(y)]\,dy\,\right|\,dx
\end{array}
$$
\bigskip\noindent{\bf A.} Switching the order of the two integrals, the term $A$ is bounded by the $L^1$-norm of the difference between $u$ and $v$:
$$
\begin{array}{rl}
A
& 
\dis 
\leq(\|u\|_{L^\infty}+\|v\|_{L^\infty})\int_\R |u(t,y)-v(t,y)|\int_\R e^{-|x-y|}(1+u_x^2(t,x))\,dx\,dy
\\
&\dis
\leq (2+E^u)(\|u\|_{L^\infty}+\|v\|_{L^\infty})\|u(t)-v(t)\|_{L^1}
\end{array}
$$
and then, by Lemma \ref{L1J}, $A\leq C(\bar u,\bar v)J(u(t),v(t))$.

\bigskip
\noindent{\bf B.}
Define 
$$
F(y)\doteq\int_{-\infty}^y (v^2(z)-v^2(\psi_{(t)}(z))\psi_{(t)}'(z))\,dz=\int_{\psi_{(t)}(y)}^{y}v^2(z)\,dz
$$
we have,
integrating by parts
$$
\begin{array}{rl}
\dis
\left|\int_\R e^{-|x-y|}\sign(x-y)F'(y)\,dy\,\right|
&
\dis
\leq 2|F(x)|+\int_R e^{-|x-y|}|F(y)|\,dy
\\
&\dis
\leq \|v\|^2_{L^\infty}|x-\psi(x)|+\int_\R e^{-|x-y|}|y-\psi(y)|\,dy
\end{array}
$$
moreover, substituting the previous expression into the term $B$ we obtain
\begin{eqnarray*}
B&\!\!\!\!\!\!\!\!\!\!&\leq 2E^{\bar v}J^{\psi_{(t)}}(u(t),v(t))+\int_\R (1+u_x^2(t,x)) \int_\R e^{-|x-y|} |y-\psi_{(t)}(y)|\,dy\,dx
\\
&\!\!\!\!\!\!\!\!\!\!&=2E^{\bar v}\left\{J^{\psi_{(t)}}(u(t),v(t))+\int_\R |y-\psi_{(t)}(y)|\int_\R(1+u_x^2(t,x))e^{-|x-y|} \,dx\,dy\right\}
\\
&\!\!\!\!\!\!\!\!\!\!&\leq 2E^{\bar v}(3+E^{\bar u})\cdot J^{\psi{(t)}}(u(t),v(t)).
\end{eqnarray*}

\bigskip
\noindent {\bf C.} Observe that since the function $y\mapsto \psi_{(t)}(y)$ is non decreasing, the quantities $x-y$ and $\psi_{(t)}(x)-\psi_{(t)}(y)$ have the same sign, and since the function $t \mapsto e^{-|t|}$ is Lipschitz continuous either in $(-\infty,0)$ or in $(0,+\infty)$ we have
$$
\begin{array}{rl}
\dis
\left|e^{-|x-y|}-e^{-|\psi_{(t)}(x)-\psi_{(t)}(y)|}\right|
&\leq \dis
e^{-\min\{|x-y|,|\psi_{(t)}(x)-\psi_{(t)}(y)|\}}\left||x-y|-|\psi_{(t)}(x)-\psi_{(t)}(y)|\right|
\\
&
\dis
\leq
e^{-\min\{|x-y|,|\psi_{(t)}(x)-\psi_{(t)}(y)|\}}\left(|x-\psi_{(t)}(x)|+|y-\psi_{(t)}(y)|\right)
\end{array}
$$ 
now 
$$
-\min\{|x-y|,|\psi_{(t)}(x)-\psi_{(t)}(y)|\}\leq -|x-y|+2C_1(t),
$$
where $C_1(t)$ is the function (\ref{stpsi}), related to the Property \ref{condue} of Definition \ref{deftranmas} $\psi_{(t)}$, then
$$
\begin{array}{rl}
C
\leq 
&\!\!\dis e^{C_1(t)}\!\!\int_\R |y-\psi(y)|
\left[v^2(t,\psi_{(t)}(y)) + \frac {v_x^2(t,\psi_{(t)}(y))}2\right]\psi_{(t)}'(y)\cdot
\\
&\dis\qquad\qquad\cdot \int_\R(1+u_x^2(t,x))e^{-|x-y|}\,dx\,dy+
+\dis 2 E^{\bar v}\int_\R(1+u_x^2(t,x))|x-\psi_{(t)}(x)|\,dx
\\
\leq 
& \!\!\dis
\left[(2+E^{u})e^{C_1(t)}(1+\|u\|_{L^\infty})+ 2 E^{\bar v}\right] J^{\psi_{(t)}}(u(t),v(t))
\end{array}
$$

\bigskip
\noindent{\bf D.} 
Here we can use the estimate given by the change in base measure. Since  
$$
\int_\R \Big|1+u_x^2(t,x)- \big(1+v_x^2(t,\psi_{(t)}(x))\big)\psi_{(t)}'(x)\Big|\, dx\leq J^{\psi_{(t)}}(u(t),v(t))\,,
$$
we obtain
$$
\begin{array}{rl}
D \leq&\!\!\!
\dis
\frac 12 \int_\R(1+u_x^2(t,x)) \int_\R e^{-|x-y|}
\left|(1+u_x^2(t,y))-(1+v_x^2(t,\psi_{(t)}(y)))\psi_{(t)}'(y)\right|  \,dy\,\,dx
\\
\quad&\!\!\!
\dis  + \frac 12 \int_\R(1+u_x^2(t,x))\cdot 
\left|\int_\R  e^{-|x-y|}\sign(x-y) \left[\psi_{(t)}'(y)-1\right]dy\, \right|\,dx
\\
\leq& \!\!\!
\dis
\int_\R(1+u_x^2(t,x))\left|(\psi_{(t)}(x)-x)-\frac{e^{-|x-y|}}2\int_\R (\psi_{(t)}(y)-y)\,dy\right|\,dx 
\\
&\!\!\!
\quad\dis+(1+ E^{\bar u})J^{\psi_{(t)}}(u(t),v(t))
\\
\leq & \!\!\!
2(2+ E^{\bar u})J^{\psi_{(t)}}(u(t),v(t))
\dis
\end{array}
$$
where in the last estimate we integrated by part as in the term $B$. 

\bigskip
\noindent The control for the terms $I_3$, $I_4$ and $I_5$ can be obtained exactly as the ones in \cite{BF2}, whom we refer the reader to.  
The previous estimates implies that there exists a smooth function $C=C^{\bar u,\bar v}(t)$ which depends only to the variable $t$ and to the initial data $\bar u,\ \bar v$ such that 
$$
\frac{d}{dt} J^\psi(u,v)\leq C^{\bar u,\bar v}(t) J^{\psi}(u,v)
$$
which yields
$$
J(u(t),v(t))\leq J(u(s), v(s)) e^{\left|\int_s^t C^{\bar u,\bar v}(\sigma)\,d\sigma\right|}\qquad {\mbox {for every $s,t\in\R$.}}
$$
\end{proof}

\section{Proof of the main theorems\label{3-6}}
Thanks to the analysis in the previous sections (\ref{3-3} for the spatially periodic case, \ref{decayinfty} and \ref{metricsection} in the whole real line),
we now all the ingredients
toward a proof of Theorems \ref{theo1} and \ref{theo2}.
Since the two cases are analogous, we prove them in the periodic case. 
The estimates in (\ref{ch1-weaktop}) follow from Lemma \ref{lem3}.
Given an initial data $\bar u\in H^1_\per$,
to construct the solution
of the Camassa-Holm equation
we consider a sequence of multi-peakons $\bar u_n$,
converging to $\bar u$ in $H^1_\per$.
Then we consider the corresponding solutions
$t\mapsto u_n(t)$, defined for all $n\geq 1$ and $t\in\R$.
This is possible because of Lemmas \ref{lem1} and \ref{lem2}.

We claim that the sequence $u_n(t)$ is Cauchy in $L^2_\per\,$.
Indeed, by Lemma \ref{lem3} and Lemma \ref{lem5},
$$
\begin{array}{rl}
\big\|u_m(t)-u_n(t)\big\|_{L^1_\per}
&\dis \leq C\cdot J\big(u_m(t),\,u_n(t)\big)
\\
&\dis \leq
C\cdot e^{\kappa |t|} \,J\big(u_m(0),\,u_n(0)\big)
\leq C^2\cdot e^{\kappa |t|}\big\|u_m(t)-u_n(t)\big\|_{H^1_\per}\,.
\end{array}
$$
Therefore, $u_n(t)\to u(t)$ in $L^1_\per$, for some function
$u:\R\mapsto H^1_\per\,$.   By interpolation, the convergence
$u_n\to u$ also holds in all spaces $L^p_\per$, $1\leq p\leq\infty$.
The continuity estimates (\ref{ch1-contu})-(\ref{ch1-stabsol}) now follow by passing to
the limit in Lemma \ref{lem5} and \ref{lem6}.

It remains to show that the limit function $u(\cdot)$
is actually a solution to the Camassa-Holm equation and
its energy $E(t)$ in (\ref{ch1-energy}) is a.e.~constant.
Toward these goals, we observe that all solutions $u_n$
are Lipschitz continuous with the same Lipschitz
constant, as maps from $\R$ into $ L^2_\per\,$. Indeed
$$
\|u_{n,t}\|_{L^2_\per}\leq \|u_n\|_{L^\infty}\cdot
\|u_{n,x}\|_{L^2_\per }+\left\| \frac 12\,e^{-|x|}\right\|_{L^2}
\cdot \left\|u_n^2+\frac{u^2_{n,x}}2\right\|_{L^1_\per}.
$$
As a consequence, the map $t\mapsto u(t)$ has uniformly bounded
$H^1_\per$ norm, and is Lipschitz continuous with values
in $L^2_\per\,$.   In particular,
$u$ is uniformly H\"older continuous as a function of $t,x$
and the convergence $u_m(t,x)\to u(t,x)$ holds uniformly
for $t$ in bounded sets.
Moreover, since $L^2_\per$ is a reflexive space,
the time derivative $u_t(t)\in L^2_\per$ is well
defined for a.e.~$t\in\R$.

We now observe that, for each $n\geq 1$,
both sides of the equality
\begin{equation}
\frac d{dt} u_n=-u_n\, u_{n,x} -P^{u_n}_x
\label{ch1-evolequinl2}
\end{equation}
are continuous as functions
from $\R$ into $L^1_\per\,$, and the identity holds
at every time $t\in\R$, with the exception of the
isolated times where a peakon interaction occurs.

At any time $t$ where no peakon interaction occur in the solution
$u_n$, we define
$\mu^{(n)}_t$ to be the measure with density
$u_n^2(t,\cdot)+\frac 12 u_{n,x}(t,\cdot)$ w.r.t.~Lebesgue measure.
By Lemmas \ref{lem5} and \ref{lem4},
the map $t\mapsto \mu^{(n)}_t$ can be extended by weak continuity
to all times $t\in\R$.
We can now redefine
\begin{equation}
P^{u_n}(t,x)\doteq \int \frac 12 e^{-|x-y|}\,d\mu_t^{(n)}(y)\,,\qquad\qquad
P^u(t,x)\doteq \int \frac 12 e^{-|x-y|}\,d\mu_t(y)\,.
\label{ch1-approxconv}
\end{equation}
where $\mu_t$ is the weak limit of the measures $\mu^{(n)}_t$.
Because of the convergence $J\big(u_n(t),\,u(t)\big)\to 0$,
by Lemma \ref{lem4} the map $t\mapsto \mu_t$ is well defined
and continuous w.r.t.~the weak topology of measures.
Using again
Lemma \ref{lem4}, we can take the limit of (\ref{ch1-evolequinl2}) as $n\to\infty$,
and obtain the identity (\ref{ch1-diffl2}), for every $t\in \R$
and $P=P^u$ defined by (\ref{ch1-approxconv}).

For each $n$, the total energy
$\mu^{(n)}_t\big(]0,1]\big)=E^{\bar u_n}$
is constant in time and converges to $E^{\bar u}$ as $n\to\infty$.
Therefore we also have
$$
\mu_t\big(]0,1]\big)=E^{\bar u}\doteq \int_0^1\big[ \bar u^2(x)+
\bar u_x^2(x)\big]\,dx\qquad\qquad t\in\R\,.
$$

To complete the proof of Theorem \ref{theo1}, it now only remains to prove
that the measure $\mu_t$ is absolutely continuous with density
$$
u^2(t,\cdot)+\frac 12  u^2_x(t,\cdot)
$$ 
w.r.t.~Lebesgue measure, for a.e.~time $t\in\R$.

\begin{figure}[ht]
\centerline{
\includegraphics[width=9cm, height=5cm]{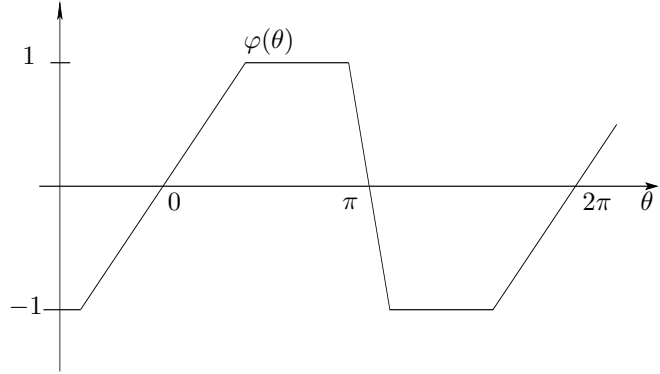}}
\caption{Definition of function $\vphi$\label{varfi}}
\end{figure}

In this direction, we recall that, by (\ref{ch1-uxquadro}), the function $w\doteq
u^2_{n,x}$
satisfies the linear transport equation with source
$$
w_t+(uw)_x=(u_n^2-P^{u_n})u_{n,x}\,.
$$
Moreover, along any characteristic curve
$t\mapsto\xi(t)$
by (\ref{ch1-nonlocalequation}) one has
\begin{equation}
\frac d{dt}\Big[2\arctan u_{n,x}\big(t,\xi(t,x)\big)\Big]
~=~\frac 2{1+u_{n,x}^2}\,
\left[ u_n^2-\frac {u_{n,x}^2}2 -P^{u_n}\right]~\leq~-\frac 12\,,
\label{ch1-decrarct}
\end{equation}
whenever $u_{n,x}^2$ is sufficiently large.
For $\ve>0$ small, consider the piecewise affine, $2\pi$-periodic function
(see fig. \ref{varfi})
$$\vp(\theta)=\left\{
\begin{array}{ccl}
 \theta &\quad&\hbox{if}\quad 0\leq\theta\leq 1\,
\\
1&\qquad&\hbox{if}\quad1\leq\theta\leq \pi-\ve\,,\\
(\pi-\theta)/\ve
&\quad&\hbox{if}\quad \pi-\ve\leq\theta\leq\pi+\ve\,,\\
-1&\quad&\hbox{if}\quad \pi+\ve\leq\theta\leq 2\pi-1\,,\\
\theta-2\pi&\quad&\hbox{if}\quad 2\pi-1\leq\theta\leq 2\pi\,.
\end{array}
\right.
$$
\psfrag{-u}{$-1$}
\psfrag{u}{$1$}
\psfrag{ko}{$0$}
\psfrag{pi}{$\pi$}
\psfrag{2pi}{$2\pi$}
\psfrag{phi}{$\vphi(\theta)$}
\psfrag{theta}{$\theta$}
\vskip 10pt
\vskip 10pt
\n and define
$$\beta_n(t)=\int_0^1 \vp\big( 2\arctan u_{n,x}(t,x)\big)\, u_{n,x}^2(t,x)
\,dx\,.$$
By (\ref{ch1-nonlocalequation}) and (\ref{ch1-decrarct}) we now have
\begin{equation}
\frac d{dt}\beta_n(t)\geq \frac 1{ 4\ve}\int_{\{2\arctan u_{n,x}
\in [\pi-\ve,
\pi]\cup [-\pi, -\pi+\ve]\}} u_{n,x}^2(t,x)\,dx -C\cdot \int_{\{2\arctan u_{n,x}
\in[-1,1]\}}u_{n,x}^2\,dx
\label{ch1-s6}
\end{equation}
for some constant $C$, independent of $\ve,n$.
Since all functions $\beta_n$ remain uniformly bounded, by (\ref{ch1-s6}) for any
time interval $[\tau,\tau']$
we obtain
\begin{equation}
\int_\tau^{\tau'}\int_{\{2\arctan u_{n,x}
\in [\pi-\ve,
\pi]\cup [-\pi, -\pi+\ve]\}} u_{n,x}^2(t,x)\,dx dt\leq
\ve C'\cdot \big(1+\tau'-\tau)\,,
\label{ch1-s7}
\end{equation}
where the constant $C'$ depends only on the $H^1_\per$
norm of the functions $u_n$, hence is uniformly valid for all
$n,\ve$.
Because of (\ref{ch1-s7}), the sequence of functions
$u_{n,x}^2$
is equi-integrable on any domain of the form $[\tau,\tau']\times [0,1]$.
Namely
\begin{equation}
\lim_{\kappa\to\infty}\int_\tau^{\tau'}
\int_{\{x\in[0,1],~ u^2_{n,x} >\kappa\}}
u_{n,x}^2(t,x)\,dxdt=0\,,
\label{ch1-s8}
\end{equation}
uniformly w.r.t.~$n$.
By Lemma \ref{lem4} we already know that $\big\|u_{n,x}^p(t)-u_x^p(t)
\big\|_{L^1_\per}\to 0$
for every fixed time $t$ and $1\leq p <2$.
Thanks to the equi-integrability condition (\ref{ch1-s8}) we now have
$$u^2_{n,x}\to u_x^2\qquad\quad \hbox{in}\quad L^1\big( [\tau,\tau']
\times [0,1]\big)\,.$$
By Fubini's theorem, this implies
$$\lim_{n\to\infty}\int_0^1 u_{n,x}^2(t,x)\,dx=\int_0^1 u_x^2(t,x)\,dx
$$
for a.e.~$t\in [\tau,\tau']$. At every  such time $t$, the measure $\mu_t$
is absolutely continuous and the definition (\ref{ch1-approxconv}) coincides with (\ref{ch1-nonlocalP}).
This completes the proof of Theorem \ref{theo1}.
\endproof
\section{Uniqueness\label{3-7}}
Before proving Theorem \ref{theo3}, we remark that the solution satisfying
all conditions in Theorem \ref{theo1} need not be unique.

\vskip 10pt
\psfrag{0}{$0$}
\psfrag{u(0)}{$u(0)$}
\psfrag{u(t)}{$u(t)$}\psfrag{u(T)}{$u(T)$}
\psfrag{g}{$b$}
\psfrag{q1}{$q_1$}
\psfrag{q2}{$q_2$}
\psfrag{utau}{$u(\tau)$}
\psfrag{utautilde}{$\tilde u(\tau)$}
\begin{figure}[ht]
\centerline{\includegraphics[width=12cm, height= 7cm]{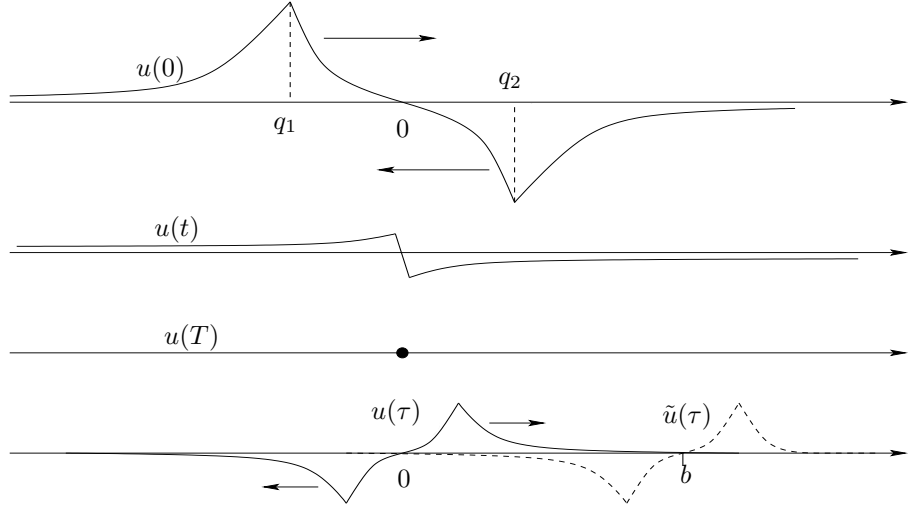}}
\caption{Two solutions of the peakon-antipeakon interaction\label{pap}}
\end{figure}
\vskip 10pt
\begin{example}
\label{ch1-ex}
\rm Let $u=u(t,x)$ be a solution containing exactly
two peakons of opposite strengths $p_1(t)=-p_2(t)$, located
at points $q_1(t)=-q_2(t)$ (see fig. \ref{pap}).  We assume that initially
$p_1(0)>0\,q_1(0)<0$.  At a finite time $T>0$, the two peakons
interact at the origin.  In particular, as $t\to T^-$ there holds
$$p_1(t)\to\infty,\qquad p_2(t)\to -\infty,\qquad  q_1(t)\to 0,\qquad
q_2(t)\to 0\,.$$
Moreover, $\big\|u(t)\big\|_{L^\infty}\to 0$, while
the measure $\mu_t$ approaches a Dirac mass at the origin.
We now have various ways to extend the solution
beyond time $T$:
$$
\tilde u(\tau,x)\doteq 0\,,
$$
\begin{equation}
u(\tau,x)=-u(\tau-T,\, -x)\,,
\label{ch1-dispari}
\end{equation}
Clearly, $\tilde u$ dissipates all the energy, and does
not satisfy the identity (\ref{ch1-energy}).  The function $u$
in (\ref{ch1-dispari}) is the one constructed by our algorithm in Section \ref{3-2}.
However, there are infinitely many other solutions that
still satisfy (\ref{ch1-energy}), for example
$$
\tilde u(\tau,x)=u(\tau, x-b)
$$
where $u$ is as in (\ref{ch1-dispari}) and $b\not= 0$.
The additional condition
in Theorem \ref{theo3} rules out all of them, because
as $\tau\to T+$, the corresponding measures $\tilde \mu_\tau$
approach a Dirac mass at the point $x=b$, not at the origin.
\end{example}
We can now give a proof of Theorem \ref{theo3}.
As a first step, we extend our distance $J$
to a larger domain $\D$, consisting of
couples $(u,\mu)$, where $u\in H^1_\per$ and $\mu$ is a
positive (spatially periodic) measure whose absolutely
continuous part has density $u^2+u_x^2$ w.r.t.~Lebesgue measure.
This extension is achieved by continuity:
$$
J\big( (u,\mu),\, (\tilde u,\tilde \mu)\big)\doteq
\liminf_{n\to\infty} J(u_n,\tilde u_n)
$$
where the infimum is taken over all couple of sequences
$(u_n,\tilde u_n)_{n\geq 1}$ such that
$$\|u_n-u\|_{L^\infty}\to 0\,,\qquad\qquad
\|\tilde u_n-\tilde u\|_{L^\infty}\to 0\,,$$
$$u^2_{n,x}\wto \mu\,,\qquad\qquad \tilde u^2_{n,x}\wto\tilde \mu\,.$$
We observe that the flow $\Phi$ constructed in Theorem \ref{theo2} can be
continuously extended to a locally Lipschitz continuous group
of transformations on the domain $\D$.

Now let $t\mapsto \tilde u(t)$ be a solution of the Cauchy problem
(\ref{ch1-equation}), (\ref{ch1-initialcond}), satisfying all the required conditions.
In particular,  the map $t\mapsto \big(\tilde u(t),\tilde \mu_t\big)$
is Lipschitz continuous w.r.t.~the distance $J$, with values
in the domain $\D$.

Calling $t\mapsto \big(\tilde u(t),\tilde \mu_t\big)
\doteq \Phi_t(\bar u, \bar u^2_x)$
the unique solution of the Cauchy problem obtained
as limit of multi-peakon approximations, we need to
show that $\tilde u(t)=u(t)$
for all $t$.  To fix the ideas, let $t>0$.
By the Lipschitz continuity of the flow,
we can use the error estimate
\begin{equation}
J\Big(\big(\tilde u(t),\tilde \mu_t\big)\,,~\big(u(t), \mu_t\big)
\Big)\leq e^{C_2 t} \int_0^t \liminf_{h\to 0}
\frac 1h\cdot J\Big( \big(\tilde u(\tau+h),\tilde \mu_{\tau+h}\big)\,,~
\Phi_h\big(\tilde u(\tau),\tilde \mu_\tau\big)\Big)\,d\tau
\label{ch1-bressanest}
\end{equation}
For a proof of (\ref{ch1-bressanest}), see \cite[pp. 25--27]{B}.
The conditions stated in Theorem \ref{theo1} now imply that, at almost
every time $\tau$,
the measure $\tilde\mu_t$ is absolutely continuous
and the integrand in (\ref{ch1-bressanest}) vanishes.  Therefore $\tilde u(t)=u(t)$
for all $t$.
\endproof
\v
We can now prove that, in multi-peakon solutions,
interactions involving exactly two peakons are the only possible ones.
\begin{cor}\label{cormax2peak} Let $t\mapsto u(t,\cdot)$ be a multi-peakon solution
of the form (\ref{ch1-mpeak}), which remains regular
on the open interval $]0,T[\,$.
Assume that at time $T>0$ an interaction occurs, say among the first
$k$ peakons, so that
$$\lim_{t\to T-} q_i(t)=\bar q\qquad\qquad i=1,\ldots,k\,.$$
Then $k=2$.
\end{cor}
\begin{proof}
We first observe that the Camassa-Holm equations
(\ref{ch1-equation}) are time reversible.  In particular, our
proof of Theorem \ref{theo3} shows that the
solution to a Cauchy problem is unique both forward and backward in time.

Now consider the data
$\big(u(T),\,\mu_T\big)\in\D$, where
$\mu_T$ is the weak limit of the measures $\mu_t$ having density
$u^2(t)+u^2_x(t)$ w.r.t.~Lebesgue measure, as $t\to T-$.
By the analysis in Section \ref{3-2}, we can construct a backward solution
of this Cauchy problem
in terms of exactly two incoming peakons.
By uniqueness, this must coincide with the given solution $u(\cdot)$
for all $t\in [0,T]$.
\end{proof}

\part{The discrete Boltzmann equation}
\chapter{Symmetry groups of differential equations\label{symmgenerator}}
In this chapter we introduce the theory of the symmetry groups applied to differential equation, which is a tool that will fits in the study of evolutionary equations. The goal of this Chapter is to develop a useful method that will explicitly determine the symmetry group for the system of discrete Boltzmann equation, which will be the starting point of the discussion of Chapter \ref{chapblowup} for the blow-up issue. The key point is to transform the equation which has an asymptotic blow-up at a time $T$ into an equation, related to a rescaling of the first equation, which approach a steady state as $\tau$ goes to infinitive (see Section \ref{symblowup}). 
\section{Group and differential equations}
The symmetry group of a system of differential equations is the largest local group of transformations acting on the independent and dependent variables with the property that it transforms solution of the system to other solution.  In the first part of this section we review a general computational method for (almost) any given system of differential equations. For more information about the application of group theory to the differential equation, we refer the reader to \cite{O,Ov}.

We start recalling some useful definition in the abstract theory of Transformation Group.  
\begin{dhef}
\label{chap:4def}
Let $\mathcal M$ be a smooth manifold. A \emph{local group of transformations} acting on $\mathcal M$ given by a (local) Lie group $G$ is the couple $(\mathscr U,\Psi)$ where 
\begin{itemize}
\item $\mathscr U$ is an open subset $\{\iota\}\times \mathcal M\subset \mathscr U \subset G\times \mathcal M$
\item $\Psi$ a smooth map $\Psi:\mathscr U\to \mathcal M$  
\end{itemize}
satisfy the properties
\begin{itemize}
\item[(a)] If $(h,x)\in \mathscr U$, $(g,\Psi(h,x))\in \mathscr U$ and $(g\cdot h,x)\in \mathscr U$ then
$$
\Psi(g,\Psi(h,x))=\Psi(g\cdot h, x).
$$
\item[(b)] For all $x\in \mathcal M$, $\Psi(\iota,x)=x$.
\item[(c)] If $(g,x)\in \mathscr U$ then $(g^{-1}, \Psi(g,x))\in \mathscr U$ and $\Psi(g^{-1}, \Psi(g,x))=x$.
\end{itemize}
For brevity, when it does not make confusion, we denote $\Psi(g,x)$ by $g\cdot x$. 
\end{dhef}
\begin{dhef}
Let $G$ be a local group of transformation acting on a manifold $\mathcal M$. A subset $\mathcal S \subset \mathcal M$ is called \emph{$G-$invariant}, and $G$ is called a \emph{symmetry group} of $\mathcal S$ if whenever $x\in \mathcal S$ and $g\in G$ are such that $g\cdot x$ is defined, then $g\cdot x\in \mathcal S$.
\end{dhef}
\begin{oss}
\rm
In our applications, as far as the differential equation is concerned, the subset $\mathcal S$ will be usually the graph of the solution of the differential equation 
\begin{equation}
\Phi(x,u, \dots, D^\alpha u,\dots )=0,
\label{eqdifgen}
\end{equation}
where i.e. set of solutions determined by the common zeros of collection of smooth functions $\Phi=(\Phi_1, \dots, \Phi_l)$, where 
$\Phi_j=\Phi_j(x,u, \dots, p^\alpha, \dots)$ depends on the variables $x$ and the unknowns and their derivatives $u, D^\alpha u, \dots$, and where, for every multi-index $\alpha=(\alpha_1,\dots,\alpha_m)$, $D^\alpha$ indicates the differential operator
$$
D^\alpha = \left(\frac{\partial}{\partial x_1}\right)^{\alpha_1}\cdots\left(\frac{\partial}{\partial x_m}\right)^{\alpha_m}
$$
In this context, it is thus useful to introduce the graph of a function  $u:\Omega \to R^n$ defined on a open set $\Omega \subset \R^m$
$$\Gamma_u=\{(x,u(x)): x\in \Omega \}$$
which is a smooth submanifold of $\R^m\times \R^n$.  
The action of a given transformation $g\in G$ maps the graph $\Gamma_u$ into the subset  $g\cdot \Gamma_u=\{(\tilde x,\tilde u)=g\cdot(x,u) :(x,u)\in \Gamma_u \}$ which is not necessarily the graph of a function $\tilde u$. However, since $G$ acts smoothly and the identity $\iota	\in G$ leaves $\Gamma_u$ unchanged, by restricting the domain $\Omega$ for every $g\in G$ near the identity the transformation $g\cdot \Gamma_u$ is the graph $\Gamma_{\tilde u}$ of a function $\tilde u$.  
\end{oss}
As an example of action on the graph of a function, let consider a vector field $\mathbf v(\xi,\eta):\mathcal M=\R^{m+n}\to \R^{m+n}$ which can be seen in local coordinate $(x_1,\dots,x_m,u_1,\dots,u_n)$ as 
$$
\mathbf v=\sum_{i=1}^{m}\xi_i(x,u)\frac{\partial}{\partial {x_i} }+\sum_{j=1}^{n}\eta_j(x,u)\frac{\partial}{\partial {u_j} }
$$
it acts on a smooth scalar function $\phi:\R^{m+n}\to \R$ as a derivation
\begin{equation}
\mathbf v\bullet \phi(x,u)= \lim_{\eps\to 0}\frac{\phi(x+\eps \xi,u+\eps\eta)-\phi(x,u)}{\eps}.
\label{derwrtv}
\end{equation}
The most important operation on vector fields is their Lie bracket or commutator. Whenever we think that two vector field $\mathbf v,\mathbf w$ act as a derivation, their Lie bracket $[\mathbf v,\mathbf w]$ is the unique vector field satisfying
\begin{equation}
\label{Liebracket}
[\mathbf v,\mathbf w]\bullet \phi = \mathbf v\bullet (\mathbf w \bullet \phi) - \mathbf w\bullet (\mathbf v \bullet \phi)\qquad \mbox{for all smooth functions $\phi$}
\end{equation}
The integral curve of the vector field $\mathbf v$ is a smooth parametrized curve $P(\theta)=(x,u)$ whose tangent vector at any point coincides with the value of $\mathbf v$ at the same point: 
$$
\frac d{d\theta} P= \mathbf v(P).
$$
Starting from a given initial data $P(0)=\bar P=(\bar x,\bar u)$ the corresponding integral curve is often denoted by the suggestive exponential notation
$$
\exp(\theta \mathbf v)\bar P.
$$
From the existence and uniqueness of solution to systems of ordinary differential equations we easily obtain the semigroup property for the \emph{flow} generated by $\mathbf v$:
\begin{eqnarray}
&&\exp(0 \mathbf v)\bar P=\bar P
\\
&&\exp[(\theta_1+\theta_2) \mathbf v]\bar P=\exp(\theta_1 \mathbf v)[\exp(\theta_2\mathbf v)\bar P].
\end{eqnarray}
From these formulas, compared with the property (a)-(b) of definition \ref{chap:4def}, we see that the flow generated by a vector field is the same as a local action of the Lie group $\R$ on the manifold $\mathcal M$ which is called a \emph{one-parameter group of transformations}. The vector field $\mathbf v$ is called the \emph{infinitesimal generator} of the action. $\mathbf v$ is also called an \emph{infinitesimal symmetry generator} for (\ref{eqdifgen}) if the map $P\mapsto \exp(\theta \mathbf v)P$ transforms the graph of a solution $u$ into the graph of another solution. 
\begin{oss}
\label{transgraph}
\rm
Recall that if $\theta$ is sufficiently small, there exists a neighborhood $V$ such that the set
$$
\Gamma_\theta=\{\exp(\theta \mathbf v)(x,u): (x,u)\in \Gamma_u\}
$$
coincides on $V$ with the graph $\Gamma_{u^\theta}$ of a smooth function $u^\theta$. By using a Taylor expansion, the flow $\exp(\theta \mathbf v)$ maps the point $(x,u)$ into the point 
$$
\exp(\theta \mathbf v)(x,u)=(x+\theta \xi(x,u)+o(\theta),u+\theta \eta(x,u) + o(\theta))
$$
therefore, the differentiation w.r.t. $\theta$ at the origin yields the useful formula
\begin{equation}
\label{condsol}
\left.\frac d{d\theta} u^\theta(x)\right|_{\theta=0} =\eta - \nabla u(x)\cdot \xi.
\end{equation}
\end{oss}
The previous formula gives a necessary condition in order to prove that a particular vector field $\mathbf v =(\xi, \eta)$ is an infinitesimal symmetry generator for the differential equation (\ref{eqdifgen}). If the function $u^\theta$ another solution to this equation, then $\Phi(x,u^\theta, \dots,D^\alpha u^\theta ,\dots)=0$. Differentiating w.r.t. $\theta$ in $0$ we get
\begin{equation}
\label{condeq}
\left. \sum_{\alpha} \frac{\partial \Phi}{\partial p^\alpha}\cdot \frac{d}{d\theta}(D^\alpha u^\theta)\right|_{\theta=0}
=\sum_{\alpha} \frac{\partial \Phi}{\partial p^\alpha}\cdot D^\alpha\left(\eta-\nabla u(x)\cdot \xi\right)=0
\end{equation}
In the following we shall prove that the previous condition is sufficient in order to construct an infinitesimal symmetry generator $\mathbf v$.  
\begin{prop}
\label{propinvfun}
Let $G$ be a connected group of transformation acting on the manifold $\mathcal M$. A smooth real-valued function $\zeta:\mathcal M \to \R$ is an invariant function for $G$ if and only if 
\begin{equation}
\mathbf v\bullet \zeta=0\qquad \mbox{for all $x\in \mathcal M$}
\label{condinvfunc}
\end{equation}
and every infinitesimal generator $\mathbf v$ of $G$.
\end{prop}
\begin{proof}
Suppose that $\zeta$ is an invariant function for $G$. According to (\ref{derwrtv}), if $x\in \mathcal M$ 
$$
\frac d{d\theta} \zeta(\exp(\theta \mathbf v)x)=\mathbf v\bullet \zeta [\exp(\theta \mathbf v)x]
$$ 
since $\zeta$ is invariant, setting $\theta =0$ it proves the necessity of (\ref{condinvfunc}).
Conversely, if (\ref{condinvfunc}) holds then $\zeta(\exp(\theta \mathbf v)x)$ is a constant for the connected subgroup $\{\exp(\theta \mathbf v)\}$ of $G_x\doteq \{ g\in G: \mbox{$g\cdot x$ is defined} \}$. But by the properties of the Lie group, every element of $G_x$ can be written as a finite product $g=\exp (\theta^1 \mathbf v_{i_1})\cdots\exp (\theta^k \mathbf v_{i_k})$ for some infinitesimal generator $v_i$ of $G$, hence $\zeta(g\cdot x)=\zeta(x)$ for all $g\in G_x$.
\end{proof}

In a similar way we can prove the following theorem which gives an infinitesimal criterion of invariance for a general equation
$$
\Phi(x)=0 \qquad x\in \mathcal M
$$
that will be useful whenever we are concerning a differential equation  
$$
\Phi(x,u, \dots, D^\alpha u, \dots)=0.
$$ 
\begin{theorem}
\label{thalg}
Let $G$ be a connected local Lie group of transformations acting on a $p-$dimensional manifold $\mathcal M$. Let $\Phi: \mathcal M\to \R^l$, $l\leq p$, define a system of equations
\begin{equation}
\Phi_\nu(x)=0 \qquad \nu=1,\dots, l
\label{sysalg}
\end{equation}
and assume that the system has maximal rank at every solution $x$ of the system, namely
$$
\rank \left(
\begin{array}{ccc}
\frac{\partial \Phi_1}{\partial x^1}(x) &
\cdots &
\frac{\partial \Phi_1}{\partial x^m}(x) 
\\
\vdots &\ddots&\vdots
\\
\frac{\partial \Phi_l}{\partial x^1}(x) &
\cdots &
\frac{\partial \Phi_l}{\partial x^m}(x)
\end{array}
\right)=l
\qquad \Phi(x)=0.
$$
Then $G$ is a symmetry group of the system if and only if 
\begin{equation}
\label{condinvsys}
\mathbf v\bullet \Phi_\nu(x)=0 \qquad \nu=1,\dots, l,\quad \Phi(x)=0
\end{equation}
for every infinitesimal generator $\mathbf v$ of $G$.
\end{theorem} 
\begin{proof}
Let $x_0$ be a solution of the system (\ref{sysalg}). As in Proposition \ref{propinvfun},  the necessary condition follows from differentiating w.r.t. $\eps$ the identities
$$
\Phi_\nu(\exp(\eps\mathbf v)x_0)=0
$$
and setting $\eps=0$. Conversely, by using the maximal rank condition, we can choose local coordinates $y=(y^1,\dots,y^m)$ such that $x_0=0$ and $\Phi$ has the simple form $\Phi(y)=(y_1,\dots,y^l)$. Let $\mathbf v$ be an infinitesimal generator  of $G$, which can be expressed in the new coordinates as
$$
\mathbf v=\xi^1(y)\frac{\partial}{\partial y^1}+\cdots+\xi^m(y)\frac{\partial}{\partial y^m}.
$$
The condition (\ref{condinvsys}) turns to be 
$$
\mathbf v(y^\nu)=\xi^\nu(y)=0 \qquad \mbox{for all $\nu = 1\dots l$}
$$
whenever $y^1=\cdots =y^l=0$. Since the flow $\phi(\theta)=\exp(\theta \mathbf v)x_0$ satisfies the system of ODE
$$
\left\{
\begin{array}{l}
\frac{d}{d\theta}\phi^i=\xi^i(\phi(\theta))
\\
\phi^i(0)=0
\end{array}
\right.\qquad i=1\dots m,
$$ 
the uniqueness of the solution yields to conclude that $\phi^i(\theta)=0$ for $\theta$ sufficiently small. $\exp(\theta \mathbf v)x_0$ is thus again a solution to $\Phi(x)=0$. As in Proposition \ref{propinvfun}, by the properties of the connected local Lie group $G$ we gain the result.  

\end{proof}
The previous theorem can be adapted for the differential equation in order to get sufficient condition for obtain an infinitesimal symmetry generator. The equation (\ref{eqdifgen}) contains not only the unknowns $u$ but also its derivatives, so in order to use Theorem \ref{thalg} we can think that the solution is a point which contains all of these functions. To do this we need to prolong the basic space representing the independent and dependent variables under consideration to a space which also represents the various partial derivatives occurring in the system.

If we consider function $u:\Omega\subset\R^m\to \R^n$, the number of derivatives of order $k$ is 
$$
p_k=
n\cdot \left( 
\begin{array}{c}
m+k-1
\\
k
\end{array}
\right)
$$ 
If $N$ is the maximum order of the derivatives involved in the differential equation (\ref{eqdifgen}), we introduce thus the \emph{$N-$th jet space} $\Omega\times U^{[N]}\doteq \Omega\times\R^m\times \R^{p_1}\times\dots \times\R^{p_N}$, whose coordinates represent all the derivatives of the function $u$ from $0$ to $N$. If $u$ is a function whose graph lies in a manifold $\mathcal M\subset \Omega\times \R^m$, we define its prolongation
$$
\pr^{(N)}u=(u,(D^{\alpha_1} u)_{|\alpha_1|=p_1}, \dots,(D^{\alpha_N} u)_{|\alpha_N|=p_N})
$$  
whose graph lies in the $N-$th jet space  $\mathcal M^{(N)}\doteq\mathcal M\times\R^{p_1}\times\dots \times\R^{p_N}$. From this point of view, a smooth solution of the given system (\ref{eqdifgen}) is a function $u(x)$ such that 
$$
\Phi_\nu(x,\pr^{(N)}u(x))=0 \qquad \nu=1,\dots, l,
$$ 
whenever $x$ lies in the domain of $u$. It means that the graph of the prolongation of $u$ must lie entirely within the subvariety of the zeroes of the system. 

Now suppose that $G$ is a local group of transformations acting on $\mathcal M$. The prolongation of $G$ as a local action group on the prolonged manifold $\mathcal M^{(N)}$ is defined so that if $g\in G$, $\pr^{(N)} g$ transforms the derivatives of $u$ into the corresponding derivatives of the transformed $g\cdot (x,u)$. To evaluate the action of the prolonged $\pr^{(N)} g$ on a couple $(x_0,u_0^{(N)})$ we simply choose a particular function $f$ whose derivatives agree, up to $N-$th order, to the point $(x_0,u_0^{(N)})$, apply the action $g$ to $f$ and then prolong $g\cdot f$. 
Last, we have to define also the prolongation of a vector field up to the order $N$. It follows by viewing it as the infinitesimal generator of the corresponding action group $\pr^{(N)}[\exp(\theta \mathbf v)]$:
$$
\pr^{(N)}\mathbf v\doteq \left.\frac d{d\theta}\right|_{\theta=0} \pr^{(N)}[\exp(\theta \mathbf v)].
$$
Writing  
$$
\mathbf v = \xi(x,u) \frac \partial{\partial x} + \eta(x,u) \frac \partial{\partial u},
$$
the general formula of such a vector field is given by \cite[Theorem 2.36]{O} and it is the following formal expression
\begin{equation}
\label{defprvf}
\pr^{(N)}\mathbf v = \mathbf v + \sum_{1\leq |\alpha|\leq N} \phi^\alpha(x, u^{(N)}) \frac \partial{\partial u^\alpha} 
\end{equation}
where
$$
\begin{array}{l}
\dis u^\alpha = D^\alpha u,
\\
\dis \phi^\alpha= D^\alpha\left(\eta-\xi\cdot \nabla u\right)+ \xi\cdot \nabla u^\alpha.
\end{array}
$$
For future use, the prolongation of the Lie bracket vector field is
\begin{equation}
\label{prolLie}
\pr^{(N)}[\mathbf v,\mathbf w]=[\pr^{(N)}\mathbf v,\pr^{(N)}\mathbf w].
\end{equation}
By applying Theorem \ref{thalg} to the equation $\Phi(x,\pr^{(N)}u(x))=0$ we obtain the following theorem, which agree with the formula (\ref{condeq})
\begin{theorem}
\label{theocond}
Suppose
$$
\Phi_\nu(x,\pr^{(N)}u(x))=0 \qquad \nu=1,\dots, l,
$$
is a system of differential equations of maximal rank defined over $\mathcal M\subset \omega\times \R^{n}$. If $G$ is a local group of transformations acting on $\mathcal M$, and 
$$
\pr^{(N)}\mathbf v\bullet\Phi_\nu(x,\pr^{(N)}u)=0,\qquad \mbox{for all $\nu=1,\dots,l$,} 
$$
for every infinitesimal generator $\mathbf v$ of $G$, then $G$ is a symmetry group of the system.\endproof
\end{theorem}
\section{Symmetries and blow-up}
\label{symblowup}
In this section we enter in deep detail of the application of local symmetry groups to the blow-up issue. For a more complete description, the reader can see \cite{B4}. Our goal is to apply the theory previously developed in order to choose an appropriate infinitesimal symmetry generator which describes the asymptotic  behaviour of a solution which has blow-up in finite time.

We shall consider an evolution problem on a Banach space $(\|\cdot\|,E)$
\begin{equation}
\label{ODEvf}
\dot x = f(x).
\end{equation}
Assume that there exists trajectories that blow-up in finite time
$$
\lim_{t\to T^-} \|x(t)\|=+\infty.
$$
Suppose that there exists a second vector field $g$ such that
\begin{itemize}
\item Trajectories of $\dot y = f(y)-g(y)$ do not blow-up. Instead, they approach a steady state $\bar y$ as time  goes to infinity.
\item There is an explicit computable transformation that maps a trajectory $s\mapsto y(s)$, $s\in [s_0,+\infty[$ into a trajectory $t\mapsto x(t)$, $t\in [t_0,T[$.
\end{itemize}
In this case we could first accurately study the asymptotic behaviour of $y(s)$ as $s\to +\infty$, and then recover information on the behaviour of the blowing-up solution $x(t)$. To implement this approach, it is clear that the auxiliary vector field $g$ must be carefully selected. 

According to the theory developed in the previous section, we shall look for condition on
\begin{equation}
\frac d{d\theta} x^\theta = g(x^\theta)
\label{ODEtrans}
\end{equation}
in order to have that $u^\theta\doteq \exp(\theta g)u$ gives another solution to the equation (\ref{ODEvf}), provided that we have existence and uniqueness of solution to the Cauchy problem (\ref{ODEtrans}).  

As the following analysis will show, the crucial assumption on $g$ is the relation
\begin{equation}
f+[f,g]\equiv 0.
\label{idLie}
\end{equation}
\begin{lemma}
Let $g$ be a vector field such that  \ref{idLie} holds. Then the vector field
\begin{equation}
\label{newvfield}
\mathbf v \doteq -t\frac \partial{\partial t} + g(x) \frac{\partial}{\partial x}
\end{equation}
is an infinitesimal generator of symmetry group. In other words, if $t \mapsto x(t)$ is a solution of (\ref{ODEvf}) then 
$$
\left\{\exp(\theta \mathbf v) (t,x(t))\,:\, t\in I\right\} 
$$
is the graph of another solution
$$
x^\theta(t) \doteq \exp (\theta g)(x(e^\theta t)).
$$
\end{lemma}
\begin{proof}
To check that $x^\theta$ is indeed a solution, we write 
$$
\frac{dx^\theta}{dt}=e^\theta {\rm \,Jac}\left[\exp(\theta g)\right] f\left(\exp(-\theta g)x^\theta\right).$$
We claim that the right hand side of the previous formula coincides with $f(x^\theta)$, and this can be done by proving 
\begin{equation}
\label{idgef}
v(\theta)\doteq {\rm \,Jac}\left[\exp(\theta g)\right] f\left(\exp(-\theta g)y\right)=e^{-\theta}f(y) \qquad \mbox{for all $y\in E$.}
\end{equation}
Trivially, $v(0)=f(y)$. As far as the derivative of $v$ is concerned, let use the Lie bracket property (see \cite{doC})
$$
[f,g]= \lim_{\eps\to 0} \frac{{\rm \,Jac}\left[\exp(\eps g)\right]g(\exp(\eps f))-g}\eps.
$$
and the hypothesis (\ref{idLie})
$$ 
\begin{array}{rl}
\frac{dv(\theta)}{d\theta}&=
\lim\limits_{\eps\to 0}\frac{{\rm \,Jac}\left[\exp[(\eps+\theta)] g)\right]f(\exp[(-\theta-\eps)g]y)-{\rm \,Jac}\left[\exp(\theta g)\right]f(\exp(-\theta g)y) }{\eps}
\\
&=-\dis {\rm \, Jac}\left[\exp(\theta g)\right]\cdot [g,f] \exp(-\theta g)y =-v(\theta).
\end{array}
$$
The function $v(\theta)$ is thus
$$
v(\theta)=e^{-\theta}v(0)=e^{-\theta}f(y).
$$
\end{proof}
\begin{theorem} {(\bf Blow-up rescaling)}
\label{blowuprescaling}
Consider two vector fields $f,g$ satisfying (\ref{idLie}). Let $y:[0,\infty[\mapsto E$ be a solution to 
$$
\dot y(s)=f(y(s))-g(y(s)).
$$
Then the function $x:[0,1[\mapsto E$ defined by 
\begin{equation}
x(t)\doteq \exp(s g)(y(s))\qquad s\doteq \ln\frac 1{1-t}
\label{ressol}
\end{equation}
is a solution of (\ref{ODEvf}).
\end{theorem}
\begin{proof}
\psfrag{ys}{$y(s)$}
\psfrag{xt}{$x(t)$}
\psfrag{yb}{$\bar y$}
\psfrag{expsg}{$\exp(sg)$}
\begin{figure}[ht]
\centerline{
\includegraphics[width=7 cm]{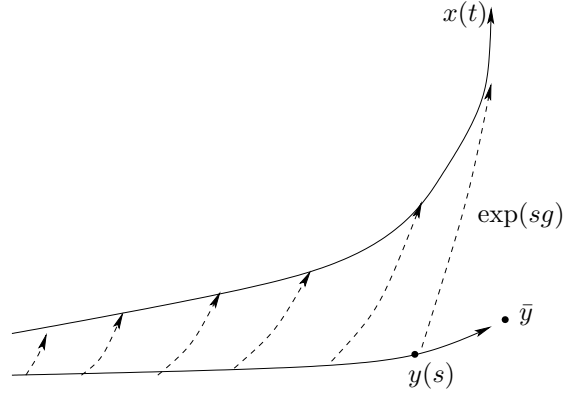}}
\caption{Connection between $y(\cdot)$ and $x(\cdot)$}
\end{figure}
Notice that
$$
g(\exp(sg)y)=\left.\frac d{d\eps} \exp(\eps g)(\exp (sg)y)\right|_{\eps= 0}={\rm\, Jac}[\exp(sg)]g(y)
$$
Let compute the derivative w.r.t. $s$ of the function $\exp(sg)y(s)$. By the previous identity and by (\ref{idgef}) we have 
$$
\begin{array}{rl}
\dis \frac{d}{ds} \exp(sg)y(s) = &
\dis g(\exp(sg)y(s))+{\rm\, Jac}[\exp(sg)]\dot y(s)
\\
= &
\dis g(\exp(sg)y(s))+{\rm\, Jac}[\exp(sg)]f(y(s))-{\rm\, Jac}[\exp(sg)]g(y(s))
\\
= &
{\rm\, Jac}[\exp(sg)]f(\exp(-sg)x)
\\
= &e^{-s}f(x)
\end{array}
$$ 
now we can check that (\ref{ressol}) is a solution to (\ref{ODEvf}).
$$
\frac d{dt}x(t)=\frac{d}{ds} \exp(sg)y(s)\cdot \frac{ds}{dt}=e^{-s}f(x)\cdot e^s = f(x).
$$
\end{proof}
We shall apply the above theory to a special case of partial differential equations. Inside the equation (\ref{eqdifgen}) we highlight the time variable $x\leadsto (t,x)$, obtaining the evolution equation
\begin{equation}
\label{evoleq}
F(x,\pr^{(N)}u)-u_t\doteq F(x,u,\dots,D^{\alpha} u, \dots) - u_t = 0
\end{equation}
where the derivatives on $D^\alpha$ involves only derivatives with respect to the spatial variable$~x$. In the following, we indicate with $[u]$ the prolonged function $u$ in the space $\mathcal M\times U^{(N)}$. Thus, the total derivative of a function
$F[u(t)]$ takes the form $\nabla F\cdot [u_t]$, where $[u_t]=(u_t,\dots,D^\alpha u_t,\dots)$.

The main tool is Theorem \ref{theocond}, which gives characterization on the symmetric vector field. From this theorem we obtain
\begin{theorem}
\label{generatorforPDE}
An evolutionary vector field 
$$
\mathbf v\doteq -t\frac \partial{\partial t}+\xi(x,u) \frac\partial{\partial x}+\phi(x,u) \frac\partial{\partial u}
$$
is an infinitesimal symmetry generator for (\ref{evoleq}) if and only if for every function $u$ one has 
\begin{equation}
\label{bracketPDE}
(F+[F,G])[u]=0
\end{equation}
where 
$$
F[u]\doteq F(x,\pr^{(N)}u),\qquad G[u]\doteq -\xi(x,u)\nabla u+\phi(x,u),
$$
\begin{equation}
\label{defFGbracket}
[F,G][u]\doteq \lim_{\eps\to 0} \frac {G[u+\eps F[u]]-G[u]-F[u+\eps G[u]]+F[u]}\eps.
\end{equation}
\end{theorem}
\begin{proof}
Given any smooth solution $u$, let $\theta$ sufficiently small and $u^\theta$ be the transformed function by the vector field $\mathbf v$, according to Remark \ref{transgraph}. As Theorem \ref{theocond} states, the thesis follows once we prove (\ref{condeq}). By formula (\ref{condsol}), substituting $(\xi,\eta)$ with $((-t,\xi),\phi )$
we have
$$
\left.\frac{d u^\theta}{d\theta}\right|_{\theta =0} =\phi +t u_t-\xi\cdot \nabla u = t u_t + G[u].
$$
Thus, formula (\ref{condeq}) becomes
$$
\begin{array}{rl}
\dis
\left.
\frac{d\Phi[u^\theta]}{d\theta}\right|_{\theta = 0} &= 
\dis
\frac{d\F[u]}{d\theta}-\frac{d u^\theta_t}{d\theta}
=
 \nabla F[u]\cdot[t u_t+G[u]]-u_t-t u_{tt}-\nabla G\cdot[u]_t
\\
&=
\dis
 \nabla F[u]\cdot[t u_t+G[u]]-F[u]-t\nabla F[u]\cdot[u_t] -\nabla G\cdot[u]_t
\\
&=
\dis 
-F[u]+\nabla F[u]\cdot [G[u]]-\nabla G[u]\cdot[F[u]].
\dis

\end{array}
$$
thus
$$
\left. \frac{d\Phi[u^\theta]}{d\theta}\right|_{\theta = 0} = 0 
$$
if and only if $[F,G][u]=-F$.
\end{proof}
\section{A group of symmetry for the discrete Boltzmann equation\label{symm}}
The work plan presented in the previous two sections fits in the study of the blow-up rate of discrete Boltzmann equation. Let consider the system of PDE (see \cite{MP})
\begin{equation}
\label{bzmnncap4}
\partial_t u_i+ \mathbf c_i \nabla_x u_i = \sum_{jk} a_{ijk} u_j u_k \qquad i= 1,\dots,l
\end{equation}
where
\begin{itemize}
\item $(t,x)\in \R\times \R^3$,
\item $\mathbf c_i\in \R^3$ plays the role of velocity,
\item $u_i(t,x)$ is the density of particles having speed $\mathbf c_i$, 
\item $a_{ijk}$ are the coefficients of the quadratic collision term, with $a_{ijk}=0$ if $j=k$.
\end{itemize}
We implement Theorem \ref{generatorforPDE} in order to recover an evolutionary vector field which generates a group of symmetry for the system (\ref{bzmnncap4}). Looking for a vector field of the form
$$
-t\frac \partial{\partial t}+\xi(x)\cdot\nabla_x +\phi(u)\cdot\frac \partial{\partial u}
$$
where
$$
\nabla_x\doteq \left(\frac \partial{\partial x_1},\frac \partial{\partial x_2}, \frac \partial{\partial x_3} \right),
\qquad \frac \partial{\partial u} \doteq \left(\frac \partial{\partial u_1},\frac \partial{\partial u_2},\dots, \frac \partial{\partial u_l} \right)
$$
we have to find condition for which (\ref{bracketPDE}) holds. Note that in this case the two maps $F$ and $G$ are
$$
(F[u])_i\doteq -c_i\cdot \nabla_x + \sum_{j,k} a_{ijk} u_j u_k\qquad (G[u])_i\doteq -\xi(x)\cdot\nabla_x u_i + \phi_i(u) \qquad i=1\dots l.
$$
\begin{theorem}
\label{rescperboltzmann}
Setting
$$
\xi(x)\doteq - x\qquad \phi_i(u)\doteq u_i,
$$
the identity 
$$
[F,G][u]+F[u]=0
$$
holds for every $u$ solution to (\ref{bzmnncap4}).
\end{theorem}
\begin{proof}
Let compute (\ref{defFGbracket}).
\begin{equation}
\label{GI}
\begin{array}{rl}
\dis \frac{(G[u+\eps F[u]]-G[u])_i}\eps
=& 
\dis 
\frac{x\cdot \nabla_x \left(u_i-\eps \mathbf c_i\cdot \nabla_x u_i + \eps \sum_{jk}a_{ijk} u_j u_k\right)}\eps
\\
&\dis 
+\frac{u_i -\eps \mathbf c_i\cdot \nabla_x u_i +\eps\sum_{jk}a_{ijk} u_j u_k - x\cdot \nabla_x u_i - u_i}\eps
\\
=&
\dis 
-x\cdot \nabla_x \left(\mathbf c_i\cdot \nabla_x u_i\right) + (F[u])_i
\\
&
\dis + \sum_{jk}a_{ijk} \left(u_k x\cdot \nabla_x u_j + u_j x\cdot \nabla_x u_k \right)
\end{array}
\end{equation}
\begin{equation}
\label{EFFE}
\begin{array}{rl}
\dis \frac{(F[u+\eps G[u]]-F[u])_i}\eps
=&
\dis
\frac{-\mathbf c_i \cdot \nabla_x u_i - \eps \mathbf c_i \nabla_x(x\cdot \nabla_x u_i)-\eps \mathbf c_i \nabla_x u_i }
\eps
\\
&\dis +\frac{\sum_{jk}a_{ijk}(u_j +\eps x\cdot \nabla_x u_j + \eps u_j)(u_k+\eps x\cdot \nabla_x u_k + \eps u_k)}\eps
\\
&\dis -\frac{-\mathbf c_i \nabla_x u_i + \sum_{jk}a_{ijk}u_j u_k}\eps
\\
=&
\dis
-\mathbf c_i\cdot \nabla_x u_i - \mathbf c_i \nabla_x(x\cdot \nabla_x u_i)
\\
&\dis +\sum_{jk}a_{ijk}(2 u_j u_k + u_j x\cdot \nabla_x u_k + u_k x\cdot \nabla_x u_j) 
+\mathcal O(\eps)
\end{array}
\end{equation}
Note that
$$
x\cdot \nabla_x(\mathbf c_i \cdot \nabla_x u_i)=\sum_j x_j \sum_k (\mathbf c_i)_k u_{x_j x_k}
$$
$$
\mathbf c_i\cdot \nabla_x(x\cdot \nabla_x u_i)=\mathbf c_i\cdot \nabla_x u_i + \sum_k (\mathbf c_i)_k \sum_j x_j  u_{x_j x_k}.
$$
Hence, the limit of the difference of the formulas (\ref{GI}), (\ref{EFFE}) yields
$$
([F,G][u])_i=-2(F[u])_i+ (F[u])_i=-(F[u])_i.
$$
\end{proof}
\begin{cor}
If $u(t,x)=(u_1(t,x),\dots,u_l(t,x))$ is a solution, then 
$$
u^\theta(t,x)=e^\theta u(e^\theta t,e^\theta x)
$$ 
is another solution to the system of PDE (\ref{bzmnncap4}).
Hence, by performing the change of variables
\begin{equation}
\label{systemchange}
\left\{
\begin{array}l
\dis\tau=\ln \frac 1{1-t}
\\
\dis\eta = e^\tau x = \frac x{1-t}
\\
\dis w_i = e^{-\tau}u_i = (1-t)u_i
\end{array}
\right.
\end{equation}
blow-up to (\ref{bzmnncap4}) occurs if a solution to
$$
(w_i)_\tau = -(\mathbf c_i+\eta)\cdot\nabla_\eta w_i + \sum_{jk} a_{ijk} w_j w_k -w_i
$$
approach a steady state as $\tau$ goes to infinity.
\end{cor}
\begin{proof}
By Theorem \ref{rescperboltzmann} the vector field
$$
\mathbf v\doteq -t\frac \partial {\partial t} - x \cdot \nabla_x + u\cdot \frac \partial {\partial u}
$$
is the generator of symmetry associated to $G$, hence fixed $\theta$, the graph of the solution is mapped by $\mathbf v$ in the graph of another solution:
$$
\Gamma^\theta=\{\exp(\theta \mathbf v)(t,x,u): u= u(t,x)\}
$$ 
that is
$$
\left\{
\begin{array}{l}
t(\theta)=e^{-\theta}t
\\
x(\theta)=e^{-\theta}x
\\
u(t(\theta),x(\theta))=e^\theta u(t,x)
\end{array}
\right.
$$
thus, the new solution is $u^\theta(t,x)=e^\theta u(e^\theta t, e^\theta x)$.
\medskip

As far as the second part of the corollary is concerned, we shall use Theorem \ref{blowuprescaling}. Suppose that there exists a solution $w=(w_1(\tau,\eta),\dots,w_l(\tau,\eta))$ to the system 
$$
\partial_\tau w = F[w]-G[w]
$$
which approach a steady state as $\tau\to \infty$. Then the transformation by $\exp (\tau \mathbf v)$ with $\tau = \ln \frac1{1-t}$ is the graph of the solution 
$$
u(t,x)=\exp(\tau G) w(\tau,\eta)
$$
which corresponds to the change of variables
$$
\left\{
\begin{array}l
\dis 
x = e^{-\tau}\eta = (1-t)\eta 
\\
\dis
u_i = e^\tau w_i = \frac{w_i}{1-t} 
\end{array}
\right.
$$
which yields (\ref{systemchange}).
\end{proof}

\chapter{The two dimensional Broadwell model\label{chapblowup}}
\section{The discrete Boltzmann equation}
Consider the simplified model of a gas whose particles 
can have only finitely many speeds, say 
$\mathbf c_1,\ldots,\mathbf c_N\in\R^n$.  Call $u_i=u_i(t,x)$ the density
of particles with speed $\mathbf c_i$.
The evolution of these densities can then be described by a semilinear
system of the form
\begin{equation}
\partial_t u_i +\mathbf c_i\cdot \nabla u_i =\sum_{j,k} a_{ijk}\, u_ju_k
\qquad\qquad i=1,\ldots,N.
\label{dbzm-eq}
\end{equation}
Here the coefficient $a_{ijk}$ measures the rate at which new
$i$-particles are created, as a result of collisions between
$j$- and $k$-particles. 
In a realistic model, these coefficients must satisfy a set of identities,
accounting for the conservation of mass, momentum and energy.

Given a continuous, bounded  initial data
\begin{equation}
u_i(0,x)=\bar u_i(x),
\label{dbzm-12}
\end{equation}
on a small time interval $t\in [0,T]$
a solution of the Cauchy problem can be constructed by the
method of characteristics.
Indeed, since the system is semilinear,
this solution is obtained as the fixed point of the
integral transformation
$$
u_i(t,x)=\bar u_i(x-c_it)+\int_0^t
\sum_{j,k} a_{ijk}\, u_ju_k \big(s,~x-\mathbf c_i(t-s)\big)\,ds\,.
$$
For sufficiently small time intervals, the existence of a unique
fixed point follows from the contraction mapping principle, 
without any assumption on the constants $a_{ijk}$.

If the initial data is suitably small, the solution remains uniformly 
bounded for all times \cite{B2}. For large initial data, on the
other hand, the global existence and stability of solutions is 
known only in the one-dimensional case \cite{B1, HT, T}. 
Since the right hand side has quadratic growth,
it might happen that the solution blows up in finite time.
Examples where the $L^\infty$ norm of the solution
becomes arbitrarily large as $t\to\infty$ are easy to construct
\cite{I}.  In the present chapter we focus on the
two-dimensional Broadwell model and examine the possibility that
blow-up actually occurs in finite time. 

Since the equations (\ref{dbzm-eq}) admit a natural symmetry group,
one can perform an asymptotic rescaling of variables
and ask whether there is a blow-up solution which, in the rescaled variables,
converges to a steady state.   This technique has been widely used to
study blow-up singularities of reaction-diffusion equations with
superlinear forcing terms \cite{GV, GK}.  See also \cite{J} for an example of
self-similar blow-up for hyperbolic conservation laws.
Our results show, however, that for the two-dimensional Broadwell model
no such self-similar blow-up solution exists. 

If blow-up occurs at a time $T$, our results imply that for times $t
\to T-$ one has
$$
\big\|u(t)\big\|_{L^\infty}~>~ \frac 15\, \frac{\ln \big|\ln(T-t)\big|}{T-t}\,.
$$
This means that the blow-up rate must be different from
the natural growth rate $\big\|u(t)
\big\|_{L^\infty} = \O(1)\cdot(T-t)^{-1}$ 
which would be obtained
in case of a quadratic equation $\dot u=C\,u^2$.

In the final section of this chapter we discuss a possible scenario for
blow-up.  The analysis highlights how carefully chosen should
be the initial data, if blow-up is ever to happen.
This suggests that finite time blow-up is a highly non-generic
phenomenon, something one would not expect to encounter 
in  numerical simulations.
\section{Coordinate rescaling}
In the following, we say that $P^*=(t^*,x^*)$ is a \emph{blow-up point}
if 
$$\limsup_{~x\to x^*,\, t\to t^*-} u_i(t,x)=\infty$$
for some $i\in\{1,\ldots,N\}$.  
Define the constant
$$C\doteq \max_i |\mathbf c_i|\,.$$
We say that
$(t^*,x^*)$ is a \emph{primary blow-up point}
if it is a blow-up point and the backward cone
$$\Gamma\doteq \big\{ (t,x)\,;~~|x-x^*|< 2C\,(t^*-t)\big\}$$
does not contain any other blow-up point.
\begin{lemma}
\label{5-lemma1}
Let $u=u(t,x)$ be a solution
of the Cauchy problem (\ref{dbzm-eq})-(\ref{dbzm-12}) with continuous initial data.
If no primary blow-up point exist, then $u$ is continuous
on the whole domain $[0,\infty[\,\times\R^n$.
\end{lemma}
\begin{proof} If $u$ is not continuous, it must be unbounded in
the neighborhood of some point. Hence some blow-up point 
exists.
Call ${\cal B}$ the set of such blow-up points. Define the function
$$\vp(x)\doteq \inf_{(\tau,\xi)\in{\cal B}}\big\{\tau+C\,|x-\xi|\big\}\,.
$$
By Ekeland's variational principle (see \cite{AE}, p.254), 
there exists a point $x^*$ such that
$$\vp(x)\geq \vp(x^*)-\frac C2\, |x-x^*|$$
for all $x\in\R^2$.  Then $P^*\doteq \big(\vp(x^*),\,x^*\big)$
is a primary blow-up point. 
\end{proof}
Let now $(t^*,x^*)$ be a primary blow-up point.
One way to study the local asymptotic behavior of
$u$ is to rewrite the system in
terms of the rescaled variables $w_i= w_i(\tau,\eta)$, defined by
\begin{equation}
\left\{
\begin{array}{rcl}
\tau &=& \dis -\ln (t^*-t),\\
\eta &=& \dis e^{\tau}x~=~\frac{x-x^*}{t^*-t}\,,
\\
w_i &=& \dis e^{-\tau}u_i~=~(t^*-t)u_i.
\end{array}
\right.
\label{dbzm-rescaled}
\end{equation}
The corresponding system of evolution equations is
\begin{equation}
\partial_\tau w_i+(\mathbf c_i+\eta)\cdot \nabla_\eta w_i
=~-w_i+\sum_{j,k} a_{ijk} \, w_jw_k
\qquad\qquad i=1,\ldots,n.
\label{dbzm-req}
\end{equation}
Any nontrivial stationary or periodic solution $w$ of (\ref{dbzm-req})
would yield a solution $u$ of (\ref{dbzm-eq}) which blows up at $(t^*,x^*)$.
On the other hand, the non-existence of such solutions
for (\ref{dbzm-req}) would suggest that finite time blow-up
for (\ref{dbzm-eq}) is unlikely.
\section{The two-dimensional Broadwell model}
\begin{figure}[ht]
\centerline{
\includegraphics[width=10cm]{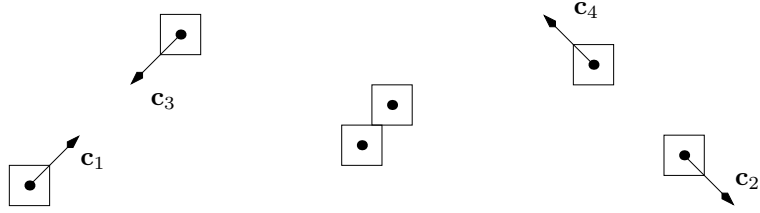}}
\caption{Moving particles with prescribed speeds\label{quadratini}}
\end{figure}

Consider a system on $\R^2$
consisting of 4 types particles (fig. \ref{quadratini}), with speeds
$$\mathbf c_1=(1,1),\qquad \mathbf c_2=(1,-1),\qquad \mathbf c_3=(-1,-1),\qquad \mathbf c_4=(-1,1).$$
\psfrag{c1}{$\mathbf c_1$}
\psfrag{c2}{$\mathbf c_2$}
\psfrag{c3}{$\mathbf c_3$}
\psfrag{c4}{$\mathbf c_4$}
\vskip 10pt
\vskip 10pt

The evolution equations are
\begin{equation}
\left\{
\begin{array}{l}
\partial_t u_1+c_1\cdot\nabla u_1=u_2u_4-u_1u_3\,,\\
\partial_t u_3+c_3\cdot\nabla u_3=u_2u_4-u_1u_3\,,\\
\partial_t u_2+c_2\cdot\nabla u_2=u_1u_3-u_2u_4\,,\\
\partial_t u_4+c_4\cdot\nabla u_4=u_1u_3-u_2u_4\,.
\end{array}
\right.
\label{dbzm-brdw}
\end{equation}
After renaming variables,
the corresponding rescaled system (\ref{dbzm-req}) takes the form
\begin{equation}
\left\{\begin{array}{l}
\partial_t w_1+(x+1)\partial_x w_1+(y+1)\partial_y w_1
=w_2w_4-w_1w_3-w_1\,,\\
\partial_t w_3+(x-1)\partial_x w_3+(y-1)\partial_y w_3
=w_2w_4-w_1w_3-w_3\,,\\
\partial_t w_2+(x+1)\partial_x w_2+(y-1)\partial_y w_2
=w_1w_3-w_2w_4-w_2\,,\\
\partial_t w_4+(x-1)\partial_x w_4+(y+1)\partial_y w_4
=w_1w_3-w_2w_4-w_4\,.\\
\end{array}
\right.
\label{dbzm-rbrdw}
\end{equation}
\v
Our first result rules out the possibility of asymptotically
self-similar blow-up solutions.  A sharper estimate will be proved later.
\begin{theorem}
\label{5-theo1}
The system (\ref{dbzm-rbrdw}) admits no 
nontrivial positive bounded solution which is constant or periodic in time.
\end{theorem}
\v
\begin{proof}  Assume 
\begin{equation}
0\leq w_i(t,x,y)\leq \kappa
\label{dbzm-pos}
\end{equation}
for all  $t,x,y$, $i=1,2$.
Choose $\ve\doteq e^{-2\kappa}/2$, so that
$$\ve<\frac 1\kappa\,,\qquad\qquad \ve e^{2\kappa x}\leq \frac 12
\quad x\in [-1,1]\,.$$
Define
$$Q_{14}(t,y)\doteq
\int_{-1}^1 \Big[\big(1-\ve e^{2\kappa x}\big) w_1(t,x,y)
+\big(1-\ve e^{-2\kappa x}\big) w_4(t,x,y)\Big]\,dx\,.$$
$$Q_{14}(t)\doteq \sup_{|y|\leq 1} Q_{14}(t,y)\,,$$

\vskip 10pt
\psfrag{(1,1)}{$(1,1)$}
\psfrag{(1, -1)}{$(1,-1)$}
\psfrag{(-1, 1)}{$(-1,1)$}
\psfrag{(-1, -1)}{$(-1,-1)$}
\psfrag{y(t)}{$y(t)$}
\psfrag{y}{$y$}
\begin{figure}[ht]
\centerline{\includegraphics[width=8cm]{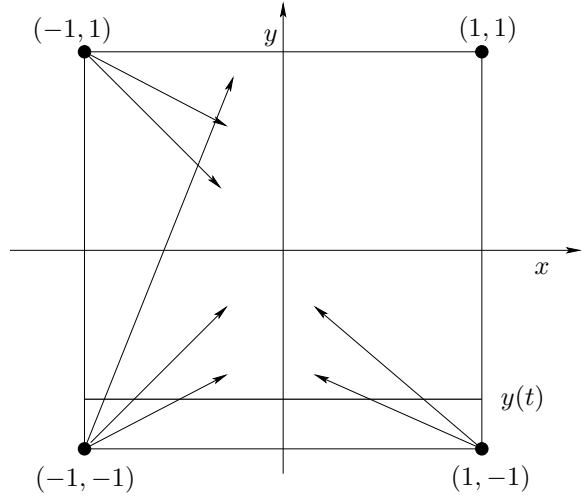}}
\caption{Interaction on a moving line\label{segmentoviaggiante}}
\end{figure}
\vskip 10pt

Restricted to any horizontal moving line 
$y=y(t)$ such that $\dot y=y+1$ (fig.\ref{segmentoviaggiante}), the equations (\ref{dbzm-rbrdw}) become
$$\begin{array}{l}
\partial_t w_1+(x+1)\partial_x w_1
=w_2w_4-w_1w_3-w_1\,,\\
\partial_t w_4+(x-1)\partial_x w_4
=w_1w_3-w_2w_4-w_4\,.\\
\end{array}
$$
A direct computation now yields
$$
\begin{array}{l}
\dis \frac d{dt}Q_{14}\big(t,y(t)\big)\\
\dis \qquad\leq
-2\ve\kappa \int_{-1}^1 \Big[e^{2\kappa x}(1+x)w_1
+e^{-2\kappa x}(1-x) w_4\Big] \,dx
\\
\qquad\qquad
\dis+\int_{-1}^1 \big(\ve
e^{2\kappa x}-\ve
e^{-2\kappa x}\big)
(w_1w_3-w_2w_4)\,dx \\
\dis \qquad \leq -\ve\kappa  \int_{-1}^1 \Big[e^{2\kappa x}(1+x)w_1
+e^{-2\kappa x}(1-x) w_4\Big] \,dx
-\ve\kappa \int_{-1}^0 e^{-2\kappa x}(1-x) w_4\,dx
\\
\dis \quad\qquad 
-\ve\kappa \int_0^1e^{2\kappa x}(1+x)w_1\,dx
+\int_{-1}^0 \ve \kappa e^{-2\kappa x} w_4\,dx+
\int_0^1 \ve \kappa e^{2\kappa x} w_1\,dx\\
\dis \qquad \leq -\ve\kappa  \int_{-1}^1 \Big[e^{2\kappa x}(1+x)w_1
+e^{-2\kappa x}(1-x) w_4\Big] \,dx\,.
\end{array}
$$
Call 
$$I(t,y)\doteq\int_{-1}^1 w_1(t,x,y)\,dx\,.$$
The definition of $\ve$ and the bound (\ref{dbzm-pos}) on $w_1$ imply 
$$
\begin{array}{rl}
\dis\int_{-1}^1 e^{2\kappa x}(1+x)w_1\,dx
&
\dis \geq 2\ve \int_{-1}^1 (1+x)w_1\,dx
\\
&\dis
\geq  2\ve\int_{-1}^{-1+I/\kappa} (1+x)\kappa\,dx\\
&\dis=\ve\,I^2/\kappa
\end{array}
$$
From this, and a similar estimate for $w_4$, we obtain
$$
\begin{array}{rl}
\dis
\int_{-1}^1 \Big[e^{2\kappa x}(1+x)w_1
+e^{-2\kappa x}(1-x) w_4\Big] \,dx
&\geq\dis
\frac \ve\kappa \left(\int_{-1}^1 w_1\,dx\right)^2+
\frac\ve\kappa \left(\int_{-1}^1 w_4\,dx\right)^2\,
\\
&\dis \geq \frac\ve{\kappa}
\frac {Q^2_{14}}2\,.
\end{array}
$$
Since $\ve<\kappa^{-1}$, this yields
\begin{equation}
\frac d{dt}Q_{14}\big(t,y(t)\big)
\leq -\frac{\ve^2}2\,Q^2_{14}\big(t,y(t)\big)
\,.
\label{dbzm-decrQ}
\end{equation}
Observing that the Cauchy problem
$$\dot z=-\frac{\ve^2}2\,z^2,\qquad\qquad z(0)=4\kappa$$
has the solution  
$$z(t)=\left( \frac 1{4\kappa}+\frac{\ve^2}2 t\right)^{-1},$$
by a comparison argument from (\ref{dbzm-decrQ}) we deduce
$$Q_{14}(t)\leq\left( \frac 1{4\kappa}+\frac{\ve^2}2 t\right)^{-1}.$$
Since
$$\int_{-1}^1\!\int_{-1}^1 w_1(t,x,y)\,dxdy\leq 4Q_{14}(t),$$
and since a similar estimate can be performed for all components
$w_i$, we conclude
\begin{equation}
\int_{-1}^1\!\int_{-1}^1 w_i(t,x,y)\,dxdy\leq 4
\left( \frac 1{4\kappa}+\frac{e^{-4\kappa}}8 t\right)^{-1}.
\label{dbzm-bound}
\end{equation}
The right hand side of (\ref{dbzm-bound}) approaches zero as $t\to\infty$.
Therefore, nontrivial constant or time-periodic $L^\infty$
solutions of (\ref{dbzm-rbrdw}) cannot exist.
\end{proof}

\section{Refined blow-up estimates}
If $(t^*,x^*)$ is a blow-up point, our analysis has 
shown that
in the rescaled coordinates $\tau,\xi$ the corresponding functions
$w_i$ must become unbounded as $\tau\to\infty$.   In this section
we refine the previous result, establishing a lower bound
for the rate at which such explosion takes place.
\begin{theorem}
\label{5-theo2}
Let $u$ be a continuous
solution of the Broadwell system (\ref{dbzm-rescaled}). Fix any point
$(t^*,x^*)$ and consider the corresponding rescaled variables
$\tau,\xi,w_i$. If
$$
\max_{|\xi_1|, |\xi_2|\leq 1}\,
w_i(\tau,\xi_1,\xi_2)\leq \theta\ln\tau\qquad\qquad i=1,2,3,4\,,
$$
for some $\theta<1/4$ and all $\tau$ sufficiently large,
then
$$
\lim_{\tau\to\infty} w_i(\tau,\xi)=0\qquad\qquad i=1,2,3,4\,,
$$
uniformly for $\xi \in\R^2$ in compact sets.
Therefore $(t^*,x^*)$ is not a blow up point.
\end{theorem}
\v
Since $w_i=(t^*-t) u_i$ and 
$\tau\doteq \big|\ln(t^*-t)\big|$, the above implies
\v
\begin{cor} If $(t^*,x^*)$ is a primary blow-up point, 
then
$$\limsup_{x\to x^*,~~t\to t^*- }~~  \big|u(t,x)\big|\cdot 
\frac{t^*-t}  {\ln\big|\ln(t^*-t)\big|}~\geq ~\frac 14\,.$$
\end{cor}
\v
\n \emph{Proof of Theorem \ref{5-theo2}.} 
\v
Let $w_i=w_i(t,x,y)$ provide a solution to the system (\ref{dbzm-rbrdw}),
with 
\begin{equation}
0\leq w_i(t,x,y)\leq \theta \,\ln t\doteq k(t)
\label{dbzm-kt}
\end{equation}
for all $t\geq t_0$ and $x,y\in [-1,1]$. The proof will
be given in two steps. First we show that the $L^1$ norm of the
components $w_i$ approaches zero as $t\to\infty$.
Then we refine the estimates, and prove that also the $L^\infty$
norm asymptotically vanishes.
\v
\n STEP 1: Integral estimates.
Consider the function
$$
Q_{14}(t,y)\doteq\int_{-1}^1 \left[\left(1- \frac{e^{2k(t)(x-1)}}2\right)
w_1(t,x,y)+\left(1- \frac{e^{-2k(t)(x+1)}}2 \right)w_4(t,x,y)\right]dx
$$
with $k(t)$ as in (\ref{dbzm-kt}).
As in the proof of Theorem \ref{5-theo1}, let $t\mapsto y(t)$ be a solution to
$\dot y=y+1$. Then
$$
\begin{array}{rl} 
\dis \frac d{{dt}}Q_{14}\big(t,\,y(t)\big)
=&\!\!\!\dis \int_{-1}^1\big[-(x-1) k'  e^{2k(t)(x-1)}w_1+(x+1) k' e^{-2k(t)(x+1)}w_4
\big]\,dx
\\ 
&\!\!\!\dis +\int_{-1}^1\left(1- \frac{{e^{2k(t)(x-1)}}}2\right)
\big[-(x+1)w_{1x}+w_2w_4-w_1w_3-w_1\big]\,dx\\ 
&\!\!\!\dis +\int_{-1}^1\left(1-\frac{{e^{-2k(t)(x+1)}}}2\right)\big[-(x-1)w_{4x}+
w_1w_3-w_2w_4-w_4\big]\,dx
\end{array}
$$
To estimate the right hand side, we notice that
$$
\begin{array}{rl}
A&\dis\!\!\! \doteq\int_{-1}^1\left(1-\frac{{e^{2k(t)(x-1)}}}2\right)
\big[(1+x)w_{1x}+w_1\big]dx\geq k(t)\int_{-1}^1(x+1)e^{2k(x-1)}w_1\,dx
\\
B&\dis\!\!\! \doteq\int_{-1}^1\left(1-\frac {{ e^{-2k(t)(x+1)} }}2\right)
\big[(x-1)w_{4x}+w_4\big]dx\geq k(t)\int_{-1}^1(1-x)e^{-2k(x+1)}w_4\,dx
\\
C&\dis\!\!\! \doteq\int_{-1}^1(w_1w_3-w_2w_4)\left(\frac{{e^{2k(x-1)}}}2-
\frac{{e^{-2k(x+1)}}}2\right)\,dx
\\
& \dis \!\!\! \leq k(t)\int_0^1\frac{{e^{2k(x-1)}}}2 w_1dx+k(t)
\int_{-1}^0
\frac{{e^{-2k(x+1)}}}2 w_4\,dx\,.
\end{array}
$$
Therefore,
$$
\begin{array}{rl}
\dis \frac d{dt} Q_{14}\big(t,\,y(t)\big)
=&\dis \int_{-1}^1\big[-(x-1) k'e^{2k(t)(x-1)}w_1+(x+1) k' e^{-2k(t)(x+1)}
w_4\big]\,dx
\\
&\qquad\qquad-A-B+C
\\
\leq
&~\dis k'(t)\int_{-1}^1\big[(1-x)e^{2k(t)(x-1)}w_1+(1+x)
e^{-2k(t)(x+1)}
w_4\big]\,dx
\\
&\dis \quad -\frac {k(t)}2 \int_{-1}^1 \big[(1+x) e^{2k(t)(x-1)}w_1+(1-x)
e^{-2k(t)(x+1)}
w_4\big]\,dx.
\end{array}
$$
If $k(t)\geq 1/2$,
we claim that the following two inequalities hold: 
\begin{equation}
\begin{array}{l}
\dis(1-x)e^{2k(t)(x-1)}\leq 1-e^{2k(t)(x-1)}\,,
\\
\\
\dis(1+x)e^{-2k(t)(x+1)}\leq 1-e^{-2k(t)(x+1)}\,.
\end{array}
\label{dbzm-banal}
\end{equation}
To prove the first inequality we need to show that
$$h_k(s)\doteq 
1-e^{2ks}+se^{2ks}\geq 0\qquad\qquad \hbox{for all}~~s\in [-2,0]\,.$$
This is clear because
$
h_k(0)=0
$
and
$$
h'_k(s)=e^{2ks}(1-2k+2ks)\leq 0\qquad\qquad s\in [-2,0]
$$
if $k\geq {1/ 2}$. Hence 
$h_k(s)$ is positive for $s\in [-2,0]$, as claimed.
The second inequality in (\ref{dbzm-banal}) is proved similarly.

When $t\geq t_0\doteq e^{1/(2\theta)}$ one has
$k(t)\geq \frac 12$ and hence
$$\frac d{dt}Q_{14}\big(t,\,y(t)\big)
\leq
k'(t)Q_{14}- \frac{k(t)}2\int_{-1}^1 [(1+x) e^{2k(t)(x-1)}w_1+(1-x)
e^{-2k(t)(x+1)}w_4]dx\,.
$$
Setting $I=\int_{-1}^1w_1dx$, we 
obtain
$$
\int_{-1}^1(1+x)e^{2k(t)(x-1)}w_1dx 
~\geq~\int_{-1}^{-1+I/k(t)}(1+x)e^{-4k(t)} k(t)\,dx 
~=~e^{-4k(t)}\frac{I^2}{2k(t)}\,.
$$
Using the above, and a similar estimate for the integral of $w_4$, we obtain
\begin{equation}
\begin{array}{l}
\dis \frac{k(t)}2 \int_{-1}^1 [(1+x) e^{2k(t)(x-1)}w_1+(1-x)e^{-2k(t)(x+1)}w_4]
dx
\\
\dis \qquad \geq\frac { e^{-4k(t)}}4\left[\left( \int_{-1}^1w_1dx\right)^2+
\left( \int_{-1}^1w_4dx\right)^2\right]
\\
\qquad \dis \geq \frac{ e^{-4k(t)}} 8\,Q^2_{14}\,.
\end{array}
\label{bzmn-weighted}
\end{equation}
Calling 
$$Q_{14}(t)=\max_{|y|\leq 1} \,Q_{14}(t,y)\,,$$
from (\ref{bzmn-weighted}) we deduce
$$
\frac d{dt}Q_{14}(t)
\leq
k'(t)Q_{14}(t)-\frac{ e^{-4k(t)}}8\,Q_{14}(t)^2.
$$
Recalling that $k(t)=\theta\ln t$ for some $0<\theta<1/4$, the previous 
differential inequality can be written as
\begin{equation}
\frac d{dt} Q_{14}
\leq  \frac \theta t Q_{14}-\frac 1{8t^{4\theta}}Q_{14}^2.
\label{dbzm-diffineq}
\end{equation}

Notice that $ Q_{14}\big(t_0,y(t_0)\big)\leq 2k(t_0)$, and define
the constant
$$
A_0\doteq  \max\big\{2k(t_0)t_0^{1-4\theta}, ~8(1-3\theta)\big\}\,.$$
Then the function
$$
z(t)\doteq A_0\,t^{4\theta-1}
$$
satisfies
\begin{equation}
\frac d{dt} z(t)\geq \frac \theta t z-\frac 1{t^{4\theta}}z^2\qquad 
\qquad z(t_0)\geq Q_{14}
(t_0,y(t_0))\,.
\label{dbzm-zeq}
\end{equation}
Comparing  (\ref{dbzm-diffineq}) with (\ref{dbzm-zeq}) we conclude
$$
Q_{14}(t)\leq z(t)\qquad\qquad t\geq t_0\,.
$$
This implies the estimate  
$$
\int_{-1}^1 w_i(t,x,y_0)\,dx \leq 2 Q_{14}(t)\leq 2 A_0\,t^{4\theta-1}
$$
for $t\geq t_0$, $i\in \{1,4\}$ and any $y_0\in [-1,1]$. 
An entirely similar
argument applied to $Q_{12}$, $Q_{23},\ldots$ yields the estimates
\begin{equation}
\int_{-1}^1 w_i(t,x,y_0)dx \leq 2 A_0t^{4\theta-1}\,,\qquad\qquad
\int_{-1}^1 w_i(t,x_0,y)dy \leq 2 A_0t^{4\theta-1}\,.
\label{dbzm-wito0}
\end{equation}
for $i=1,2,3,4$, $x_0,y_0\in [-1,1]$ and $t\geq t_0$.
\vs

\n STEP 2: Pointwise estimates. 
Using the integral bounds (\ref{dbzm-wito0}), we now seek a uniform bound of the
form
\begin{equation}
w_i(t,x,y)\leq C_0
\label{dbzm-infty}
\end{equation}
for some constant $C_0$ and all $x,y\in [-1,\,1]$, $t>0$.  
\v
To prove (\ref{dbzm-infty}), let $t\mapsto x(t)$, $t\mapsto y(t)\big)\in [-1,1]$ 
be solutions of
$$\dot x=x+1\,,\qquad\qquad \dot y=y+1\,.$$
Call 
$$A(t)\doteq \int_{x(t)}^1 (w_1+w_4)\big(t,x,y(t)\big)\,dx\,.$$
From our previous  estimates (\ref{dbzm-wito0}) it trivially follows 
$$
A(t)\leq 4 A_0t^{4\theta-1}.
$$
The time derivative of
$A(t)$ is computed as
$$
\begin{array}{rl}
\dis \frac{dA}{dt}=&\dis -(x(t)+1)(w_1+w_4)(t,x(t),y(t))
\\
&+\dis
\int_{x(t)}^1
\big[\partial_t w_1+(y(t)+1)\partial_y w_1+\partial_tw_4+(y(t)+1)
\partial_y 
w_4\big]\,dx
\\
=&\dis -\big(x(t)+1\big)(w_1+w_4)\big(t,x(t),y(t)\big)
\\
&\dis
-\int_{x(t)}^1\big[w_1+(x+1)
\partial_x w_1+w_4+(x-1)\partial_xw_4\big]\,dx
\\
=&\dis -(x(t)+1)(w_1+w_4)\big(t,x(t),y(t)\big)-\int_{x(t)}^1(w_1+w_4)dx
\\
&\quad\dis -2w_1\big(t,1,y(t)\big)+(x(t)+1)w_1\big(t,x(t),y(t)\big)+
\int_{x(t)}^1w_1\,dx
\\
&\quad \dis -2w_4\big(t,1,y(t)\big)+(x(t)-1)w_4\big(t,x(t),y(t)\big)+
\int_{x(t)}^1w_4\,dx
\\
\leq&\dis  \big[x(t)-1-(x(t)+1)\big]w_4\big(t,x(t),y(t)\big)~=~-2w_4\big(
t,x(t),y(t)\big)
\,.
\end{array}
$$
This implies 
\begin{equation}
w_4\big(t,x(t),y(t)\big)\leq -\frac 12 \frac{dA}{dt}\,.
\label{dbzm-wminAprimo}
\end{equation}
The total derivative of $w_1$ along a characteristic line is now
given by
$$
\begin{array}{rl}
\dis \frac d{dt}w_1\big(t,x(t),y(t)\big)
=&
\dis 
w_2w_4-w_1w_3-w_1\leq w_2w_4-w_1\leq \frac 12 w_2\left(\frac {-dA}{dt}
\right)-w_1
\\
\leq&
\dis -w_1+\frac{k(t)} 2\left(\frac{-dA}{dt}\right)\,.
\end{array}
$$
In turn,  for $t\geq t_0$ this yields the inequality
\begin{equation}
\begin{array}{rl}
\dis w_1\big(t,x(t),y(t)\big)
&\dis \leq e^{-(t-t_0)}
\left[w_1(t_0)+\int_{t_0}^te^{s-t_0}k(s)(-A'(s))ds\right]
\\
&\dis \leq e^{-(t-t_0)} \big[w_1(t_0)+A(t_0)k(t_0)\big]
\\
&\quad\dis +e^{-(t-t_0)}
\int_{t_0}^t A(s)\big(e^{s-t_0}k(s)\big)'\, ds.
\end{array}
\label{dbzm-Apesato}
\end{equation}
The first term on the right hand side of (\ref{dbzm-Apesato})
approaches zero exponentially fast.
Concerning the second, we have
$$e^{-(t-t_0)}
\int_{t_0}^t A(s)\,e^{s-t_0}\big(k(s)+k'(s)\big)\,ds
\leq \int_{t_0}^t e^{-(t-s)}\,2A_0 s^{4\theta-1}\,\left(\theta
\ln s+
\frac\theta s\right)\,ds\,.
$$
This also approaches zero as $t\to\infty$.
Repeating the same computations for all components, we conclude that
for some time $t_1$ sufficiently large there holds
\begin{equation}
\label{wlimitato}
w_i(t_1,x,y) < \frac 12\qquad\qquad \hbox{for all}~x,y\in [-1,1]\,.
\end{equation}
By continuity, the inequalities in (\ref{wlimitato}) remain valid for all
$x,y$ in a slightly larger square, say
$[-1-\epsilon,~1+\epsilon]$.
For $t\geq t_1$ we now define
$$M(t)\doteq \max\Big\{ w_i(t,x,y)\,;~~i=1,2,3,4,~~~x,y\in
[-1-e^{t-t_1}\epsilon\,,~1+e^{t-t_1}\epsilon]\Big\}\,.$$
From the equations (\ref{dbzm-rbrdw}) and (\ref{wlimitato}) it now follows
$$\frac d {dt} M(t)\leq -M(t)+M^2(t)\leq \frac{M(t)}2\,,\qquad\qquad M(t_1)
\leq \frac 12\,.$$
$$M(\tau)\leq \left[ 1+e^{\tau-t_1}\left(\frac 1{M_1}-1\right)\right]^{-1}
\leq e^{-\tau}\cdot e^{t_1}\qquad\qquad\hbox{for all}~~ \tau\geq t_1.$$
Returning to the original variables $u_i=e^\tau w_i$, this yields
$$u_i\leq e^{t_1}$$
in a whole neighborhood of the point $P^*=(t^*, x^*)$.  Hence 
$P^*$ is not a blow-up point.
\endproof
\section{A tentative blow-up scenario}
For a solution of the rescaled equation (\ref{dbzm-brdw}), the total mass
$$m(t)\doteq \int_{-1}^1\int_{-1}^1 \sum_{i=1}^4 w_i(t,x,y)\, dxd y$$
may well become unbounded as $t\to\infty$.
On the other hand, the one-dimensional 
integrals along horizontal or vertical segments 
decrease monotonically.  
Namely, if $t\mapsto y(t) $ satisfies $\dot y=y+1$, then
$$\frac d{dt}\int_{-1}^1 \big[w_1(t,x,y(t))+w_4(t,x,y(t))\big]\,dx \leq 
0\,.
$$
Similarly, if $\dot x=x-1$, then 
$$\frac d{dt}\int_{-1}^1 \big[w_3(t,x(t),y)+w_4(t,x(t),y)\big]\,dy 
\leq 0\,.$$
Analogous estimates hold for
the sums $w_1+w_2$ and $w_2+w_3$.
Therefore, a bound on the initial data
$$w_i(0,x,y)\in [0, M]\qquad\qquad\hbox{for all}~~x,y\in [-1,1]\,,$$
yields uniform integral bounds on the line integrals of all
components:
$$
\int_{-1}^1 w_i(t,x,y)\,dx\leq 4M\,,\qquad\qquad
\int_{-1}^1 w_i(t,x,y)\,dy\leq 4M\,.
$$
\psfrag{(1,1)}{$(1,1)$}
\psfrag{(1,-1)}{$(1,-1)$}
\psfrag{(-1,1)}{$(-1,1)$}
\psfrag{(-1,-1)}{$(-1,-1)$}
\psfrag{P1}{$P1$}
\psfrag{P2}{$P2$}
\psfrag{P3}{$P3$}
\psfrag{P5}{$P5$}
\psfrag{Q}{$Q$}
\psfrag{Q'}{$Q'$}
\vskip 10pt
\begin{figure}[ht]
\centerline{\includegraphics[width=8cm]{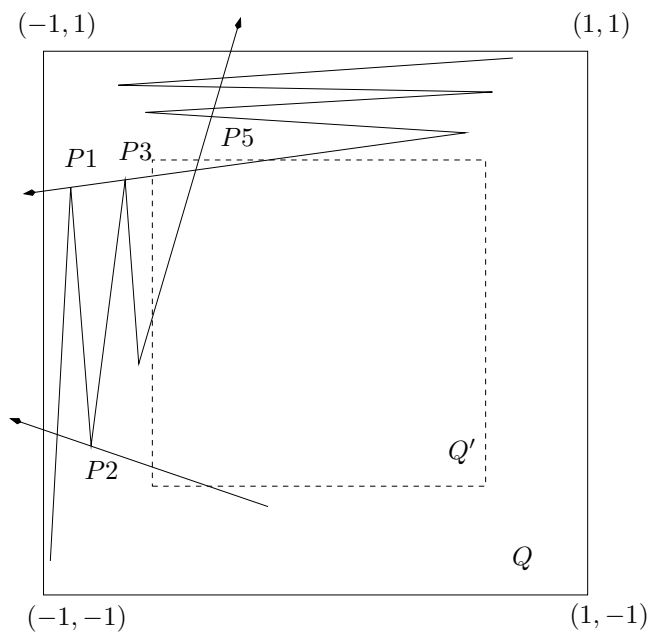}}
\caption{A possible interaction between particles\label{interpart}}
\end{figure}
\vskip 10pt

If finite time blow-up is to occur, the mass which is initially 
distributed along each horizontal or vertical segment
must concentrate itself within a very small region, thus forming
a narrow packet of particles with increasingly high density.
A possible scenario is illustrated in fig.~\ref{interpart}. 
A packet of 1-particles is initially located at $P_1$. 
In order to contribute to blow-up, this packet must remain
within the unit square $Q$. At $P_2$ these 1-particles interact
with 3-particles and produce a packet of 4-particles.
In turn, at $P_3$ these interact with 2-particles 
and produce again a packet of 1-particles.
After repeated interactions, the packet of alternatively 
1- and 4-particles eventually enters within the smaller square $Q'$.
After this time, it interacts with a packet of 2-particles at $P_5$
(transforming it into a packet of 1-particles)
and eventually exits from the domain $Q$.

To help intuition, it is convenient to describe a packet as
being ``young''
until it enters the smaller square $Q'$, and ``old'' afterwards.
To maintain a young packet inside $Q$, one needs the presence of
old packets interacting with it near the points
$P_2,P_3,P_4\ldots$ ~  On the other hand,
after it enters $Q'$, our packet can in turn be used to
hit another young packet, say at $P_5$, and preventing it from
leaving the domain $Q$.

As $t\to\infty$, the density of the packets must approach infinity.
One thus expects that most of the mass will be concentrated
along a finite number of one-dimensional curves.
Say, the packet of alternatively 1- and 4-particles should be
located along a moving curve $\gamma_{14}(t,\theta)$, where $\theta$
is a parameter along the curve.  The time evolution of such a curve
is of course governed by the equations

$$\frac \partial{\partial t}\gamma_{14}=c_1\qquad\hbox{or}\qquad
\frac \partial{\partial t}\gamma_{14}=c_4$$
depending on whether $\gamma_{14}(t,\theta)$ consists of 1- or 4-particles.
The presence of interactions impose 
highly nonlinear constraints on
these curves.   For example,
the interaction occurring in $P_5$
at time $t$ implies
the crossing of the two curves $\gamma_{14}$ and $\gamma_{12}$,
namely
$$\gamma_{14}(t,\theta)=\gamma_{12}(t,\tilde \theta)=P_5$$
for some parameter values $\theta,\tilde\theta$.
The complicated geometry of these curves resulting from the
above constraints
has not been analyzed. 

\addcontentsline{toc}{chapter}{Bibliography}

\end{document}